\documentclass[a4paper,12pt]{article}
\usepackage{amsthm}
\usepackage{amsmath}
\allowdisplaybreaks[4]%
\usepackage{bm}%
\usepackage{amssymb}
\usepackage{tabularx}
\usepackage{indentfirst}
\usepackage{enumerate}
\usepackage{dsfont}
\usepackage{graphicx}
\usepackage{subfigure}
\usepackage{epsfig}
\usepackage{graphics}
\usepackage{cases}
\usepackage[compress]{cite}
\usepackage{txfonts}
\usepackage{geometry}
 \usepackage{epstopdf}
 \usepackage{lineno}
 \usepackage[justification=centering]{caption}
 \usepackage{color}
 \usepackage{hyperref}
\hypersetup{
colorlinks=true,
linkcolor=blue,
anchorcolor=blue,
citecolor=blue}
\newtheorem{theorem}{Theorem}[section]
\newtheorem{lemma}{Lemma}[section]

\newtheorem{exm}{Example}[section]

\newtheorem{remark}{Remark}[section]
\newtheorem{definition}{Definition}[section]

\numberwithin{equation}{section}

\begin{document}

\thispagestyle{empty}

\title{On an impulsive faecal-oral model in a moving\\ infected environment \thanks{The first author is supported by Postgraduate Research \& Practice Innovation Program of Jiangsu Province (KYCX24\_3711) and the second author acknowledges the support of the National Natural Science Foundation of China (No. 12271470).}}

\date{\empty}

\author{Qi Zhou$^1$, Zhigui Lin$^1 \thanks{Corresponding author. Email: zglin@yzu.edu.cn (Z. Lin).}$, and Michael Pedersen$^2$
\\
{\small 1 School of Mathematical Science, Yangzhou University, Yangzhou 225002, China}\\
{\small 2 Department of Applied Mathematics and Computer Science,}\\
{\small Technical University of Denmark, DK 2800, Lyngby, Denmark}
}
 \maketitle
\begin{quote}
\noindent
{\bf Abstract.}{\footnotesize\small~This paper develops an impulsive faecal-oral model with free boundary to in order to understand how the exposure to a  periodic disinfection and expansion of the  infected region together influences the spread of faecal-oral diseases. We first check that this impulsive model has a unique globally nonnegative classical solution. The principal eigenvalues of the corresponding periodic eigenvalue problem at the initial position and infinity are defined as $\lambda^{\vartriangle}_{1}(h_{0})$ and $\lambda^{\vartriangle}_{1}(\infty)$, respectively. They both depend on the impulse intensity $1-G'(0)$ and expansion capacities $\mu_{1}$ and $\mu_{2}$. The possible long time dynamical behaviours of the model are next explored in terms of  $\lambda^{\vartriangle}_{1}(h_{0})$ and $\lambda^{\vartriangle}_{1}(\infty)$: if $\lambda^{\vartriangle}_{1}(\infty)\geq 0$, then the diseases are vanishing; if $\lambda^{\vartriangle}_{1}(\infty)<0$ and $\lambda^{\vartriangle}_{1}(h_{0})\leq 0$, then the disease are spreading; if $\lambda^{\vartriangle}_{1}(\infty)<0$ and $\lambda^{\vartriangle}_{1}(h_{0})> 0$, then for any given $\mu_{1}$, there exists a $\mu_{0}$ such that spreading happens as $\mu_{2}\in( \mu_{0},+\infty)$, and vanishing happens as $\mu_{2}\in(0, \mu_{0})$. Finally, numerical examples are presented  to corroborate the correctness of the obtained theoretical findings and to further understand the influence of an impulsive intervention and expansion capacity on the spreading of the diseases. Our results show that both the increase of impulse intensity and the decrease of expansion capacity have a positive contribution to the prevention and control of the diseases.
}


\end{quote}

\section{Introduction}
Infectious diseases have consistently been a major threath to  human health, and the repeated outbreaks of infectious diseases in history have brought great disasters to human livelihoods and national economies. As early as 2nd century Anno Domini, the plague pandemic in the Roman Empire caused sharp population decline and economic deterioration, leading to the collapse of the empire \cite{Ma-Zhou-Wang}. Since 1816, cholera has remained a major public health threat in areas with limited access to clean water and sanitation services \cite{Xu-Zou-Dent}. In the year of 2010 alone, it is estimated that cholera affects 3-5 million people and causes 100,000-130,000 deaths in the world \cite{Tian-Wang}. COVID-19 has been globally prevalent during the last five years, and Hay and his collaborators estimated that 15.9 million people died globally from the COVID-19 pandemic in 2020 and 2021 combined \cite{Hay-collaborators}.

Each kind of infectious disease usually has its own specific transmission route, such as contact transmission \cite{Shen-Lai-Zhang}, faecal-oral transmission \cite{Zhou-Lin-Santos}, droplet transmission \cite{Brankston-Gitterman-Hirji}, vector-borne transmission \cite{Rocklov-Dubrow}, blood transmission \cite{Hewitt-Ijaz-Brailsford}, vertical transmission \cite{Cao-Krogstad-Korber}, and so on. The faecal-oral route, as a common route of transmission, is where pathogenic bacteria are excreted in the faeces of patients or carriers, and then ingested by susceptible individuals through various ways \cite{Korn}. Enteric viruses like norovirus, rotavirus and astrovirus have long been accepted as spreading in the population through faecal-oral transmission \cite{Ghosh-Kumar-Santiana}.
Examples of common infectious diseases caused by faecal-oral transmission are cholera, hand-foot-mouth diseases, and so on \cite{Zhou-Lin-Santos}. Additionally, numerous references suggested that COVID-19 may also be spread by this route, see, for instance, \cite{Xu-Li-Zhu, Zhang-Castro-Mccune, Heller-Mota-Greco}. Thus, the faecal-oral transmission route has received extensive attention from a number of  mathematicians, epidemiologists, and medical scientists.

Mathematical models have become a highly effective and important tool for understanding the transmission mechanisms and designing the control measures of infectious diseases \cite{Zhou-Li-Hu-Zhang}. During the previous few decades, a number of faecal-oral models have been established to study the development of the infectious diseases transmitted by the faecal-oral route, see, for instance, \cite{Capasso-Paveri-Fontana,Capasso-Maddalena-1, Capasso-Maddalena-3}. In the summer of 1973, cholera  occurred in the Mediterranean region of Europe. To describe its spread, Capasso and Paveri-Fontana first proposed the ordinary differential equation (ODE) model \cite{Capasso-Paveri-Fontana}. Notice that this ODE model ignores the heterogeneity of space and the movement of infected individuals and pathogenic bacteria, and hence fails to show the spatial pattern of the diseases. In addition, the outbreak predictions can have severe deviations from reality. To obtain a  more  realistic model, Capasso and Maddalena modified the ODE model into a partial differential equation (PDE) model with Neumann boundary conditions \cite{Capasso-Maddalena-1} and Robin boundary conditions \cite{Capasso-Maddalena-3}, and investigated the long time dynamical behaviours of the models.

Notice that the models mentioned above are all formulated on a fixed spatial region. But when there is an outbreak of diseases transmitted by the faecal-oral route, it is very significant to know how fast the infected region is moving. To answer this question, Wang and Du \cite{Wang-Du, Du-Ni-Wa} developed the following free boundary model:
\begin{eqnarray}\label{Wang-Du}
\left\{
\begin{array}{ll}
u_{t}=d_{1}\Delta u-a_{11}u+a_{12}v,\; &\, t>0, x\in(g(t),h(t)), \\[2mm]
v_{t}=d_{2}\Delta v-a_{22}v+f(u),\; &\, t>0, x\in(g(t),h(t)), \\[2mm]
u=0, v=0,\; &\, t>0, x\in\{g(t), h(t)\},\\[2mm]
g'(t)=-\mu_{1}u_{x}(t,g(t))-\mu_{2}v_{x}(t,g(t)),\; &\, t>0,\\[2mm]
h'(t)=-\mu_{1}u_{x}(t,h(t))-\mu_{2}v_{x}(t,h(t)),\; &\, t>0,\\[2mm]
g(0)=-h_{0}, u(0,x)=u_{0}(x),\; &\,x\in[-h_{0},h_{0}],\\[2mm]
h(0)=h_{0},~~ v(0,x)=v_{0}(x),\; &\,x\in[-h_{0},h_{0}]
\end{array} \right.
\end{eqnarray}
with suitable initial functions $u_{0}(x)$ and $v_{0}(x)$, where $g(t)$ and $h(t)$ are the moving left and right boundaries of the infected region at time $t$. The meaning of all the notations in model \eqref{Wang-Du} is listed in \autoref{Tab:01}
\begin{table}[!htb]
\centering
\caption{List of notations and their meanings in model \eqref{Wang-Du}}
\vspace{0.1cm}
\label{Tab:01}
\small
\begin{tabular}{ll}
\hline
 Notation& Biological meaning\\
\hline
$u(t,x)$ & Spatial density of the pathogenic bacteria at time $t$ and position $x$\\
$v(t,x)$ & Spatial density of the infected individuals at time $t$ and position $x$\\
$d_{1}$  & Diffusion coefficient of the pathogenic bacteria\\
$d_{2}$  & Diffusion coefficient of the infected individuals\\
$a_{11}$  & Unit decreasing rate of the pathogenic bacteria\\
$a_{22}$  & Unit decreasing rate of the infected individuals\\
$f$           & Infection rate of human population caused by the pathogenic bacteria                \\
$a_{12}$  & Unit growth rate of the pathogenic bacteria caused by the infected individuals\\
$\mu_{1}$   &Expansion capacity of the pathogenic bacteria \\
$\mu_{2}$   &Expansion capacity of the infected individuals \\
$h_{0}$   &Length of the right boundary of the initial infection region\\
\hline
\end{tabular}
\end{table}

In model \eqref{Wang-Du}, the fourth and fifth equations stipulate that the expansion rate of the infected region is inversely proportional to a linear combination of the spatial gradients of the spatial density $u(t,x)$ of the pathogenic bacteria and the spatial density $v(t,x)$ of the infected individuals at the front. It is a particular example of the well-known Stefan condition, which has been used to model a number of problems. For instance, it can be used to describe the spreading of a new or invasive species \cite{Du-Lin}, the healing of wounds \cite{Chen-Friedman}, the growth of tumours \cite{Chen-Friedman-1}, the melting of ice \cite{Rubinstein}, the infiltration of vapours \cite{Merz-Rybka}, and the oxygen intake in muscles \cite{Crank}, to mention a few.

In addition to the moving infection environment, the diseases are usually affected by human beings in the spreading process. For example, the periodic spraying of disinfectant liquids leads to a sharp reduction in the spatial density of pathogenic bacteria within the environment for a short period of time. This prevention and control measure usually influences or even alters the dynamic behaviour of the diseases. However, classical PDE models do not describe this human intervention phenomenon. It is well known that the most prominent feature of impulse differential models is that they can adequately take into account the effects  of instantaneous bursts phenomena on the state of the models. Therefore, researchers have turned their attention to the impulse differential models.

The theoretical investigation of impulsive differential equations began with the 1960s work of Mil'man and Myshkis \cite{Mil'man-Myshkis}. In the field of biomathematics, Lewis and Li proposed  impulsive reaction diffusion equation models describing a seasonal birth pulse in the species, and showed that how the pulse affects the species' spread and persistence \cite{Lewis-Li}. Based on the work of \cite{Lewis-Li}, there have been a number of studies on impulse reaction diffusion equation models. We refer to \cite{Fazly-Lews-Wang} for  higher dimensional models extended from \cite{Lewis-Li}, \cite{Fazly-Lews-Wang-1} for a hybrid impulsive reaction-advection-diffusion model, \cite{Wu-Zhao} for an impulsive integro-differential model and a further research with spatial heterogeneity in \cite{Wu-Zhao-1}, \cite{Li-Zhao-Cheng} for a mosquito model with periodic evolution domain, \cite{Meng-Lin-Pedersen-1} for a competition model in stream environments, \cite{Zhang-Yi-Chen} for a model with shifting environments, and references therein. However, to the best of our knowledge, no research has considered both the implementation of an impulsive intervention and the expansion of the infected environment up to now.

This paper develops an impulsive faecal-oral model in a moving infected environment by incorporating the impulse intervention into model \eqref{Wang-Du}.
Because of the introduction of the impulse intervention, the theoretical analysis of the model is more challenging compared to \cite{Wang-Du}. Specifically, the solution no longer has the half-flow property, thus some energy methods are no longer available; impulsive conditions (to be described in detail later) should also be considered when constructing the upper and lower solutions; the steady state is governed by a more complex parabolic problem.

With the introduction of the impulse intervention, it is natural to ask whether some results in \cite{Wang-Du} like the well-posedness, the spreading-vanishing dichotomy, and so on still hold? How does the pulse intervention influence the moving front of the infected environment? Additionally, how the implementation of the impulse intervention and the expansion of the infected region together influence or even alter the spread of faecal-oral diseases. These are the motivations of this paper. The main contributions of this work are shown as follows:
\begin{itemize}
\item{An impulsive faecal-oral model in a moving infected environment is formulated by taking the impulsive intervention into account, and the well-posedness of the model is proved.}
\item{A spreading-vanishing dichotomy for this model is proved in terms of the principal eigenvalue $\lambda^{\vartriangle}_{1}(\infty)$. Based on this, sharp criteria for spreading and vanishing are given.}
\end{itemize}

The organization of the rest of this paper is as follows. In the next section, we formulate an impulsive faecal-oral model with free boundary, and then prove
that it has a unique globally nonnegative classical solution. \autoref{Section-3} recalls some main results for the impulsive faecal-oral model on a fixed region, and proves that the diseases are either spreading or vanishing. This section also establishes the criteria that determine spreading and vanishing, which are governed by the principal eigenvalues $\lambda^{\vartriangle}_{1}(h_{0})$ and $\lambda^{\vartriangle}_{1}(\infty)$, and the expansion capacities $\mu_{1}$ and $\mu_{2}$, or the initial functions $u_{0}(x)$ and $v_{0}(x)$. By numerical simulation, \autoref{Section-4} confirms the correctness of the theoretical findings, and  visually shows the impact of the impulsive intervention and the moving infection environment on the spread of the diseases.
Finally, a brief discussion and conclusion is given in \autoref{Section-5}.

\section{Model derivation and well-posedness}\label{Section-2}
In this section we first develop an impulsive faecal-oral model with free boundary, and then prove that this impulsive model has a unique globally nonnegative classical solution.

As mentioned in the introduction, the evolution of faecal-oral transmitted diseases is not a continuous process and cannot be described simply by classical differential equations. Instead, impulsive differential equations serve as a more appropriate way to describe it. By including an impulsive intervention into model \eqref{Wang-Du}, the impulsive faecal-oral model with free boundary is described by the following equations:
\begin{eqnarray}\label{Zhou-Lin}
\left\{
\begin{array}{ll}
u_{t}=d_{1}\Delta u-a_{11}u+a_{12}v,\; &\, t\in((k\tau)^{+}, (k+1)\tau], x\in(g(t),h(t)), \\[2mm]
v_{t}=d_{2}\Delta v-a_{22}v+f(u),\; &\, t\in((k\tau)^{+}, (k+1)\tau], x\in(g(t),h(t)), \\[2mm]
u=0, v=0,\; &\, t\in(k\tau, (k+1)\tau], ~~x\in\partial(g(t),h(t)),\\[2mm]
u((k\tau)^{+},x)=G(u(k\tau,x)), \; &\, x\in(g(k\tau), h(k\tau)), \\[2mm]
v((k\tau)^{+},x)=v(k\tau,x), \; &\,  x\in(g(k\tau), h(k\tau)), \\[2mm]
g'(t)=-\mu_{1}u_{x}(t,g(t))-\mu_{2}v_{x}(t,g(t)),\; &\, t\in(k\tau, (k+1)\tau],\\[2mm]
h'(t)=-\mu_{1}u_{x}(t,h(t))-\mu_{2}v_{x}(t,h(t)),\; &\, t\in(k\tau, (k+1)\tau],\\[2mm]
g(0)=-h_{0}, u(0,x)=u_{0}(x),\; &\,x\in[-h_{0},h_{0}],\\[2mm]
h(0)=h_{0},~~~ v(0,x)=v_{0}(x),\; &\,x\in[-h_{0},h_{0}], ~k=0,1,2,3,4,\cdots,
\end{array} \right.
\end{eqnarray}
where $\tau$ represents the time span between two neighboring pulse interventions, and $u((k\tau)^{+},x)$ and $v((k\tau)^{+},x)$ denote the right limit of $u(t,x)$ and $v(t,x)$ at $t=k\tau$, respectively. We always take $k=0,1,2,\cdots$ unless otherwise specified. Here, the other symbols have the same meaning as in model \eqref{Wang-Du}.

This model captures the situation in which a natural developmental stage and a human intervention stage occur alternately. During the natural spreading phase, the spatial densities of pathogenic bacteria and infected individuals are described by continuous equations, see the first two lines of model \eqref{Zhou-Lin}.
At the end of the human intervention phase, the spatial densities are represented by impulse equations, see the fourth and fifth lines of model \eqref{Zhou-Lin}.
As this paper investigates the long-time dynamical behaviors of model \eqref{Zhou-Lin}, $1-G'(0)$ is used to denote impulse intervention intensity.

In order to simplify the notation, we write, $\rho_{1}(u)$ and $\rho_{2}(u)$ if present, as $f(u)$ and $G(u)$ respectively. We make the following assumptions
about the initial  functions $u_{0}$ and $v_{0}$, the growth function $f(u)$ and the  impulsive function $G(u)$:
\begin{itemize}
\item[(A1)] Initial functions $u_{0}(x)$ and $v_{0}(x)$ satisfy
\begin{eqnarray*}
\left\{
\begin{array}{l}
u_{0}\in\mathbb{C}^{2}([-h_{0},h_{0}]), ~u_{0}(\pm h_{0})=0~\text{and}~u_{0}(x)>0~\text{in}~(-h_{0},h_{0}), \\[2mm]
v_{0}\in\mathbb{C}^{2}([-h_{0},h_{0}]), ~v_{0}(\pm h_{0})=0~\text{and}~v_{0}(x)>0~~\text{in}~(-h_{0},h_{0});
\end{array} \right.
\end{eqnarray*}
\item[(A2)] Growth function $f(u)$ satisfy
\begin{eqnarray*}
\left\{
\begin{array}{l}
f\in \mathbb{C}^{1}([0,\infty)), ~f(0)=0~\text{and}~f'(u)>0~\text{for}~u\in[0, \infty), \\[2mm]
\frac{f(u)}{u}~\text{is~strongly~decreasing}~\text{for}~u>0~\text{and}~\lim\limits_{u\rightarrow+\infty}\frac{f(u)}{u}<\frac{a_{11}a_{22}}{a_{12}};
\end{array} \right.
\end{eqnarray*}
\item[(A3)] Impulsive function $G(u)$ satisfy
\begin{eqnarray*}
\left\{
\begin{array}{l}
G\in \mathbb{C}^{2}([0,\infty)), G'(0)>G(0)=0~\text{and}~G(u)>0,G'(u)\geq0~\text{for}~u>0, \\[2mm]
\frac{G(u)}{u}~\text{is~strongly~decreasing}~\text{for}~u>0~\text{and}~0<\frac{G(u)}{u}<1~\text{for}~u>0;
\end{array} \right.
\end{eqnarray*}
\item[(A4)] For each $i\in\{1,2\}$, the constants $H_{i}>0$, $\varpi_{i}>0$, and $\kappa_{i}>1$ can be found such that $\rho_{i}(u)\geq \rho'_{i}(0)u-H_{i}u^{\kappa_{i}}$ for $0\leq u\leq\varpi_{i}$.
\end{itemize}

It should be mentioned that (A1) and (A2) are conventional assumptions about the initial functions $u_{0}(x)$ and $v_{0}(x)$ and the growth function $f(u)$, respectively, see, for example, \cite{Capasso-Maddalena-3, Wang-Du} and the references therein. (A3) is a natural assumption about the impulse function $G(u)$, however (A4) is a technical assumption needed in the construction of the lower solution to the impulse problem, which was first introduced by Lewis and Li \cite{Lewis-Li}. In \cite{Zhou-Lin-Santos, Li-Zhao-Cheng}, some impulse and growth functions that satisfy the assumption conditions are given for the interested reader. In the sequel, Assumptions (A1)-(A4) always hold without further statement.

Now, an impulsive faecal-oral model with free boundary is developed. To achieve a realistic scenareo, the global classical solution of the newly developed model must exist and be unique and nonnegative. This problem will be studied below.

In \cite{Wang-Du}, Wang and Du gave some properties of the global solution of model \eqref{Wang-Du}. Here we will mainly address the difficulties caused by the introduction of the impulse intervention in the process of proving the well-posedness of the solution to model \eqref{Zhou-Lin}.
To conveniently present the solution of the impulsive differential equation, we introduce the following notations:
\begin{equation}
\begin{aligned}
\mathbb{PC}^{\mu}_{t}:=\mathbb{PC}^{\mu}_{t}\big(\mathbb{R}^{+}\big)&:=\big\{u:u\in \cap_{k\in\mathbb{N}}\mathbb{C}^{\mu}((k\tau, (k+1)\tau])\big\};\\
\mathbb{PC}_{t,x}:=\mathbb{PC}_{t,x}\big(\mathbb{R}^{+}\times[g(t), h(t)]\big)&:=\big\{u:u\in \cap_{k\in\mathbb{N}}\mathbb{C}\big((k\tau, (k+1)\tau]\times [g(t), h(t)] \big)\big\};\\
\mathbb{PC}^{p,q}_{t,x}:=\mathbb{PC}^{p,q}_{t,x}\big(\mathbb{R}^{+}\times[g(t), h(t)]\big)&:=\big\{u:u\in \cap_{k\in\mathbb{N}}\mathbb{C}^{p,q}((k\tau, (k+1)\tau]\times [g(t), h(t)])\big\};
\nonumber
\end{aligned}
\end{equation}
where $\mu$, $p$, $q>0$. Based on the work \cite{Wang-Du}, we next present the existence and uniqueness of a local solution to model \eqref{Zhou-Lin}.
\begin{lemma}\label{lemma 2-1}
For arbitrary given initial functions $u_{0}(x)$ and $v_{0}(x)$ and arbitrary $\alpha\in(0,1)$, there exists a $T\in(0,\tau)$ such that
model \eqref{Zhou-Lin} has a unique solution
\begin{equation*}
(u,v,g,h)\in\big[\mathbb{PC}^{1+\alpha/2,2+\alpha}_{t,x}\big]^{2}\times\big[\mathbb{PC}^{1+(1+\alpha)/2}_{t}\big]^{2}.
\end{equation*}
\begin{proof}
Since $u_{0}(x)$ and $v_{0}(x)$ and $G(u)$ satisfy Assumptions (A1) and (A3), we can be obtain that
\begin{eqnarray*}
\left\{
\begin{array}{l}
u(0^{+}, x)\in\mathbb{C}^{2}([-h_{0},h_{0}]), ~u(0^{+}, \pm h_{0})=0~\text{and}~u(0^{+}, x)>0~\text{in}~(-h_{0},h_{0}), \\[2mm]
v(0^{+}, x)\in\mathbb{C}^{2}([-h_{0},h_{0}]), ~v(0^{+}, \pm h_{0})=0~\text{and}~v(0^{+}, x)>0~~\text{in}~(-h_{0},h_{0}).
\end{array} \right.
\end{eqnarray*}
Taking $u(0^{+}, x)$ and $v(0^{+}, x)$ as new initial functions, the conclusion can then be obtained by employing similar techniques as in \cite[Theorem 2.1]{Wang-Du}. This ends the proof.
\end{proof}
\end{lemma}
In order to prove that the local solution can be extended to the global one, the following estimates are given.
\begin{lemma}\label{lemma 2-2}
Let $(u,v,g,h)$ be a solution of model \eqref{Zhou-Lin} defined for $t\in\mathbb{R}^{+}$. Then there exist constants $C_{1}$, $C_{2}$ and $C_{3}$
independent of $k$ such that
\begin{equation}\label{ws-0}
\begin{aligned}
0&<-g'(t), h'(t)< C_{1} \text{~for~}t\in\mathbb{R}^{+},\\
0<u(t,x)<C_{2}&,~0<v(t,x)< C_{3} \text{~for~}t\in\mathbb{R}^{+}\text{~and~}x\in(g(t), h(t)).
\end{aligned}
\end{equation}
\begin{proof}
In order to facilitate the reading, the proof of this lemma is divided into two parts.

\textbf{(1) The case of $0\bm{<}\bm{t}$$\bm{\leq}\tau$}

Taking $u(0^{+}, x)$ and $v(0^{+}, x)$ as initial functions, it follows from Assumption (A2) and the strong maximum principle for cooperative systems (see \cite[Theorem 13 in Chapter 3]{Protter-Weinberger}) that $u, v>0$ for $t\in(0,\tau]$ and $x\in(g(t), h(t))$. Then, using the Hopf boundary lemma
for cooperative systems (see \cite[Theorem 14 in Chapter 3]{Protter-Weinberger}), one can obtain that
\begin{equation*}
u_{x}(t,h(t)), v_{x}(t, h(t)), -u_{x}(t,g(t)), -v_{x}(t,g(t))<0
\end{equation*}
for $t\in(0,\tau]$. This means that $-g'(t), h'(t)>0$ for $t\in(0,\tau]$.

From Assumption (A2), it follows that for a fixed sufficiently small positive constant $\epsilon_{0}$, there exists a sufficiently large constant
\begin{equation*}
M_{0}>\max\bigg\{\max\limits_{\bar{\Omega}}u_{0}(x), \frac{a_{12}+\epsilon_{0}}{a_{11}}\max\limits_{\bar{\Omega}}v_{0}(x)\bigg\}
\end{equation*}
such that when $M>M_{0}$, we have
\begin{equation*}
f(M)< \frac{a_{11}a_{22}}{a_{12}+\epsilon_{0}}M.
\end{equation*}
Define $C_{2}:=M$, $C_{3}:=\frac{a_{11}}{a_{12}+\epsilon_{0}}M$, $N_{1}:=\max\limits_{[0, \tau)\times\Omega}u(t,x)$, and $N_{2}:=\max\limits_{0\leq u\leq\max(C_{2},N_{1})}f'(u)$.
Let
\begin{equation}\label{WS-1}
\zeta=(C_{2}-u)e^{-(a_{12}+N_{2})t}~\text{and}~\eta=(C_{3}-v)e^{-(a_{12}+N_{2})t}.
\end{equation}
By standard calculation, it follows that
\begin{equation}\label{WS-2}
\begin{array}{ll}
&\zeta_{t}-d_{1}\Delta \zeta+(a_{11}+a_{12}+N_{2})\zeta-a_{12}\eta\\
=&\Big[-u_{t}-(C_{2}-u)(a_{12}+N_{2})+d_{1}\Delta u-a_{12}(C_{3}-v)\\
&+(a_{11}+a_{12}+N_{2})(C_{2}-u)\Big]e^{-(a_{12}+N_{2})t}\\
=&\Big[a_{11}C_{2}-a_{12}C_{3}\Big]e^{-(a_{12}+N_{2})t}>0
\end{array}
\end{equation}
and
\begin{equation}\label{WS-3}
\begin{array}{ll}
&\eta_{t}-d_{2}\Delta \eta+(a_{22}+a_{12}+N_{2})\eta-f'(\vartheta)\zeta\\
=&\Big[-v_{t}-(C_{3}-v)(a_{12}+N_{2})+d_{2}\Delta v-f'(\vartheta)(C_{2}-u)\\
&+(a_{22}+a_{12}+N_{2})(C_{3}-v)\Big]e^{-(a_{12}+N_{2})t}\\
=&\Big[a_{22}C_{3}-f(C_{2})\Big]e^{-(a_{12}+N_{2})t}>0,
\end{array}
\end{equation}
where $\min(C_{2}, u)\leq\vartheta\leq\max(C_{2}, u)$ for $(t,x)\in(0, \tau]\times(g(t),h(t))$. Thus, by the assumption that $G(u)>0$ for $u>0$ in (A3) and \eqref{WS-1}-\eqref{WS-3}, it can be obtained that
\begin{eqnarray*}
\left\{
\begin{array}{ll}
\zeta_{t}>d_{1}\Delta \zeta-(a_{11}+a_{12}+N_{2})\zeta+a_{12}\eta,\; &\, t\in(0^{+}, \tau], x\in(g(t),h(t)), \\[2mm]
\eta_{t}>d_{2}\Delta \eta-(a_{22}+a_{12}+N_{2})\eta+f'(\vartheta)\zeta,\; &\, t\in(0^{+}, \tau], x\in(g(t),h(t)), \\[2mm]
\zeta(t,x)>0, \eta(t,x)>0,\; &\, t\in(0^{+}, \tau], x\in\{g(t), h(t)\},\\[2mm]
\zeta(0^{+},x)>0, \eta(0^{+},x)>0, \; &\, x\in[-h_{0}, h_{0}].
\end{array} \right.
\end{eqnarray*}
Then, the strong maximum principle for cooperative systems yields that $\zeta, \eta> 0$ for $t\in(0, \tau]$ and $x\in(g(t),h(t))$.
This implies that $u<C_{2}$ and $v< C_{3}$ for $(t,x)\in(0, \tau]\times (g(t),h(t))$.

Next, we prove that $h'(t), -g'(t)< C_{1}$ for $t\in(0, \tau]$. Define
\begin{equation*}
\Omega_{\Lambda}=\{(t,x):t\in(0^{+}, \tau], ~x\in(h(t)-\Lambda^{-1},h(t))\},
\end{equation*}
and
\begin{equation*}
\phi=\varphi=C\cdot\omega:=C[2\Lambda(h(t)-x)-\Lambda^{2}(h(t)-x)^{2}],
\end{equation*}
where constants $\Lambda$ and $C$ are to be specified later. By taking $\Lambda:=\max\Big\{\sqrt{\frac{a_{12}}{2d_{1}}}, \sqrt{\frac{f'(0)}{2d_{2}}}\Big\}$
and choosing $C>\max\{C_{2},C_{3}\}$, it follows that
\begin{eqnarray*}
\left\{
\begin{array}{ll}
\phi_{t}-d_{1}\Delta\phi\geq2d_{1}C\Lambda^{2}\geq a_{12}v\geq u_{t}-d_{1}\Delta u,\; &\, (t,x)\in\Omega_{\Lambda}, \\[2mm]
\varphi_{t}-d_{2}\Delta\varphi\geq2d_{2}C\Lambda^{2}\geq f(u)\geq v_{t}-d_{2}\Delta v,\; &\, (t,x)\in\Omega_{\Lambda}.
\end{array} \right.
\end{eqnarray*}
By enlarging $C$ further if necessary, we have
\begin{eqnarray*}
\left\{
\begin{array}{ll}
u(t, h(t)-\Lambda^{-1})\leq C\leq \phi(t, h(t)-\Lambda^{-1}),\; &\, t\in(0^{+}, \tau], \\[2mm]
v(t, h(t)-\Lambda^{-1})\leq C\leq \varphi(t, h(t)-\Lambda^{-1}),\; &\, t\in(0^{+}, \tau], \\[2mm]
\phi(t, h(t))=\varphi(t, h(t))=0,\; &\, t\in(0^{+}, \tau], \\[2mm]
u(t, h(t)-\Lambda^{-1})\leq C\leq \phi(t, h(t)-\Lambda^{-1}),\; &\, t\in(0^{+}, \tau], \\[2mm]
u(0^{+},x)\leq \phi(0^{+},x), v(0^{+},x)\leq \varphi(0^{+},x),\; &\, x\in[h_{0}-\Lambda^{-1}, h_{0}].
\end{array} \right.
\end{eqnarray*}
Applying the maximum principle to $\phi-u$ and $\varphi-v$, it can be obtained that $u\leq\phi$ and $v\leq\varphi$ for $(t,x)\in
\overline{\Omega}_{\Lambda}$. This implies that
\begin{equation*}
\begin{split}
h'(t)&=-\mu_{1}u_{x}(t,h(t))-\mu_{2}v_{x}(t,h(t))\leq -\mu_{1}\phi_{x}(t,h(t))-\mu_{2}\varphi_{x}(t,h(t))\\
&<2\Lambda C(\mu_{1}+\mu_{2}+1):=C_{1}.
\end{split}
\end{equation*}
A similar method can prove that $-g'(t)< C_{1}$.

\textbf{(2) The case of $\bm{t}$$\bm{>}\bm{\tau}$}

In the light of Case (1) and the assumption that $0<\frac{G(u)}{u}<1$ for $u>0$ in (A3), $u(\tau^{+}, x)$ and $v(\tau^{+}, x)$ as new initial functions satisfy the conditions $u(\tau^{+}, x)<C_{2}$ and $v(\tau^{+}, x)<C_{3}$ for $x\in(g(t),h(t))$, and $u(\tau^{+}, x)\leq\phi(\tau^{+},x)$ and $v(\tau^{+}, x)\leq\varphi(\tau^{+},x)$ for $x\in(h(\tau)-\Lambda^{-1}, h(\tau))$. Repeating the procedure in Case (1) yields that $0<-g'(t), h'(t)< C_{1}$ for $t\in (\tau, 2\tau]$, and
$0<u<C_{2}$ and $0<v< C_{3}$ for $(t,x)\in (\tau, 2\tau]\times (g(t), h(t))$. Step by step, \eqref{ws-0} can be obtained. This completes the proof.
\end{proof}
\end{lemma}
Next, we are ready to extend the solution obtained in \autoref{lemma 2-1} to all $t>0$.
\begin{theorem}\label{theorem 2-1}
For any initial function $(u_{0}(x), v_{0}(x))$, model \eqref{Zhou-Lin} has a unique nonnegative solution $(u,v,g,h)$ for all $t\in\mathbb{R}^{+}$. Furthermore,
\begin{equation*}
(u,v,g,h)\in\big[\mathbb{PC}^{1,2}_{t,x}\cap \mathbb{PC}_{t,x}\big]^{2}\times \big[\mathbb{PC}^{1}_{t}\cap \mathbb{C}((0,\infty))\big]^{2}.
\end{equation*}
\begin{proof}
Suppose $t_{\max}$ is the maximum existence time of the solution. \autoref{lemma 2-1} yields that $t_{\max}>0$. Next we proceed  by contradiction to show that $t_{\max}=\infty$. Assume that $t_{\max}<\infty$. This assumption implies that there exists $k_{0}\in\mathbb{N}$ such that $t_{\max}\in((k_{0}\tau)^{+}, (k_{0}+1)\tau]$. We consider the following two cases.

(1) When $t_{\max}\in((k_{0}\tau)^{+}, (k_{0}+1)\tau)$, $(u((k_{0}\tau)^{+},x),v((k_{0}\tau)^{+},x))$ is viewed as a new initial vector-valued function.
From \autoref{lemma 2-2} and \cite[Theorem 1.1]{Wang-Du}, one can obtain that
\begin{equation*}
(u, v)\in[\mathbb{C}^{1,2}((k_{0}\tau, (k_{0}+1)\tau]\times [g(t), h(t)])]^{2},
\end{equation*}
which is a contradiction to the definition of $t_{\max}$.

(2) When $t_{\max}=(k_{0}+1)\tau$, we have that
\begin{equation*}
(u, v)\in[\mathbb{C}^{1,2}((k_{0}\tau, (k_{0}+1)\tau)\times [g(t), h(t)])]^{2}.
\end{equation*}
From  \cite[Theorem 1.1]{Wang-Du}, it follows that
\begin{equation*}
(u((k_{0}+1)\tau, x), v((k_{0}+1)\tau, x))\in[\mathbb{C}^{2} [g((k_{0}+1)\tau), h((k_{0}+1)\tau)]]^{2}.
\end{equation*}
Taking $(u(((k_{0}+1)\tau)^{+},x),v(((k_{0}+1)\tau)^{+},x))$ as a new initial function and repeating the procedure in (1) yields that the time interval in which the solution exists is $[0,(k_{0}+2)\tau ]$, which is also a contradiction to the definition of $t_{\max}$.

In summary, model \eqref{Zhou-Lin} admits a unique nonnegative global classical solution. This ends the proof.
\end{proof}
\end{theorem}
\section{\bf Spreading and vanishing}\label{Section-3}
This section first introduces some results of the corresponding model on a fixed region. Based on these results, the spreading-vanishing dichotomy
is established and some sufficient conditions governing spreading or vanishing are provided.
\subsection{\bf Results for the fixed region model}\label{Section-3.1}
Here, we introduce some main results for the impulsive faecal-oral model on a fixed region, which will be used in the sequel.

First we introduce the impulsive faecal-oral model on the one-dimensional region $(s_{1},s_{2})$,l
\begin{eqnarray}\label{3-1}
\left\{
\begin{array}{ll}
u_{t}=d_{1}\Delta u-a_{11}u+a_{12}v,\; &\, t\in((k\tau)^{+}, (k+1)\tau], x\in(s_{1},s_{2}), \\[2mm]
v_{t}=d_{2}\Delta v-a_{22}v+f(u),\; &\, t\in((k\tau)^{+}, (k+1)\tau], x\in(s_{1},s_{2}), \\[2mm]
u=0, v=0,\; &\, t\in(k\tau, (k+1)\tau], ~~~~~x\in\{s_{1},s_{2}\},\\[2mm]
u((k\tau)^{+},x)=G(u(k\tau,x)), \; &\, x\in(s_{1},s_{2}), \\[2mm]
v((k\tau)^{+},x)=v((k\tau),x), \; &\,  x\in(s_{1},s_{2}), ~~~~~~~k=0,1,2,\cdots
\end{array} \right.
\end{eqnarray}
with initial functions $u(0,x)=u_{0}(x)$ and $v(0,x)=v_{0}(x)$ for $x\in[s_{1}, s_{2}]$ and the corresponding periodic eigenvalue problem
\begin{eqnarray}\label{3-2}
\left\{
\begin{array}{ll}
\phi_{t}=d_{1}\Delta \phi-a_{11}\phi+a_{12}\psi+\lambda \phi,\; &\, t\in(0^{+}, \tau], x\in(s_{1},s_{2}), \\[2mm]
\psi_{t}=d_{2}\Delta \psi-a_{22}\psi+f'(0)\phi+\lambda \psi,\; &\, t\in(0^{+}, \tau], x\in(s_{1},s_{2}), \\[2mm]
\phi(t,x)=0, \psi(t,x)=0,\; &\, t\in[0, \tau],~~x\in\{s_{1},s_{2}\},\\[2mm]
\phi(0^{+},x)=G'(0)\phi(0,x), \psi(0^{+},x)=\psi(0,x), \; &\,  x\in(s_{1},s_{2})
\end{array} \right.
\end{eqnarray}
with periodic conditions $\phi(0,x)=\phi(\tau,x)$ and $\psi(0,x)=\psi(\tau,x)$ for $x\in [s_{1}, s_{2}]$. Since model \eqref{3-1} is a special form of model \textcolor{blue}{(2.1)} in \cite{Zhou-Lin-Santos}, the periodic eigenvalue problem \eqref{3-2} has a principal eigenvalue $\lambda^{\vartriangle}_{1}$ with a strongly positive eigenfunction pair $(\phi, \psi)$.
In order to emphasize the dependence of $\lambda^{\vartriangle}_{1}$  on $G'(0)$ and $(s_{1},s_{2})$, we also write $\lambda^{\vartriangle}_{1}=\lambda^{\vartriangle}_{1}(G'(0),(s_{1},s_{2}))$. Next, some properties of the principal eigenvalue $\lambda^{\vartriangle}_{1}$ and threshold-type results on the global dynamic of model \eqref{3-1} are provided.
\begin{lemma}\label{lemma 3-1}
The assertions below are valid.
\begin{enumerate}
\item[$(1)$]
If $s_{1}<s_{2}$, then we have
\begin{equation*}
\lambda^{\vartriangle}_{1}(G'(0),(s_{1},s_{2}))=\lambda^{\vartriangle}_{1}(G'(0),(0,s_{2}-s_{1})).
\end{equation*}
\item[$(2)$]
If $G'_{1}(0)<G'_{2}(0)$, then we have
\begin{equation*}
\lambda^{\vartriangle}_{1}(G_{2}'(0),(s_{1},s_{2}))<\lambda^{\vartriangle}_{1}(G'_{1}(0),(s_{1},s_{2})).
\end{equation*}
\item[$(3)$]
If $s_{1}<s_{2}$, then we have
\begin{equation*}
\lambda^{\vartriangle}_{1}(G'(0),(0,s_{2}))<\lambda^{\vartriangle}_{1}(G'(0),(0,s_{1})).
\end{equation*}
\item[$(4)$]
If $\lambda^{\vartriangle}_{1}\geq0$, then the solution of model \eqref{3-1} satisfies
\begin{equation*}
\lim\limits_{t\rightarrow+\infty}\|u(t, \cdot)\|_{\mathbb{C}[s_{1},s_{2}]}=\lim\limits_{t\rightarrow+\infty}\|v(t, \cdot)\|_{\mathbb{C}[s_{1},s_{2}]}=0.
\end{equation*}
\item[$(5)$]
If $\lambda^{\vartriangle}_{1}<0$, then the solution of model \eqref{3-1} satisfies
\begin{equation*}
\lim\limits_{m\rightarrow+\infty}(u(t+m\tau,x),v(t+m\tau,x))=(U(t,x), V(t,x))
\end{equation*}
uniformly for $(t,x)\in[0,\infty)\times[s_{1},s_{2}]$. Here $(U, V)$ stands for the one and only positive periodic solution of the problem
\begin{eqnarray}\label{SFC-1}
\left\{
\begin{array}{ll}
u_{t}=d_{1} \Delta u-a_{11}u+a_{12}v,\; &\, t\in(0^{+}, \tau], x\in(s_{1},s_{2}), \\[2mm]
v_{t}=d_{2} \Delta v-a_{22}v+f(u),\; &\, t\in(0^{+}, \tau], x\in(s_{1},s_{2}), \\[2mm]
u(t,x)=0, v(t,x)=0,\; &\, t\in[0, \tau],~~ x\in\{s_{1},s_{2}\}, \\[2mm]
u(0,x)=u(\tau,x), v(0,x)=v(\tau,x),\; &\,x\in[s_{1},s_{2}],\\[2mm]
u(0^{+},x)=G(u(0,x)), v(0^{+},x)=v(0,x), \; &\, x\in(s_{1},s_{2}).
\end{array} \right.
\end{eqnarray}
\end{enumerate}
\begin{proof}
By applying a translational transformation to the space variable $x$, assertion (1) can be proved directly. Assertions (2, 3, 5) are straightforward results of Theorem 5.2(1), Theorem 5.2(2), Theorem 4.3 in \cite{Zhou-Lin-Santos}, respectively.

Next, we focus on proving assertion (4). When $\lambda^{\vartriangle}_{1}>0$, the desired conclusion is a direct result of \cite[Theorem 4.1]{Zhou-Lin-Santos}.
Here, the proof that the conclusion holds when $\lambda^{\vartriangle}_{1}=0$ is given. Firstly we show that when $\lambda^{\vartriangle}_{1}=0$, the problem \eqref{SFC-1} has a unique nonnegative periodic solution $(0,0)$. To achieve this, we introduce the adjoint problem corresponding to problem \eqref{3-2}, that is,
\begin{eqnarray}\label{SFC-2}
\left\{
\begin{array}{ll}
-\frac{\partial \phi^{T}}{\partial t}=d_{1}\Delta \phi^{T}-a_{11}\phi^{T}+f'(0)\psi^{T}+\mu \phi^{T},\; &\, t\in(0^{+}, \tau], x\in(s_{1},s_{2}), \\[2mm]
-\frac{\partial \psi^{T}}{\partial t}=d_{2}\Delta \psi^{T}-a_{22}\psi^{T}+a_{12}\phi^{T}+\mu \psi^{T},\; &\, t\in(0^{+}, \tau], x\in(s_{1},s_{2}), \\[2mm]
\phi^{T}=0, \psi^{T}=0,\; &\, t\in[0, \tau],~~ x\in\{s_{1},s_{2}\}, \\[2mm]
\phi^{T}(0,x)=\phi^{T}(\tau,x), \psi^{T}(0,x)=\psi^{T}(\tau,x),\; &\,x\in[s_{1},s_{2}],\\[2mm]
\phi^{T}(0^{+},x)=\frac{1}{G'(0)}\phi^{T}(0,x), \psi^{T}(0^{+},x)=\psi^{T}(0,x), \; &\, x\in(s_{1},s_{2}).
\end{array} \right.
\end{eqnarray}

Obviously, $(0,0)$ is a solution to problem \eqref{SFC-1}. On the contrary, suppose that $(u,v)$ is a nonnegative periodic solution to problem \eqref{SFC-1} and satisfies the condition $u, v\not\equiv 0$ for $(t,x)\in(0^{+}, \tau]\times(s_{1},s_{2})$. By multiplying the first equations of problem \eqref{SFC-1} and problem \eqref{SFC-2} by $\phi^{T}$ and $u$, respectively, it follows that
\begin{eqnarray}\label{SFC-3}
\left\{
\begin{array}{ll}
\frac{\partial u}{\partial t}\phi^{T}=d_{1}\Delta u\phi^{T}-a_{11}u\phi^{T}+a_{12}v\phi^{T},\; &\, t\in(0^{+}, \tau], x\in(s_{1},s_{2}), \\[2mm]
\frac{\partial \phi^{T}}{\partial t}u=-d_{1}\Delta \phi^{T}u+a_{11}\phi^{T}u-f'(0)\psi^{T}u-\mu \phi^{T}u,\; &\, t\in(0^{+}, \tau], x\in(s_{1},s_{2}),
\end{array}
\right.
\end{eqnarray}
Then, integrating each side of the first and second equations of \eqref{SFC-3} on $(0^{+}, \tau]\times(s_{1},s_{2})$, respectively,
and adding the obtained results give that
\begin{equation}\label{SFC-4}
\begin{aligned}
\int_{s_{1}}^{s_{2}}\int^{\tau}_{0^{+}}\Big(\frac{\partial u}{\partial t}\phi^{T}+\frac{\partial\phi^{T}}{\partial t} u\Big)dtdx=&d_{1}\int_{s_{1}}^{s_{2}}\int^{\tau}_{0^{+}}\big(\Delta u\phi^{T}-\Delta \phi^{T}u\big)dtdx-\mu\int_{s_{1}}^{s_{2}}\int^{\tau}_{0^{+}}\phi^{T}udtdx\\
=&\int_{s_{1}}^{s_{2}}\int^{\tau}_{0^{+}}\big(a_{12}v\phi^{T}-f'(0)\psi^{T}u\big)dtdx.
\end{aligned}
\end{equation}
By using Fubini's theorem, \cite[Theorem5.1]{Zhou-Lin-Santos}, condition $\lambda^{\vartriangle}_{1}=0$, and the formula of integration by parts,
\eqref{SFC-4} can be simplified to
\begin{equation}\label{SFC-5}
\begin{aligned}
\int_{s_{1}}^{s_{2}}\int^{\tau}_{0^{+}}\Big(\frac{\partial u}{\partial t}\phi^{T}+\frac{\partial\phi^{T}}{\partial t} u\Big)dtdx=&\int_{s_{1}}^{s_{2}}\int^{\tau}_{0^{+}}\big(a_{12}v\phi^{T}-f'(0)\psi^{T}u\big)dtdx.
\end{aligned}
\end{equation}
A similar procedure can yield
\begin{equation}\label{SFC-6}
\begin{aligned}
\int_{s_{1}}^{s_{2}}\int^{\tau}_{0^{+}}\Big(\frac{\partial v}{\partial t}\psi^{T}+\frac{\partial\psi^{T}}{\partial t} v\Big)dtdx=&\int_{s_{1}}^{s_{2}}\int^{\tau}_{0^{+}}\big(f(u)\psi^{T}-a_{12}v\phi^{T}\big)dtdx.
\end{aligned}
\end{equation}
Adding \eqref{SFC-5} and \eqref{SFC-6} gives that
\begin{equation*}
\begin{aligned}
\int_{s_{1}}^{s_{2}}\big(u\phi^{T}\big)\big|_{0^{+}}^{\tau}dx=\int_{s_{1}}^{s_{2}}\int^{\tau}_{0^{+}}\big(f(u)\psi^{T}-f'(0)\psi^{T}u\big)dtdx.
\end{aligned}
\end{equation*}
With the use of the period conditions, impulse conditions and Assumption (A3), it follows that
\begin{equation}\label{SFC-7}
\begin{aligned}
\int_{s_{1}}^{s_{2}}\big(u\phi^{T}\big)\big|_{0^{+}}^{\tau}dx=\int_{s_{1}}^{s_{2}}\big(u(0,x)\phi^{T}(0,x)-\frac{G(u(0,x))}{G'(0)}\phi^{T}(0,x)\big)dx\geq0.
\end{aligned}
\end{equation}
On the other hand, Assumption (A3) yields that
\begin{equation*}
\begin{aligned}
\int_{s_{1}}^{s_{2}}\int^{\tau}_{0^{+}}\big(f(u)\psi^{T}-f'(0)\psi^{T}u\big)dtdx=\int_{s_{1}}^{s_{2}}\int^{\tau}_{0^{+}}\psi^{T}\big(f(u)-f'(0)u\big)< 0,
\end{aligned}
\end{equation*}
which show a contradiction with \eqref{SFC-7}. Therefore, the assumption is not valid, i.e., $(0,0)$ is the one and only nonnegative
periodic solution to problem \eqref{SFC-1} when $\lambda^{\vartriangle}_{1}=0$.

Notice that $(C_{2}, C_{3})$ is a supersolution of problem \eqref{SFC-1}, where $(C_{2}, C_{3})$ is defined in \autoref{lemma 2-2}. Next, choosing $(\overline{u}, \overline{v})$ as the initial iteration and using the iteration rule in \cite[(4.3)]{Zhou-Lin-Santos}, the monotonically decreasing sequence $\big\{(\overline{u}^{(i)}, \overline{v}^{(i)})\big\}_{i=0}^{n}$ can be obtained. Using the same method that was used for proving  \cite[Theorem 4.3]{Zhou-Lin-Santos}, it follows that
\begin{equation*}
\begin{aligned}
\limsup\limits_{k\rightarrow+\infty}\big(u(t+k\tau,x),v(t+k\tau,x)\big)\leq\big(U(t,x), V(t,x)\big)
\end{aligned}
\end{equation*}
for all $(t,x)\in[0,\infty)\times[s_{1}, s_{2}]$. Here $(U, V)$ stands for the one and only nonnegative periodic solution to problem \eqref{SFC-1}. Therefore, it follows that
\begin{equation*}
\lim\limits_{t\rightarrow+\infty}\|u(t, \cdot)\|_{\mathbb{C}[s_{1},s_{2}]}=\lim\limits_{t\rightarrow+\infty}\|v(t, \cdot)\|_{\mathbb{C}[s_{1},s_{2}]}=0.
\end{equation*}
This completes the proof of assertion (4).
\end{proof}
\end{lemma}
\subsection{\bf Spreading-vanishing dichotomy}\label{Section-3.2}
This subsection shows that the infected individuals and pathogenic bacteria either spread or vanish when time tends to infinity.

With the help of \autoref{lemma 2-2},
\begin{equation*}
\lim\limits_{t\rightarrow\infty}g(t)=g_{\infty}\in[-\infty,-h_{0})\text{~and~}\lim\limits_{t\rightarrow\infty}h(t)=h_{\infty}\in(h_{0}, +\infty]
\end{equation*}
are always well-defined. To start with, the following definition about spreading and vanishing is provided.
\begin{definition}
The infected individuals and pathogenic bacteria are spreading if
\begin{equation*}
h_{\infty}-g_{\infty}=\infty\text{~and~}\lim\limits_{t\rightarrow\infty}\|u(t, \cdot)\|_{\mathbb{C}(\Omega)}+\lim\limits_{t\rightarrow\infty}\|v(t, \cdot)\|_{\mathbb{C}(\Omega)}>0
\end{equation*}
for any $\Omega\subset\subset (g_{\infty}, h_{\infty})$, and vanishing if
\begin{equation*}
h_{\infty}-g_{\infty} \leq \infty\text{~and~}\lim\limits_{t\rightarrow\infty}\|u(t, \cdot)\|_{\mathbb{C}[g(t), h(t)]}+\lim\limits_{t\rightarrow\infty}\|v(t, \cdot)\|_{\mathbb{C}[g(t), h(t)]}=0.
\end{equation*}
\end{definition}
In order to estimate $(u,v,g,h)$  conveniently later, the comparison principle is now presented.
\begin{lemma}\label{lemma 3-2}
Suppose that
\begin{equation*}
(\overline{u},\overline{v},\overline{g},\overline{h})\in\Big[\mathbb{PC}^{1,2}_{t,x}\big(\mathbb{R}^{+}\times[\overline{g}(t), \overline{h}(t)]\big)\cap \mathbb{PC}_{t,x}\big(\mathbb{R^{+}}\times[\overline{g}(t), \overline{h}(t)]\big)\Big]^{2}\times \Big[\mathbb{PC}^{1}_{t}\cap \mathbb{C}(\mathbb{R}^{+})\Big]^{2},
\end{equation*}
and that
\begin{eqnarray*}
\left\{
\begin{array}{ll}
\overline{u}_{t}\geq d_{1}\Delta \overline{u}-a_{11}\overline{u}+a_{12}\overline{v},\; &\, t\in((k\tau)^{+}, (k+1)\tau], x\in(\overline{g}(t), \overline{h}(t)), \\[2mm]
\overline{v}_{t}\geq d_{2}\Delta \overline{v}-a_{22}\overline{v}+f(\overline{u}),\; &\, t\in((k\tau)^{+}, (k+1)\tau], x\in(\overline{g}(t),
\overline{h}(t)), \\[2mm]
\overline{u}=0,  \overline{v}=0,\; &\, t\in(k\tau, (k+1)\tau], ~~x\in\partial(\overline{g}(t),  \overline{h}(t)),\\[2mm]
\overline{u}((k\tau)^{+},x)\geq G(\overline{u}(k\tau,x)), \; &\, x\in(\overline{g}(k\tau),  \overline{h}(k\tau)), \\[2mm]
 \overline{v}((k\tau)^{+},x)\geq \overline{v}((k\tau),x), \; &\,  x\in(\overline{g}(k\tau),  \overline{h}(k\tau)), \\[2mm]
\overline{g}'(t)\leq-\mu_{1}\overline{u}_{x}(t,\overline{g}(t))-\mu_{2}\overline{v}_{x}(t,\overline{g}(t)),\; &\, t\in(k\tau, (k+1)\tau],\\[2mm]
\overline{h}'(t)\geq-\mu_{1}\overline{u}_{x}(t,\overline{h}(t))-\mu_{2}\overline{v}_{x}(t,\overline{h}(t)),\; &\, t\in(k\tau, (k+1)\tau],\\[2mm]
\overline{g}(0)\leq-h_{0}, u_{0}(x)\leq\overline{u}(0,x),\; &\,x\in[-h_{0},h_{0}],\\[2mm]
\overline{h}(0)\geq h_{0},~~~ v_{0}(x)\leq\overline{v}(0,x),\; &\,x\in[-h_{0},h_{0}], k=0,1,2,3,4,\cdots.
\end{array} \right.
\end{eqnarray*}
Then, the solution $(u,v,g,h)$ of model \eqref{Zhou-Lin} satisfies
\begin{eqnarray*}
\left\{
\begin{array}{l}
g(t)\geq\overline{g}(t), ~~~~~~~~\overline{h}(t)\geq h(t),~~~~~~~\text{~for~}t\in(0,\infty), \\[2mm]
\overline{u}(t, x)\geq u(t,x), \overline{v}(t, x)\geq v(t,x),\text{~for~}(t,x)\in(0,\infty)\times[g(t), h(t)].
\end{array} \right.
\end{eqnarray*}
\begin{proof}
Since $(u_{0}(x),v_{0}(x)) \leq(\overline{u}(0,x),\overline{v}(0,x))$ for $x\in[-h_{0},h_{0}]$, it follows from the assumption $G'(u)\geq 0$ for $u\geq 0$ in (A3) that
\begin{equation*}
(u(0^{+}, x), v(0^{+}, x) )\leq(\overline{u}(0^{+},x),\overline{v}(0^{+},x) ), x\in[-h_{0},h_{0}].
\end{equation*}
Then, \cite[Lemma 2.3]{Wang-Du} yields that
\begin{eqnarray}\label{3-3}
\left\{
\begin{array}{l}
\Gamma(t)=h\kappa(t), \kappa(t)=1+\sigma-\frac{\sigma}{2}e^{-\gamma t}, t\geq0, \\[2mm]
(u,v )\leq(\overline{u}, \overline{v}), (t,x)\in(0^{+},\tau]\times[g(t), h(t)].
\end{array} \right.
\end{eqnarray}
Next, choosing $u|_{t=\tau}$ and $v|_{t=\tau}$ as new initial functions, it can be obtained through a similar procedure that \eqref{3-3} holds for the time interval $(\tau^{+}, 2\tau]$. Step by step, \eqref{3-3} holds for all $t\in(0,\infty)$. This ends the proof.
\end{proof}
\end{lemma}
\begin{remark}\label{Remark 3-1}
The quadruple $(\overline{u},\overline{v},\overline{g},\overline{h})$ in \autoref{lemma 3-2} is usually called an upper solution of model \eqref{Zhou-Lin}. Similarly, a lower solution $(\underline{u},\underline{v},\underline{g},\underline{h})$ of model \eqref{Zhou-Lin} can also be defined by reversing all the inequalities, and there is an analogue of \autoref{lemma 3-2} for lower solutions.
\end{remark}

It is well known that the method of separation of variables is a very fundamental and common method for solving partial differential equations. This method can separate the terms of the equation containing each variable, thus separating the original equation into several simpler ordinary differential equations containing only one variable. Now, we apply the method of separation of variables to the eigenvalue problem \eqref{3-2}. Let
\begin{equation*}
\begin{aligned}
(\phi(t,x), \psi(t,x))=X(x)(\Phi(t), \Psi(t)), (t,x)\in[0,\tau]\times[-h,h],
\end{aligned}
\end{equation*}
where $\big(\lambda_{0}, X(x)\big)$ is the principal eigenpair of the eigenvalue problem
\begin{eqnarray*}
\left\{
\begin{array}{ll}
\Delta X(x)=-\lambda X(x)\; &\, \text{~in~}(-h,h), \\[2mm]
X(x)=0\; &\, \text{~on~}\{-h,h\},
\end{array}
\right.
\end{eqnarray*}
and $\Phi(t)$ and $\Psi(t)$ satisfy the following equations
\begin{eqnarray}\label{SFC-8}
\left\{
\begin{array}{ll}
\Phi'(t)=[\lambda-d_{1}\lambda_{0}-a_{11}]\Phi(t)+a_{12}\Psi(t),\; &\,  t\in(0^{+}, \tau], \\[2mm]
\Psi'(t)=f'(0)\Phi(t)+[\lambda-d_{2}\lambda_{0}-a_{22}]\Psi(t),\; &\,  t\in(0^{+}, \tau], \\[2mm]
(\Phi(0),\Psi(0)) =(\Phi(\tau),\Psi(\tau)),\; &\,\\[2mm]
(\Phi(0^{+}), \Psi(0^{+}))=(G'(0)\Phi(0), \Psi(0)).\; &\,
\end{array} \right.
\end{eqnarray}
It is clear that $\Phi(t)$ and $\Psi(t)$ will depend on the length of the initial region $[-h,h]$. To reflect this dependency, we also write $\Phi(t)=\Phi_{h}(t)$ and $\Psi(t)=\Psi_{h}(t)$. Based on this, some estimates of $\Phi_{h}(t)$ and $\Psi_{h}(t)$ are given next.
\begin{lemma}\label{lemma 3-3}
For any $h\geq h_{0}$, there exist positive constants $\alpha_{1}$, $\alpha_{2}$, $\beta_{1}$, and $\beta_{2}$ that are independent of $h$, such that
\begin{eqnarray*}
\begin{array}{l}
\alpha_{1}\leq \Phi_{h}(0), ~~\Phi_{h}(t)\leq \alpha_{2}, ~~~~t\in[0,\tau], \\[2mm]
\beta_{1}\leq \Psi_{h}(0), ~~\Psi_{h}(t)\leq \beta_{2},~~~~~t\in[0,\tau].
\end{array}
\end{eqnarray*}

\begin{proof}
Motivated by the work \cite{pu-lin-lou}, we first solve for $\Phi(t)$ and $\Psi(t)$. By the first two equations of \eqref{SFC-8}, it follows that
\begin{equation*}
\Bigg(
\begin{matrix}
   \Phi'_{h}(t)\\
   \Psi'_{h}(t) \\
  \end{matrix}
\Bigg)=\Bigg(
\begin{matrix}
   \lambda-d_{1}\lambda_{0}(h)-a_{11}&a_{12}\\
   f'(0)&   \lambda-d_{2}\lambda_{0}(h)-a_{22}\\
  \end{matrix}
\Bigg)
\Bigg(
\begin{matrix}
   \Phi_{h}(t) \\
   \Psi_{h}(t) \\
  \end{matrix}
\Bigg):=\mathbf{A}
\Bigg(
\begin{matrix}
   \Phi_{h}(t) \\
   \Psi_{h}(t) \\
  \end{matrix}
\Bigg).
\end{equation*}
Then, it follows from the characteristic equation $|\mathbf{A}-\kappa E|=0$ that
\begin{equation*}
\kappa_{1,2}=\frac{2\lambda-(d_{1}+d_{2})\lambda_{0}(h)-a_{11}-a_{22}\pm\sqrt{[a_{22}+(d_{2}-d_{1})\lambda_{0}(h)-a_{11}]^{2}+4a_{12}f'(0)}}{2}:=\lambda+c_{1,2}.
\end{equation*}
A simple calculation gives that
\begin{eqnarray}\label{SFC-13}
\begin{array}{l}
2\sqrt{a_{12}f'(0)}\leq c_{1}-c_{2}\leq 2\sqrt{a_{12}f'(0)}+a_{22}+a_{11}+(d_{2}+d_{1})\lambda_{0}(h_{0}), \\[2mm]
0<a_{11}+d_{1}\lambda_{0}(h)+c_{1}=-\big(a_{22}+d_{2}\lambda_{0}(h)+c_{2}\big) \\[2mm]
~~~~~~~~~~~~~~~~~~~~~~~~~~~~~~~~~~~~\leq a_{11}+a_{22}+(d_{2}+d_{1})\lambda_{0}(h_{0})+\sqrt{a_{12}f'(0)}.
\end{array}
\end{eqnarray}
The linearly independent eigenvectors $(x_{11}, x_{12})$ and $(x_{21}, x_{22})$ associated with the eigenvalues $\kappa_{1}$ and $\kappa_{2}$ satisfy
\begin{equation*}
\big(
\begin{matrix}
   x_{i1}&x_{i2}\\
  \end{matrix}
\big)
\Bigg(
\begin{matrix}
   \lambda-d_{1}\lambda_{0}(h)-a_{11}-\kappa_{i}&a_{12}\\
   f'(0)&   \lambda-d_{2}\lambda_{0}(h)-a_{22}-\kappa_{i}\\
  \end{matrix}
\Bigg)=
\big(
\begin{matrix}
   0&0\\
  \end{matrix}
\big)
\end{equation*}
for $i=1,2$, and a simple calculation yields that
\begin{equation*}
(x_{11}, x_{12})=\big(f'(0), a_{11}+d_{1}\lambda_{0}(h)-\lambda+\kappa_{1}\big)=\big(f'(0), a_{11}+d_{1}\lambda_{0}(h)+c_{1}\big)
\end{equation*}
and
\begin{equation*}
(x_{21}, x_{22})=\big(a_{22}+d_{2}\lambda_{0}(h)-\lambda+\kappa_{2}, a_{12}\big)=\big(a_{22}+d_{2}\lambda_{0}(h)+c_{2}, a_{12}\big).
\end{equation*}

Next, we consider the algebraic equations
\begin{equation*}
\Bigg(
\begin{matrix}
   e^{\kappa_{1}t} \\
   ke^{\kappa_{2}t} \\
  \end{matrix}
\Bigg)=
\Bigg(
\begin{matrix}
   x_{11}&x_{12}\\
   x_{21}&x_{22}\\
  \end{matrix}
\Bigg)
\Bigg(
\begin{matrix}
   \Phi(t)\\
   \Psi(t) \\
  \end{matrix}
\Bigg):=\mathbf{B}
\Bigg(
\begin{matrix}
   \Phi(t)\\
   \Psi(t) \\
  \end{matrix}
\Bigg),
\end{equation*}
and a simple calculation yields that
\begin{equation}\label{SFC-12}
(\Phi(t), \Psi(t))=\Bigg(\frac{a_{12}e^{\kappa_{1}t}-\big(a_{11}+d_{1}\lambda_{0}(h)+c_{1}\big)ke^{\kappa_{2}t}}{|\mathbf{B}|},
\frac{f'(0)ke^{\kappa_{2}t}-\big(a_{22}+d_{2}\lambda_{0}(h)+c_{2}\big)e^{\kappa_{1}t}}{|\mathbf{B}|}\Bigg),
\end{equation}
where
\begin{equation*}
\begin{aligned}
|\mathbf{B}|=&a_{12}f'(0)-\big(a_{11}+d_{1}\lambda_{0}(h)+c_{1}\big)\big(a_{22}+d_{2}\lambda_{0}(h)+c_{2}\big)\\
=&a_{12}f'(0)+\big(a_{22}+d_{2}\lambda_{0}(h)+c_{2}\big)^{2}>0.
\end{aligned}
\end{equation*}
By applying the last three conditions of \eqref{SFC-8}, one can obtain that
\begin{equation}\label{SFC-9}
\begin{cases}
\begin{aligned}
&~~~a_{12}-\big(a_{11}+d_{1}\lambda_{0}(h)+c_{1}\big)k=G'(0)\big[a_{12}e^{\kappa_{1}\tau}-\big(a_{11}+d_{1}\lambda_{0}(h)+c_{1}\big)ke^{\kappa_{2}\tau}\big],\\
&f'(0)k-\big(a_{22}+d_{2}\lambda_{0}(h)+c_{2}\big)=f'(0)ke^{\kappa_{2}\tau}-\big(a_{22}+d_{2}\lambda_{0}(h)+c_{2}\big)e^{\kappa_{1}\tau}.
\end{aligned}
\end{cases}
\end{equation}
For simplicity, we denote
\begin{equation*}
\begin{aligned}
&~~~~~~~~~~y=e^{\kappa_{1}\tau},~ n_{11}=a_{12},  ~n_{12}=a_{11}+d_{1}\lambda_{0}(h)+c_{1},~n_{13}=G'(0)a_{12},\\
&~n_{21}=f'(0),~n_{22}=f'(0)e^{(\kappa_{2}-\kappa_{1})\tau}, ~n_{23}=G'(0)\big(a_{11}+d_{1}\lambda_{0}(h)+c_{1}\big)e^{(\kappa_{2}-\kappa_{1})\tau},
\end{aligned}
\end{equation*}
and then \eqref{SFC-9} can be rewritten as
\begin{equation}\label{SFC-10-1}
\begin{cases}
\begin{aligned}
&n_{11}-n_{12}k=(n_{13}-n_{23}k)y,\\
&n_{12}+n_{21}k=(n_{12}+n_{22}k)y.
\end{aligned}
\end{cases}
\end{equation}
Since $\Phi_{h}(t)$ is less than or equal to zero on $[0, \tau]$ when $k=n_{13}/n_{23}$ and $\Psi_{h}(t)$ is less than or equal to zero on $[0, \tau]$ when $k=-(n_{12}/n_{22})$, \eqref{SFC-10-1} can be  written as
\begin{equation}\label{SFC-10}
\begin{cases}
\begin{aligned}
&y=(n_{11}-n_{12}k)/(n_{13}-n_{23}k),\\
&y=(n_{12}+n_{21}k)/(n_{12}+n_{22}k).
\end{aligned}
\end{cases}
\end{equation}
By considering the image of \eqref{SFC-10}, we search for solutions to \eqref{SFC-10}.
\begin{figure}[!httb]
  \centering
  \includegraphics[scale=0.75]{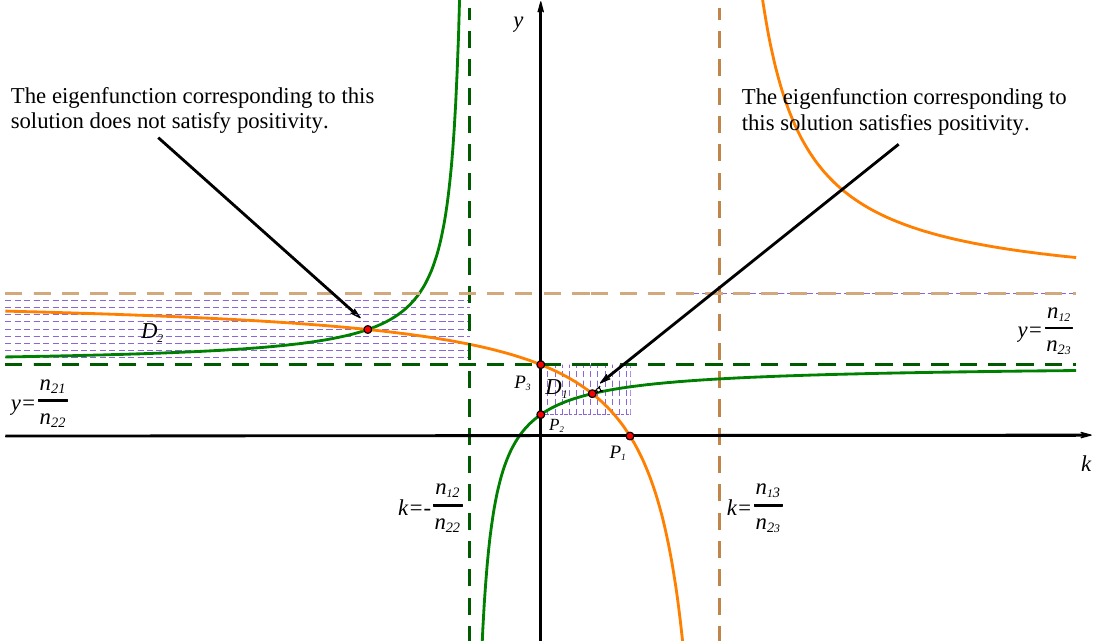}
    \caption{Distribution graph of solutions to \eqref{SFC-10}. }\label{figure 1}
\end{figure}
Let a solution where the corresponding eigenfunctions satisfy the positivity condition be $(k_{0}, y_{0})$. Then, it follows that
\begin{equation}\label{SFC-11}
\begin{aligned}
0<k_{0}<n_{11}/n_{12}~\text{and}~1<y_{0}<n_{21}/n_{22}.
\end{aligned}
\end{equation}
With the help of \eqref{SFC-13}, \eqref{SFC-12}, and \eqref{SFC-11}, one can obtain that
\begin{equation*}
\begin{aligned}
\Phi(t)&=\frac{e^{\kappa_{1}t}\big[a_{12}-\big(a_{11}+d_{1}\lambda_{0}(h)+c_{1}\big)k_{0}e^{(\kappa_{2}-\kappa_{1})t}\big]}{a_{12}f'(0)
+\big(a_{11}+d_{1}\lambda_{0}(h)+c_{1}\big)^{2}}\\
&\leq \frac{e^{\kappa_{1}t}}{f'(0)}\leq \frac{e^{(c_{1}-c_{2})t}}{f'(0)}\leq \frac{e^{\big[2\sqrt{a_{12}f'(0)}+a_{22}+a_{11}+(d_{2}+d_{1})\lambda_{0}(h_{0})\big]\tau}}{f'(0)}:=\alpha_{2},
\end{aligned}
\end{equation*}
\begin{equation*}
\begin{aligned}
\Psi(t)&=\frac{f'(0)k_{0}e^{\kappa_{2}t}+\big(a_{11}+d_{1}\lambda_{0}(h)+c_{1}\big)e^{\kappa_{1}t}}{a_{12}f'(0)+\big(a_{11}+d_{1}\lambda_{0}(h)+c_{1}\big)^{2}}\\
&\leq\frac{a_{12}f'(0)e^{\kappa_{2}t}+\big(a_{11}+d_{1}\lambda_{0}(h)+c_{1}\big)^{2}e^{\kappa_{1}t}}{\big[a_{11}+d_{1}\lambda_{0}(h)+c_{1}\big]
\big[a_{12}f'(0)+\big(a_{11}+d_{1}\lambda_{0}(h)+c_{1}\big)^{2}\big]}\\
&\leq \frac{e^{\kappa_{1}t}}{a_{11}+d_{1}\lambda_{0}(h)+c_{1}}\leq \frac{f'(0)\alpha_{2}}{a_{11}} :=\beta_{2},
\end{aligned}
\end{equation*}
and
\begin{equation*}
\begin{aligned}
\Psi(0)&=\frac{f'(0)k_{0}+\big(a_{11}+d_{1}\lambda_{0}(h)+c_{1}\big)}{a_{12}f'(0)+\big(a_{11}+d_{1}\lambda_{0}(h)+c_{1}\big)^{2}}\\
&\geq\frac{a_{11}+d_{1}\lambda_{0}(h)+c_{1}}{a_{12}f'(0)+\big(a_{11}+d_{1}\lambda_{0}(h)+c_{1}\big)^{2}}\\
&\geq\frac{a_{11}}{a_{12}f'(0)+\big(a_{11}+a_{22}+(d_{2}+d_{1})\lambda_{0}(h_{0})+\sqrt{a_{12}f'(0)}\big)^{2}} :=\beta_{1}.
\end{aligned}
\end{equation*}

Next, we prove that there exists a positive constant
\begin{equation*}
\epsilon_{0}=\min\bigg\{\frac{a_{12}}{2}, a_{12}G'(0)\Big(1-e^{-2\sqrt{a_{12}f'(0)}\tau}\Big)\bigg\}
\end{equation*}
such that
\begin{equation}\label{SFC-14}
0<k_{0}<(n_{11}-\epsilon_{0})/n_{12}:=k_{1}.
\end{equation}
By standard calculations it follows that
\begin{equation*}
\begin{aligned}
&\frac{n_{12}+n_{21}k_{1}}{n_{12}+n_{22}k_{1}}-\frac{n_{11}-n_{12}k_{1}}{n_{13}-n_{23}k_{1}}=\frac{n^{2}_{12}+n_{21}(n_{11}-\epsilon_{0})}{n^{2}_{12}
+n_{22}(n_{11}-\epsilon_{0})}-\frac{n_{12}\epsilon_{0}}{n_{12}n_{13}-n_{23}(n_{11}-\epsilon_{0})}\\
>&\frac{n^{2}_{12}+n_{21}(n_{11}-\frac{n_{11}}{2})}{n^{2}_{12}+n_{22}(n_{11}-\frac{n_{11}}{2})}-\frac{n_{12}\epsilon_{0}}{n_{12}n_{13}-n_{23}n_{11}}
>1-\frac{\epsilon_{0}}{a_{12}G'(0)[1-e^{(\kappa_{2}-\kappa_{1})\tau}]} \geq 0.
\end{aligned}
\end{equation*}
This implies that \eqref{SFC-14} is valid. Due to \eqref{SFC-14}, one can now get that
\begin{equation*}
\begin{aligned}
\Phi(0)&=\frac{a_{12}-\big(a_{11}+d_{1}\lambda_{0}(h)+c_{1}\big)k_{0}}{a_{12}f'(0)+\big(a_{11}+d_{1}\lambda_{0}(h)+c_{1}\big)^{2}}\geq \frac{\epsilon_{0}}{a_{12}f'(0)+\big(a_{11}+d_{1}\lambda_{0}(h)+c_{1}\big)^{2}}\\
&\geq \frac{\epsilon_{0}}{a_{12}f'(0)+\big(a_{11}+a_{22}+(d_{2}+d_{1})\lambda_{0}(h_{0})+\sqrt{a_{12}f'(0)}\big)^{2}}:=\alpha_{1}.
\end{aligned}
\end{equation*}
This ends the proof.
\end{proof}

\end{lemma}the threshold that governs the extinction and persistence of an infectious disease

The principal eigenvalue is usually the threshold that governs the extinction and persistence of an infectious disease. Hence, it plays a significant role in the research of epidemics. For the  impulse faecal-oral model \eqref{3-1} in a fixed region, we define
\begin{equation*}
\lambda^{\vartriangle}_{1}(G'(0),(-\infty,+\infty))=\lim\limits_{-s_{1}, s_{2}\rightarrow\infty}\lambda^{\vartriangle}_{1}(G'(0),(s_{1},s_{2})).
\end{equation*}
In addition, we also write $\lambda^{\vartriangle}_{1}(\infty)=\lambda^{\vartriangle}_{1}(G'(0),(-\infty,+\infty))$ and $\lambda^{\vartriangle}_{1}(h_{0})=\lambda^{\vartriangle}_{1}(G'(0),(-h_{0},h_{0}))$.
For the corresponding model \eqref{Zhou-Lin} with free boundary, the principal eigenvalue is defined by
\begin{equation*}
\lambda^{\triangledown}_{1}(t):=\lambda^{\vartriangle}_{1}(G'(0),(g(t),h(t))).
\end{equation*}
It virtue of \autoref{lemma 2-2} and \autoref{lemma 3-1}\textcolor{blue}{(1,3)}, it follows that $\lambda^{\triangledown}_{1}(t)$ is strongly decreasing with respect to $t$. Therefore,
\begin{equation*}
\lim\limits_{t\rightarrow\infty}\lambda^{\triangledown}_{1}(t)=\lambda^{\triangledown}_{1}(\infty)\in(-\infty, \lambda^{\triangledown}_{1}(0))
\end{equation*}
is always well-defined. In order to emphasize the dependence of $\lambda^{\triangledown}_{1}$ on $G'(0)$, $\mu_{1}$, and $\mu_{2}$, we also write $\lambda^{\triangledown}_{1}=\lambda^{\triangledown}_{1}(G'(0), \mu_{1}, \mu_{2})$. Next, we demonstrate that $\lambda^{\triangledown}_{1}(t)$ is an important index of model \eqref{Zhou-Lin}.
\begin{theorem}\label{theorem 3-1}
If $\lambda^{\vartriangle}_{1}(G'(0),(-\infty,+\infty))>0$, then the infected individuals and pathogenic bacteria are vanishing.
\begin{proof}
To begin with, construct the following functions
\begin{eqnarray}\label{3-4}
\begin{array}{l}
~~~~~~~~~~~~ ~~~~~\Gamma(t):=h\kappa(t),~\kappa(t):=1+\sigma-\frac{\sigma}{2}e^{-\gamma t}, ~~~~t\geq0, \\[2mm]
\overline{u}(t,x):=Me^{-\gamma t}\phi(t, \frac{hx}{\Gamma(t)}), \overline{v}(t,x):=Me^{-\gamma t}\psi(t,\frac{hx}{\Gamma(t)}),t\geq 0, x\in[-\Gamma(t), \Gamma(t)],
\end{array}
\end{eqnarray}
where the values of the positive constants $h(>h_{0})$, $\gamma$, $\sigma$, and $M$ are chosen subsequently, and $(\phi, \psi)$ is a strongly positive principal eigenfunction pair for the periodic eigenvalue problem \eqref{3-2} with $(s_{1}, s_{2})=(-h, h)$. From the method of separation of variables, it follows that
$(\phi(t,x),\psi(t,x))=X_{h}(x)(\Phi_{h}(t),\\\Psi_{h}(t))$, where $\big(\lambda_{0}, X_{h}(x)\big)$ is the principal eigenpair of the eigenvalue problem
\begin{eqnarray*}
\left\{
\begin{array}{ll}
\Delta X_{h}(x)=-\lambda X_{h}(x),\; &\, x\in (-h, h), \\[2mm]
X_{h}(x)=0,\; &\, x\in \{-h, h\}.
\end{array}
\right.
\end{eqnarray*}
Simple calculations yields that $\lambda_{0}=\frac{\pi^{2}}{4h^{2}}$ and $X_{h}(x)=\cos(\frac{\pi}{2h}x)$ for $x\in[-h, h]$.

Next, we proceed to show that the quadruple $(\overline{u},\overline{v},-\Gamma(t),\Gamma(t))$ is a supersolution of model \eqref{Zhou-Lin}.  With the help of \autoref{lemma 3-3}, the conclusion that $\lambda^{\vartriangle}_{1}(G'(0),(-h, h))$ is strongly decreasing w.r.t. $h$, and the assumption that $\frac{f(u)}{u}$ is strongly decreasing for $u>0$ in (A2), it follows that
\begin{equation}\label{3-5}
\begin{aligned}
&\overline{u}_{t}-d_{1}\Delta \overline{u}+a_{11}\overline{u}-a_{12}\overline{v}\\
\geq&~Me^{-\gamma t}\big[-\gamma\phi+\phi_{t}-hx\Gamma^{-2}(t)\Gamma_{t}(t)\phi_{x}-d_{1}\kappa^{-2}(t)\Delta\phi+a_{11}\phi-a_{12}\psi\big] \\
\geq&~Me^{-\gamma t}\big[\lambda^{\vartriangle}_{1}(G'(0),(-h,h))\phi+d_{1}\Delta\phi-\gamma\phi-d_{1}\kappa^{-2}\Delta\phi\big] \\
\geq&~Me^{-\gamma t}\phi\big[\lambda^{\vartriangle}_{1}(G'(0),(-\infty,+\infty))-\gamma-d_{1}\lambda_{0}(1-\kappa^{-2})\big]\\
\geq&~Me^{-\gamma t}\phi\big[\lambda^{\vartriangle}_{1}(G'(0),(-\infty,+\infty))-\gamma-d_{1}\lambda_{0}(h_{0})(1-\kappa^{-2})\big]
\end{aligned}
\end{equation}
and
\begin{equation}\label{3-6}
\begin{aligned}
&\overline{v}_{t}-d_{2}\Delta \overline{v}+a_{22}\overline{v}-f(\overline{u})\\
\geq&~Me^{-\gamma t}\big[-\gamma\psi+\psi_{t}-hx\Gamma^{-2}(t)\Gamma_{t}(t)\psi_{x}-d_{2}\kappa^{-2}(t)\Delta\psi+a_{22}\psi-f'(0)\phi\big] \\
\geq&~Me^{-\gamma t}\big[\lambda^{\vartriangle}_{1}(G'(0),(-h,h))\psi+d_{2}\Delta\psi-\gamma\phi-d_{2}\kappa^{-2}\Delta\psi\big] \\
\geq&~Me^{-\gamma t}\psi\big[\lambda^{\vartriangle}_{1}(G'(0),(-\infty,+\infty))-\gamma-d_{2}\lambda_{0}(1-\kappa^{-2})\big]\\
\geq&~Me^{-\gamma t}\psi\big[\lambda^{\vartriangle}_{1}(G'(0),(-\infty,+\infty))-\gamma-d_{2}\lambda_{0}(h_{0})(1-\kappa^{-2})\big]
\end{aligned}
\end{equation}
for $t\in((k\tau)^{+}, (k+1)\tau]$ and $x\in(-\Gamma(t), \Gamma(t))$. Select
\begin{equation*}
\gamma:=\frac{\lambda^{\vartriangle}_{1}(G'(0),(-\infty,+\infty))}{2}.
\end{equation*}
Noting that $\lim\limits_{\sigma\rightarrow 0^{+}}(1-\kappa^{-2})=0$ and $1-\kappa^{-2}>0$ for $\sigma>0$, the positive constant $\sigma$ can be chosen such that
\begin{equation*}
1-\kappa^{-2}<\frac{\lambda^{\vartriangle}_{1}(G'(0),(-\infty,+\infty))}{2\lambda_{0}(h_{0})\max\big\{d_{1}, d_{2}\big\}}.
\end{equation*}
Therefore, it follows from \eqref{3-5} and \eqref{3-6} that $\overline{u}_{t}-d_{1}\Delta \overline{u}+a_{11}\overline{u}-a_{12}\overline{v}\geq 0$ and $\overline{v}_{t}-d_{2}\Delta \overline{v}+a_{22}\overline{v}-f(\overline{u})\geq 0$ for $(t, x)\in ((k\tau)^{+}, (k+1)\tau]\times (-\Gamma(t), \Gamma(t))$.

With the help of Assumption (A2), one can obtain that
\begin{equation*}
\begin{aligned}
&\overline{u}((k\tau)^{+},x)\geq Me^{-\gamma k\tau}\phi((k\tau)^{+}, \frac{hx}{\Gamma(k\tau)})\\
\geq&~G'(0)Me^{-\gamma k\tau}\phi(k\tau, \frac{hx}{\Gamma(k\tau)}) \geq \overline{u}(k\tau,x)G'(0)\\
\geq&~G(\overline{u}(k\tau,x))
\end{aligned}
\end{equation*}
and
\begin{equation*}
\begin{aligned}
&\overline{v}((k\tau)^{+},x)\geq Me^{-\gamma k\tau}\psi((k\tau)^{+}, \frac{hx}{\Gamma(k\tau)})\\
\geq&~Me^{-\gamma k\tau}\phi(k\tau, \frac{hx}{\Gamma(k\tau)}) \geq \overline{v}(k\tau,x)
\end{aligned}
\end{equation*}
for $x\in(-\Gamma(k\tau), \Gamma(k\tau))$. For the boundary conditions, one can obtain that
\begin{equation*}
\overline{u}(t, \pm\Gamma(t) )=\overline{v}(t, \pm\Gamma(t) )=0\text{~for~}t>0.
\end{equation*}
Noticing that $X_{h}(x)=\cos(\frac{\pi}{2h}x)$ for $x\in[-h, h]$, the sufficiently large $M$ can be chosen such that
\begin{equation*}
\overline{u}(0, x)=M\Phi_{h}(0)X_{h}(\frac{2x}{2+\sigma})\geq M\alpha_{1}X_{h_{0}}(\frac{2x}{2+\sigma})\geq u_{0}(x)
\end{equation*}
and
\begin{equation*}
\overline{v}(0, x)=M\Psi_{h}(0)X_{h}(\frac{2x}{2+\sigma})\geq M\beta_{1}X_{h_{0}}(\frac{2x}{2+\sigma})\geq v_{0}(x)
\end{equation*}
for $x\in [-h_{0}, h_{0}]$, where $\vartheta_{1}(\vartheta=\alpha,\beta)$ is defined in \autoref{lemma 3-3}. Finally, choose
\begin{equation*}
h:=\max\Bigg\{h_{0}, \sqrt{\frac{2M\pi}{\gamma\sigma(2+\sigma)}(\mu_{1}\alpha_{2}+\mu_{2}\beta_{2})}\Bigg\},
\end{equation*}
where $\alpha_{2}$ and $\beta_{2}$ are defined in \autoref{lemma 3-3}. Further calculations now give that
\begin{equation*}
\begin{aligned}
&h[\Gamma'(t)+\mu_{1}\overline{u}_{x}(t, \Gamma(t))+\mu_{2}\overline{v}_{x}(t, \Gamma(t))]\\
\geq&~h\Big[h\frac{\sigma}{2}\gamma e^{-\gamma t}+\mu_{1}Me^{-\gamma t}\phi_{x}(t, h)\kappa^{-1}(t)+\mu_{2}Me^{-\gamma t}\psi_{x}(t, h)\kappa^{-1}(t)\Big]\\
\geq&~e^{-\gamma t}\Big[\frac{\sigma\gamma h^{2}}{2} -\frac{\mu_{1}\pi M}{2}\alpha_{2}\kappa^{-1}(t)-\frac{\mu_{2}\pi M}{2}\beta_{2}\kappa^{-1}(t)\Big]\\
\geq&~e^{-\gamma t}\Big[\frac{\sigma\gamma h^{2}}{2} -\frac{\pi M}{2+\sigma}(\mu_{1}\alpha_{2}+\mu_{2}\beta_{2})\Big]\geq 0
\end{aligned}
\end{equation*}
and
\begin{equation*}
\begin{aligned}
&h[-\Gamma'(t)+\mu_{1}\overline{u}_{x}(t,-\Gamma(t))+\mu_{2}\overline{v}_{x}(t, -\Gamma(t))]\\
\leq&~h\Big[-h\frac{\sigma}{2}\gamma e^{-\gamma t}+\mu_{1}Me^{-\gamma t}\phi_{x}(t, -h)\kappa^{-1}(t)+\mu_{2}Me^{-\gamma t}\psi_{x}(t, -h)\kappa^{-1}(t)\Big]\\
\leq&~e^{-\gamma t}\Big[-\frac{\sigma\gamma h^{2}}{2} +\frac{\mu_{1}\pi M}{2}\alpha_{2}\kappa^{-1}(t)+\frac{\mu_{2}\pi M}{2}\beta_{2}\kappa^{-1}(t)\Big]\\
\leq&~e^{-\gamma t}\Big[-\frac{\sigma\gamma h^{2}}{2} +\frac{\pi M}{2+\sigma}(\mu_{1}\alpha_{2}+\mu_{2}\beta_{2})\Big]\leq 0
\end{aligned}
\end{equation*}
for $t>0$. We can now complete the proof that the quadruple $(\overline{u},\overline{v},-\Gamma(t),\Gamma(t))$ is a supersolution of model \eqref{Zhou-Lin}.

Since the quadruple $(\overline{u},\overline{v},-\Gamma(t),\Gamma(t))$ is the supersolution, it follows from \autoref{lemma 3-2} that
\begin{eqnarray}\label{3-7}
\left\{
\begin{array}{l}
g(t)\geq-\Gamma(t), ~~~~~\Gamma(t)\geq h(t) ,~~~~~~~t\in\mathbb{R}^{+}, \\[2mm]
\overline{u}(t, x)\geq u(t,x), \overline{v}(t, x)\geq v(t,x),(x,t)\in\mathbb{R}^{+}\times[g(t), h(t)].
\end{array} \right.
\end{eqnarray}
In addition, it is easy to obtain from \eqref{3-4} that
\begin{equation*}
\lim\limits_{t\rightarrow+\infty}\Gamma(t)=h(1+\sigma)<+\infty
\end{equation*}
and
\begin{equation*}
\lim\limits_{t\rightarrow\infty}\Big[\|\overline{u}(t, \cdot)\|_{\mathbb{C}[g(t), h(t)]}+\|\overline{v}(t, \cdot)\|_{\mathbb{C}[g(t), h(t)]}\Big]=0.
\end{equation*}
This combined with \eqref{3-7} yields that
\begin{equation*}
h_{\infty}-g_{\infty}<2h(1+\sigma)<+\infty
\end{equation*}
and
\begin{equation*}
\lim\limits_{t\rightarrow\infty}\Big[\|u(t, \cdot)\|_{\mathbb{C}[g(t), h(t)]}+\|v(t, \cdot)\|_{\mathbb{C}[g(t), h(t)]}\Big]=0.
\end{equation*}
This means that the infected individuals and pathogenic bacteria are vanishing. This proof is now completed.
\end{proof}
\end{theorem}

\autoref{theorem 3-1} presents the long time dynamical behaviour of the impulsive faecal-oral model \eqref{Zhou-Lin} when $\lambda^{\vartriangle}_{1}(G'(0),(-\infty,+\infty))>0$. Upcoming \autoref{lemma 3-4} deals with the case when $\lambda^{\vartriangle}_{1}(G'(0),(-\infty,+\infty))=0$.
\begin{lemma}\label{lemma 3-4}
If $\lambda^{\vartriangle}_{1}(G'(0),(-\infty,+\infty))=0$, then the infected individuals and pathogenic bacteria are extinct, that is,
\begin{equation*}
\lim\limits_{t\rightarrow+\infty}\|u(t, \cdot)\|_{\mathbb{C}[g(t), h(t)]}=\lim\limits_{t\rightarrow+\infty}\|v(t, \cdot)\|_{\mathbb{C}[g(t), h(t)]}=0.
\end{equation*}
\begin{proof}
Since the proof is similar to \autoref{lemma 3-1}\textcolor{blue}{$(4)$}, only sketches are given here. For more details, see the proof of  \autoref{lemma 3-1}\textcolor{blue}{$(4)$}.

By using the condition $\lambda^{\vartriangle}_{1}(G'(0),(-\infty,+\infty))=0$, one can obtain that
\begin{eqnarray*}
\left\{
\begin{array}{ll}
U_{t}=a_{12}V-a_{11}U, ~~~~~~t\in(0^{+}, \tau], \\[2mm]
V_{t}=f(U)-a_{22}V, ~~~~~~t\in(0^{+}, \tau], \\[2mm]
U(0)=U(\tau), V(0)=V(\tau),\\[2mm]
U(0^{+})=G(U(0)), V(0^{+})=V(0)
\end{array} \right.
\end{eqnarray*}
has a unique nonnegative periodic solution $(0,0)$. Denote the solution to the problem
\begin{eqnarray*}
\left\{
\begin{array}{ll}
\overline{u}_{t}=a_{12}\overline{v}-a_{11}\overline{u},~~~~~~~~~~~~~ t\in((k\tau)^{+}, (k+1)\tau], \\[2mm]
\overline{v}_{t}=f(\overline{u})-a_{22}\overline{v},~~~~~~~~~~~~~ t\in((k\tau)^{+}, (k+1)\tau], \\[2mm]
\overline{u}((k\tau)^{+})=G(\overline{u}(k\tau)), \overline{v}((k\tau)^{+})=\overline{v}(k\tau), \\[2mm]
\overline{u}(0)=\|u_{0}(\cdot)\|_{\mathbb{C}[-h_{0}, h_{0}]}, \overline{v}(0)=\|v_{0}(\cdot)\|_{\mathbb{C}[-h_{0}, h_{0}]}
\end{array} \right.
\end{eqnarray*}
by $(\overline{u}, \overline{v})$. With the help of \autoref{lemma 3-2}, it follows that $(u,v)\leq(\overline{u}(t),\overline{v}(t))$ for $t\geq 0$ and $x\in [g(t), h(t)]$. Next, the upper and lower solutions method implies that
\begin{equation*}
\lim\limits_{t\rightarrow+\infty}[\overline{u}(t)+\overline{v}(t)]=0.
\end{equation*}
This means that
\begin{equation*}
\lim\limits_{t\rightarrow+\infty}\big[\|u(t, \cdot)\|_{\mathbb{C}([g(t), h(t)])}+\|v(t, \cdot)\|_{\mathbb{C}([g(t), h(t)])}\big]=0.
\end{equation*}
This ends the proof.
\end{proof}
\end{lemma}

Next, we discuss the asymptotic stability of model \eqref{Zhou-Lin} when $\lambda^{\vartriangle}_{1}(G'(0),(-\infty,+\infty))<0$. Before discussing the dynamic behavior, it is necessary to specify how the fronts $g(t)$ and $h(t)$ move.
\begin{lemma}\label{lemma 3-5}
If $\lambda^{\vartriangle}_{1}(G'(0),(-\infty,+\infty))<0$, and $-\infty<g_{\infty}$ or $h_{\infty}<+\infty$,  then $h_{\infty}-g_{\infty}<+\infty$.
\begin{proof}
We only need to prove that if $\lambda^{\vartriangle}_{1}(G'(0),(-\infty,+\infty))<0$ and $h_{\infty}<+\infty$,  then $h_{\infty}-g_{\infty}<+\infty$.
This is because a similar method can prove the other case. Based on the proof of contradiction, suppose that $h_{\infty}-g_{\infty}=+\infty$. Since $\lambda^{\vartriangle}_{1}(G'(0),(-\infty,+\infty))<0$, it follows from \autoref{lemma 3-1}\textcolor{blue}{$(1, 3)$} that there exists
a sufficiently large positive integer $k_{0}$ such that $\lambda^{\triangledown}_{1}(k_{0}\tau)<0$ and $\lambda^{\triangledown}_{1}(t)<\lambda^{\triangledown}_{1}(k_{0}\tau)<0$ for $t>k_{0}\tau$. For the sake of readability, the next proof is
divided into three parts.

\textbf{(1) The associated eigenvalue problem}

Consider the eigenvalue problem
\begin{eqnarray}\label{3-8}
\left\{
\begin{array}{ll}
\phi_{t}=d_{1}\Delta \phi+\epsilon\phi_{x}-a_{11}\phi+a_{12}\psi+\lambda \phi,\; &\, (t,x)\in(0^{+}, \tau]\times(g(k_{0}\tau), h(k_{0}\tau)), \\[2mm]
\psi_{t}=d_{2}\Delta \psi+\frac{d_{2}\epsilon}{d_{1}}\psi_{x}-a_{22}\psi+f'(0)\phi+\lambda \psi,\; &\, (t,x)\in(0^{+}, \tau]\times(g(k_{0}\tau), h(k_{0}\tau)), \\[2mm]
(\phi(t,x), \psi(t,x))=(0,0),\; &\, (t,x)\in[0, \tau]\times\{g(k_{0}\tau), h(k_{0}\tau)\},\\[2mm]
(\phi(0,x), \psi(0,x))=(\phi(\tau,x),\psi(\tau,x)),\; &\,x\in[g(k_{0}\tau), h(k_{0}\tau)],\\[2mm]
(\phi(0^{+},x), \psi(0^{+},x))=(G'(0)\phi(0,x), \psi(0,x)), \; &\,  x\in(g(k_{0}\tau), h(k_{0}\tau)).
\end{array} \right.
\end{eqnarray}
For a sufficiently small $\epsilon$, problem \eqref{3-8} has a principal eigenvalue $\lambda^{\epsilon}_{1}$ with $\lambda^{\epsilon}_{1}<0$,
and the corresponding eigenfunction pair $(\phi^{\epsilon}, \psi^{\epsilon})$ can be represented as
\begin{equation*}
\begin{aligned}
&\phi^{\epsilon}=\Phi(t)e^{-\frac{\epsilon}{2d_{1}}x}\cos\Big[\frac{\pi}{h(k_{0}\tau)-g(k_{0}\tau)}\Big(x-\frac{h(k_{0}\tau)+g(k_{0}\tau)}{2}\Big)\Big],  \\
&\psi^{\epsilon}=\Psi(t)e^{-\frac{\epsilon}{2d_{1}}x}\cos\Big[\frac{\pi}{h(k_{0}\tau)-g(k_{0}\tau)}\Big(x-\frac{h(k_{0}\tau)+g(k_{0}\tau)}{2}\Big)\Big],
\end{aligned}
\end{equation*}
where $\Phi(t)$ and $\Psi(t)$ satisfy \eqref{SFC-8} with $\lambda_{0}(h_{0})$ replaced by $\frac{\epsilon^{2}}{4d^{2}_{1}}+\big(\frac{\pi}{h(k_{0}\tau)-g(k_{0}\tau)}\big)^{2}$. For the sake of the later proof, we stipulate that
\begin{equation}\label{3-9}
\begin{aligned}
&\phi^{\epsilon}_{x}=-\frac{\epsilon}{2d_{1}}\phi^{\epsilon}-\frac{\pi}{h(k_{0}\tau)-g(k_{0}\tau)}e^{-\frac{\epsilon}{2d_{1}}x}\sin\Big[\frac{\pi}{h(k_{0}\tau)
-g(k_{0}\tau)}\Big(x-\frac{h(k_{0}\tau)+g(k_{0}\tau)}{2}\Big)\Big]\Phi(t),  \\
&\psi^{\epsilon}_{x}=-\frac{\epsilon}{2d_{1}}\psi^{\epsilon}-\frac{\pi}{h(k_{0}\tau)-g(k_{0}\tau)}e^{-\frac{\epsilon}{2d_{1}}x}\sin\Big[\frac{\pi}{h(k_{0}\tau)
-g(k_{0}\tau)}\Big(x-\frac{h(k_{0}\tau)+g(k_{0}\tau)}{2}\Big)\Big]\Psi(t),
\end{aligned}
\end{equation}
and
\begin{equation}\label{3-10}
\begin{aligned}
&\phi^{\epsilon}_{xx}=-\frac{\epsilon}{d_{1}}\phi^{\epsilon}_{x}-\Big[\frac{\epsilon^{2}}{4d^{2}_{1}}+\Big(\frac{\pi}{h(k_{0}\tau)
-g(k_{0}\tau)}\Big)^{2}\Big]\phi^{\epsilon},  \\
&\phi^{\epsilon}_{xx}=-\frac{\epsilon}{d_{1}}\psi^{\epsilon}_{x}-\Big[\frac{\epsilon^{2}}{4d^{2}_{1}}+\Big(\frac{\pi}{h(k_{0}\tau)
-g(k_{0}\tau)}\Big)^{2}\Big]\psi^{\epsilon}.
\end{aligned}
\end{equation}
With the help of \eqref{3-9}, it follows that a constant $x_{0}$ can be found such that
\begin{equation*}
\begin{aligned}
\phi^{\epsilon}_{x}(t,x_{0} )=\psi^{\epsilon}_{x}(t,x_{0} )=0,
\end{aligned}
\end{equation*}
\begin{equation*}
\begin{aligned}
\phi^{\epsilon}_{x},~ \psi^{\epsilon}_{x}>0 \text{~for~} t>0 \text{~and~} x\in[g(k_{0}\tau), x_{0}),
\end{aligned}
\end{equation*}
and
\begin{equation*}
\begin{aligned}
\phi^{\epsilon}_{x}, ~\psi^{\epsilon}_{x}<0 \text{~for~} t>0 \text{~and~} x\in(x_{0}, h(k_{0}\tau)].
\end{aligned}
\end{equation*}

\textbf{(2) A lower solution of $\bm{(u,v)}$ for $\bm{t\geq k_{0}\tau}$ and $\bm{x\in[g(k_{0}\tau), h(k_{0}\tau)]}$}

Here $\epsilon=0$, the principal eigenvalue of problem \eqref{3-8} is denoted by $\lambda^{\triangledown}_{1}$, and the corresponding eigenfunction is denoted by $(\phi, \psi)$. With the help of the method of separation of variables, a sufficiently large $C$ can be found such that
\begin{equation}\label{3-13}
\begin{aligned}
\frac{1}{C}\leq\frac{\psi}{\phi}\leq C\text{~for~} t\geq 0\text{~and~} x\in\big[g(k_{0}\tau), h(k_{0}\tau)\big].
\end{aligned}
\end{equation}
Since $\lambda^{\triangledown}_{1}<0$, there exists a sufficiently small positive number $\alpha$ such that
\begin{equation}\label{3-14}
\begin{aligned}
\lambda^{\triangledown}_{1}+\alpha<0\text{~and~}a_{12}C\big(1-e^{(\lambda^{\triangledown}_{1}+\alpha)\tau}\big)-\alpha<0.
\end{aligned}
\end{equation}
Now we construct the following functions
\begin{eqnarray*}
\underline{v}(t,x)=
\epsilon_{0} e^{(\lambda^{\triangledown}_{1}+\alpha)\tau}\psi(t,x),\; &\, t\geq 0, ~x\in[g(k_{0}\tau), h(k_{0}\tau)],
\end{eqnarray*}
\begin{eqnarray*}
\underline{u}(t,x)=
\left\{
\begin{array}{ll}
\epsilon_{0} \phi(t,x),\; &\, t=k\tau,~ x\in[g(k_{0}\tau), h(k_{0}\tau)], \\[2mm]
\epsilon_{0} e^{(\lambda^{\triangledown}_{1}+\alpha)\tau} \phi(t,x),\; &\, t=(k\tau)^{+}, ~x\in[g(k_{0}\tau), h(k_{0}\tau)], \\[2mm]
\epsilon_{0} e^{(\lambda^{\triangledown}_{1}+\alpha)((k+1)\tau-t)}\phi(t,x),\; &\, t\in((k\tau)^{+}, (k+1)\tau], ~x\in[g(k_{0}\tau), h(k_{0}\tau)],
\end{array} \right.
\end{eqnarray*}
where the positive constant $\epsilon_{0}$ will be determined later.

Next, we proceed to show that the pair $(\underline{u},\underline{v})$ is a lower solution of $(u,v)$ for $t\geq k_{0}\tau$ and $x\in\big[g(k_{0}\tau), h(k_{0}\tau)\big]$. For $(t,x)\in((k\tau)^{+}, (k+1)\tau]\times\big(g(k_{0}\tau), h(k_{0}\tau)\big)$, one can obtain from Assumption (A4) that
\begin{equation*}
\begin{aligned}
& \underline{u}_{t}-d_{1}\Delta \underline{u}+a_{11}\underline{u}-a_{12}\underline{v}\\
\leq&~\epsilon_{0} e^{(\lambda^{\triangledown}_{1}+\alpha)((k+1)\tau-t)}\Big[\phi_{t}-(\lambda^{\triangledown}_{1}+\alpha)\phi-d_{1}\Delta\phi +a_{11} \phi-a_{12} e^{-(\lambda^{\triangledown}_{1}+\alpha)(k\tau-t)}\psi\Big] \\
\leq&~\epsilon_{0} e^{(\lambda^{\triangledown}_{1}+\alpha)((k+1)\tau-t)}\Big[a_{12}\psi\big(1-e^{-(\lambda^{\triangledown}_{1}+\alpha)(k\tau-t)}\big)-\alpha\phi\Big] \\
\leq&~\epsilon_{0} e^{(\lambda^{\triangledown}_{1}+\alpha)((k+1)\tau-t)}\phi\Big[a_{12}C\big(1-e^{(\lambda^{\triangledown}_{1}+\alpha)\tau}\big)-\alpha\Big]\\
\leq&~0,
\end{aligned}
\end{equation*}
and
\begin{equation*}
\begin{aligned}
& \underline{v}_{t}-d_{2}\Delta \underline{v}+a_{22}\underline{v}-f(\underline{u})\\
\leq&~\epsilon_{0} e^{(\lambda^{\triangledown}_{1}+\alpha)\tau}\Big[a_{22} \psi-f'(0)e^{(\lambda^{\triangledown}_{1}+\alpha)(k\tau-t)}\phi+H_{1}\big[\epsilon_{1} e^{(\lambda^{\triangledown}_{1}+\alpha)\tau}\big]^{\kappa_{1}-1}\big[e^{(\lambda^{\triangledown}_{1}+\alpha)(k\tau-t)}\phi\big]^{\kappa_{1}}\\
&+\psi_{t}-d_{2}\Delta\psi \Big] \\
\leq&~\epsilon_{0} e^{(\lambda^{\triangledown}_{1}+\alpha)\tau}\Big[H_{1}\big[\epsilon_{0} e^{(\lambda^{\triangledown}_{1}+\alpha)\tau}\big]^{\kappa_{1}-1}\big[e^{(\lambda^{\triangledown}_{1}+f'(0)\alpha)(k\tau-t)}\phi\big]^{\kappa_{1}}
+\lambda^{\triangledown}_{1}\psi\\
&+f'(0)\phi\big(1-e^{(\lambda^{\triangledown}_{1}+\alpha)(k\tau-t)}\big) \Big]\leq 0
\end{aligned}
\end{equation*}
hold when $\epsilon_{0}$  is sufficiently small. For $x\in\big[g(k_{0}\tau), h(k_{0}\tau)\big]$, one can again obtain from Assumption (A4) that
\begin{equation*}
\begin{aligned}
&\underline{u}\big((k\tau)^{+}, x\big)- G\big(\underline{u}(k\tau,x)\big)=\epsilon_{0} e^{(\lambda^{\triangledown}_{1}+\alpha)\tau} \phi((k\tau)^{+},x)-G\big(\epsilon_{0} \phi(k\tau,x)\big)\\
\leq&~\epsilon_{0}G'(0) e^{(\lambda^{\triangledown}_{1}+\alpha)\tau} \phi(k\tau,x)-\epsilon_{0}G'(0)\phi(k\tau,x)+H_{2}\big(\epsilon_{0}\phi(k\tau,x)\big)^{\kappa_{2}} \\
\leq &~\epsilon_{0}\phi(k\tau,x)\Big[G'(0)\big(e^{(\lambda^{\triangledown}_{1}+\alpha)\tau}-1\big)+H_{2}\big(\epsilon_{0}\phi(k\tau,x)\big)^{\kappa_{2}-1}\Big]\\
\leq &~0,
\end{aligned}
\end{equation*}
and
\begin{equation*}
\begin{aligned}
&\underline{v}\big((k\tau)^{+}, x\big)- \underline{v}(k\tau,x)=\epsilon_{0} e^{(\lambda^{\triangledown}_{1}+\alpha)\tau} \psi((k\tau)^{+},x)-\epsilon_{0} e^{(\lambda^{\triangledown}_{1}+\alpha)\tau} \psi(k\tau,x)\\
=&~\epsilon_{0} e^{(\lambda^{\triangledown}_{1}+\alpha)\tau} \psi(k\tau,x)-\epsilon_{0} e^{(\lambda^{\triangledown}_{1}+\alpha)\tau} \psi(k\tau,x)\leq0
\end{aligned}
\end{equation*}
hold when $\epsilon_{0}$  is sufficiently small. Additionally, such  a sufficiently small $\epsilon_{0}$ can be chosen such that
\begin{equation*}
\begin{aligned}
\underline{u}(k_{0}\tau,x)\leq u(k_{0}\tau,x)\text{~and~}\underline{v}(k_{0}\tau,x)\leq v(k_{0}\tau,x)\text{~for~}x\in\big[g(k_{0}\tau), h(k_{0}\tau)\big].
\end{aligned}
\end{equation*}
Finally, it's easy to obtain that
\begin{equation*}
\begin{aligned}
(\underline{u}, \underline{v})=(0,0)\text{~for~} (t,x)\in[0,\infty)\times\big\{g(k_{0}\tau), h(k_{0}\tau)\big\}.
\end{aligned}
\end{equation*}
Therefore, it follows from \cite[lemma 4.1]{Zhou-Lin-Santos} that
\begin{equation}\label{3-11}
\begin{aligned}
(\underline{u}, \underline{v})\leq (u, v)\text{~for~} (t,x)\in[k_{0}\tau,\infty)\times\big[g(k_{0}\tau), h(k_{0}\tau)\big].
\end{aligned}
\end{equation}

\textbf{(3) A lower solution of model \eqref{Zhou-Lin}}

The principal eigenvalue of problem \eqref{3-8} is denoted by $\lambda^{\epsilon}_{1}$, and the corresponding eigenfunction is denoted by $(\phi^{\epsilon}, \psi^{\epsilon})$. Let
\begin{equation*}
\begin{aligned}
\tilde{x}=x_{0}+\frac{h(k_{0}\tau)-x_{0}}{h(t)-x_{0}}(x-x_{0}),~x\in[x_{0}, h(t)].
\end{aligned}
\end{equation*}
Because $h_{\infty}<+\infty$, a sufficiently large $k_{1}(\geq k_{0})$ can be found such that when $t\geq k_{1}\tau$,
\begin{eqnarray}\label{3-12}
h'(t)<\min\bigg\{\epsilon\frac{d_{2}(h(k_{0}\tau)-x_{0})}{d_{1}(h_{\infty}-x_{0})},\epsilon\frac{h(k_{0}\tau)-x_{0}}{h_{\infty}-x_{0}}\bigg\}.
\end{eqnarray}
Now we construct the following functions
\begin{eqnarray*}
\tilde{v}(t,x)=
\epsilon_{1} e^{(\lambda^{\epsilon}_{1}+\alpha)\tau}\psi^{\epsilon}(t,\tilde{x}),\; &\, (t,x)\in [0, \infty)\times[x_{0}, h(t)],
\end{eqnarray*}
\begin{eqnarray*}
\tilde{u}(t,x)=
\left\{
\begin{array}{ll}
\epsilon_{1} \phi^{\epsilon}(t,\tilde{x}),\; &\, t=k\tau,~ x\in[x_{0}, h(t)], \\[2mm]
\epsilon_{1} e^{(\lambda^{\epsilon}_{1}+\alpha)\tau} \phi^{\epsilon}(t,\tilde{x}),\; &\, t=(k\tau)^{+}, ~x\in[x_{0}, h(t)], \\[2mm]
\epsilon_{1} e^{(\lambda^{\epsilon}_{1}+\alpha)((k+1)\tau-t)}\phi^{\epsilon}(t,\tilde{x}),\; &\, t\in((k\tau)^{+}, (k+1)\tau], ~x\in[x_{0}, h(t)],
\end{array} \right.
\end{eqnarray*}
where $\alpha$ is defined in \eqref{3-14} where $\lambda^{\triangledown}_{1}$ replaced by $\lambda^{\epsilon}_{1}$, and positive constant $\epsilon_{1}(\leq\epsilon_{0})$ will be determined later.

Next, we proceed to show that the pair $(\tilde{u},\tilde{v})$ is a lower solution of model \eqref{Zhou-Lin} for $t\geq k_{1}\tau$ and $x\in[x_{0}, h(t)]$. For $t\in((k\tau)^{+}, (k+1)\tau]\bigcap ((k_{1}\tau)^{+},+\infty)$ and $x\in(x_{0}, h(t))$, it follows from \eqref{3-8}, \eqref{3-10}, \eqref{3-12} and Assumption (A4) that
\begin{align*}
& \tilde{u}_{t}-d_{1}\Delta \tilde{u}+a_{11}\tilde{u}-a_{12}\tilde{v}\\
\leq&~\epsilon_{1} e^{(\lambda^{\epsilon}_{1}+\alpha)((k+1)\tau-t)}\bigg[\phi^{\epsilon}_{t}
-\frac{(x-x_{0})(h(k_{0}\tau)-x_{0})h'(t)}{(h(t)-x_{0})^{2}}\phi^{\epsilon}_{\tilde{x}}+a_{11} \phi^{\epsilon}\\
&~-d_{1}\Big(\frac{h(k_{0}\tau)-x_{0}}{h(t)-x_{0}}\Big)^{2}\phi^{\epsilon}_{\tilde{x}\tilde{x}}-a_{12}e^{-(\lambda^{\epsilon}_{1}
+\alpha)(k\tau-t)}\psi^{\epsilon}-(\lambda^{\epsilon}_{1}+\alpha)\phi^{\epsilon}\bigg] \\
\leq&~\epsilon_{1} e^{(\lambda^{\epsilon}_{1}+\alpha)((k+1)\tau-t)}\bigg[-\alpha \phi^{\epsilon}+\Big[\epsilon-\frac{h'(t)(x-x_{0})(h(k_{0}\tau)-x_{0})}{(h(t)-x_{0})^{2}}\Big]\phi^{\epsilon}_{\tilde{x}}\\
&~+d_{1}\Big[1-\Big(\frac{h(k_{0}\tau)-x_{0}}{h(t)-x_{0}}\Big)^{2}\Big]\phi^{\epsilon}_{\tilde{x}\tilde{x}}+a_{12}\Big[1-e^{-(\lambda^{\epsilon}_{1}
+\alpha)(k\tau-t)}\Big]\psi^{\epsilon}\bigg] \\
\leq&~\epsilon_{1} e^{(\lambda^{\epsilon}_{1}+\alpha)((k+1)\tau-t)}\bigg[-d_{1}\Big[1-\Big(\frac{h(k_{0}\tau)-x_{0}}{h(t)-x_{0}}\Big)^{2}\Big]\Big[\frac{\epsilon^{2}}{4d^{2}_{1}}
+\Big(\frac{\pi}{h(k_{0}\tau)-g(k_{0}\tau)}\Big)^{2}\Big]\phi^{\epsilon}\\
&~+\Big[\epsilon\Big(\frac{h(k_{0}\tau)-x_{0}}{h(t)-x_{0}}\Big)^{2}-\frac{h'(t)(x-x_{0})(h(k_{0}\tau)-x_{0})}{(h(t)-x_{0})^{2}}\Big]\phi^{\epsilon}_{\tilde{x}}\bigg] \\
&~-\alpha \phi^{\epsilon}+\Big[\epsilon-\frac{h'(t)(x-x_{0})(h(k_{0}\tau)-x_{0})}{(h(t)-x_{0})^{2}}\Big]\phi^{\epsilon}_{\tilde{x}}
+a_{12}\Big[1-e^{(\lambda^{\epsilon}_{1}+\alpha)\tau}\Big]\psi^{\epsilon}\bigg] \\
\leq&~\epsilon_{1} e^{(\lambda^{\epsilon}_{1}+\alpha)((k+1)\tau-t)}\bigg[-\alpha \phi^{\epsilon}+a_{12}\Big[1-e^{(\lambda^{\epsilon}_{1}+\alpha)\tau}\Big]\psi^{\epsilon}\bigg] \\
\leq&~\epsilon_{1} e^{(\lambda^{\epsilon}_{1}+\alpha)((k+1)\tau-t)}\phi^{\epsilon}\bigg[a_{12}\big[1-e^{(\lambda^{\epsilon}_{1}
+\alpha)\tau}\big]C-\alpha \bigg] \\
\leq&~0,
\end{align*}
and
\begin{align*}
\begin{aligned}
& \tilde{v}_{t}-d_{2}\Delta \tilde{v}+a_{22}\tilde{v}-f(\tilde{u})\\
\leq&~\epsilon_{1} e^{(\lambda^{\epsilon}_{1}+\alpha)\tau}\bigg[a_{22} \psi^{\epsilon}-f'(0)e^{(\lambda^{\epsilon}_{1}+\alpha)(k\tau-t)}\phi^{\epsilon}+H_{1}\big[\epsilon_{1} e^{(\lambda^{\epsilon}_{1}+\alpha)\tau}\big]^{\kappa_{1}-1}\big[e^{(\lambda^{\epsilon}_{1}+\alpha)(k\tau-t)}\phi^{\epsilon}\big]^{\kappa_{1}}\\
&+\psi^{\epsilon}_{t}-\frac{h'(t)(x-x_{0})(h(k_{0}\tau)-x_{0})}{(h(t)-x_{0})^{2}}\psi^{\epsilon}_{\tilde{x}}-d_{2}\Big(\frac{h(k_{0}\tau)
-x_{0}}{h(t)-x_{0}}\Big)^{2}\psi^{\epsilon}_{\tilde{x}\tilde{x}} \bigg] \\
\leq&~\epsilon_{1} e^{(\lambda^{\epsilon}_{1}+\alpha)\tau}\bigg[\lambda^{\epsilon}_{1}\psi^{\epsilon}+f'(0)\Big[1-e^{(\lambda^{\epsilon}_{1}
+\alpha)(k\tau-t)}\Big]\phi^{\epsilon}+H_{1}\big[\epsilon_{1} e^{(\lambda^{\epsilon}_{1}+\alpha)\tau}\big]^{\kappa_{1}-1}\big[e^{(\lambda^{\triangledown}_{1}+f'(0)\alpha)(k\tau-t)}\phi\big]^{\kappa_{1}}\\
&+d_{2}\Big[1-\Big(\frac{h(k_{0}\tau)-x_{0}}{h(t)-x_{0}}\Big)^{2}\Big]\psi^{\epsilon}_{\tilde{x}\tilde{x}}+\Big[\epsilon\frac{d_{2}}{d_{1}}
-\frac{h'(t)(x-x_{0})(h(k_{0}\tau)-x_{0})}{(h(t)-x_{0})^{2}}\Big]\psi^{\epsilon}_{\tilde{x}} \bigg]\\
\leq&~\epsilon_{1} e^{(\lambda^{\epsilon}_{1}+\alpha)\tau}\bigg[\lambda^{\epsilon}_{1}\psi^{\epsilon}+H_{1}\big[\epsilon_{1} e^{(\lambda^{\epsilon}_{1}+\alpha)\tau}\big]^{\kappa_{1}-1}\big[e^{(\lambda^{\epsilon}_{1}+\alpha)(k\tau-t)}\phi\big]^{\kappa_{1}}
+\Big[\epsilon\frac{d_{2}}{d_{1}}\Big(\frac{h(k_{0}\tau)-x_{0}}{h(t)-x_{0}}\Big)^{2}\\
&-\frac{h'(t)(x-x_{0})(h(k_{0}\tau)-x_{0})}{(h(t)-x_{0})^{2}}\Big]\psi^{\epsilon}_{\tilde{x}} -d_{2}\Big[1-\Big(\frac{h(k_{0}\tau)-x_{0}}{h(t)-x_{0}}\Big)^{2}\Big]\Big[\frac{\epsilon^{2}}{4d^{2}_{1}}+\Big(\frac{\pi}{h(k_{0}\tau)
-g(k_{0}\tau)}\Big)^{2}\Big]\psi^{\epsilon} \bigg]\\
\leq&~0
\end{aligned}
\end{align*}
hold when $\epsilon_{1}$  is sufficiently small. From Assumption (A4), one can then obtain that
\begin{equation*}
\begin{aligned}
&\tilde{u}\big((k\tau)^{+}, x\big)- G\big(\tilde{u}(k\tau,x)\big)=\epsilon_{1} e^{(\lambda^{\epsilon}_{1}+\alpha)\tau} \phi^{\epsilon}((k\tau)^{+},\tilde{x})-G\big(\epsilon_{1} \phi^{\epsilon}(k\tau,\tilde{x})\big)\\
\leq&~\epsilon_{1}G'(0) e^{(\lambda^{\epsilon}_{1}+\alpha)\tau} \phi^{\epsilon}(k\tau,\tilde{x})-\epsilon_{1}G'(0)\phi^{\epsilon}(k\tau,\tilde{x})+H_{2}\big(\epsilon_{1}\phi^{\epsilon}(k\tau,\tilde{x})\big)^{\kappa_{2}} \\
\leq &~\epsilon_{1}\phi^{\epsilon}(k\tau,\tilde{x})\Big[G'(0)\big(e^{(\lambda^{\epsilon}_{1}+\alpha)\tau}-1\big)+
H_{2}\big(\epsilon_{1}\phi^{\epsilon}(k\tau,\tilde{x})\big)^{\kappa_{2}-1}\Big]\\
\leq &~0,
\end{aligned}
\end{equation*}
and
\begin{equation*}
\begin{aligned}
&\tilde{v}\big((k\tau)^{+}, x\big)-\tilde{v}(k\tau,x)=\epsilon_{1} e^{(\lambda^{\epsilon}_{1}+\alpha)\tau} \psi^{\epsilon}((k\tau)^{+},\tilde{x})-\epsilon_{1} e^{(\lambda^{\epsilon}_{1}+\alpha)\tau} \psi^{\epsilon}(k\tau,\tilde{x})\\
=&~\epsilon_{1} e^{(\lambda^{\epsilon}_{1}+\alpha)\tau} \psi^{\epsilon}(k\tau,\tilde{x})-\epsilon_{1} e^{(\lambda^{\epsilon}_{1}+\alpha)\tau} \psi^{\epsilon}(k\tau,\tilde{x})\\
=&~0
\end{aligned}
\end{equation*}
hold on $[x_{0}, h(t)]$ when $\epsilon_{1}$  is sufficiently small. For the initial conditions, there exists a sufficiently small positive constant $\epsilon_{1}$ such that
\begin{equation*}
\begin{aligned}
\tilde{u}(k_{1}\tau,x)\leq u(k_{1}\tau,x)\text{~and~}\tilde{v}(k_{1}\tau,x)\leq v(k_{1}\tau,x)\text{~for~}x\in[x_{0}, h(t)].
\end{aligned}
\end{equation*}
By applying \eqref{3-11}, it follows
\begin{equation*}
\begin{aligned}
(\underline{u}(t,x_{0}), \underline{v}(t,x_{0}))\leq (u(t,x_{0}), v(t,x_{0}))
\end{aligned}
\end{equation*}
when $t\geq k_{1}\tau$. Then, a sufficiently small $\epsilon_{1}$ can be obtained such that the left boundary condition holds. At the end,
it is easily obtained that the right boundary condition also holds. In summary, the pair $(\tilde{u},\tilde{v})$ is a lower solution of
model \eqref{Zhou-Lin} for $t\geq k_{1}\tau$ and $x\in[x_{0}, h(t)]$.

With the help of \cite[lemma 4.1]{Zhou-Lin-Santos}, it follows
\begin{equation*}
\begin{aligned}
(\tilde{u}(t,x), \tilde{v}(t,x))\leq (u(t,x), v(t,x))
\end{aligned}
\end{equation*}
for $(t,x)\in[k_{1}\tau, \infty)\times [x_{0}, h(t)]$. This combined with the equations for $\tilde{u}$ and $\tilde{v}$ yields that
\begin{equation}\label{3-15}
\begin{aligned}
\Big(\lim\limits_{t\rightarrow +\infty}u_{x}(t,h(t)), \lim\limits_{t\rightarrow +\infty}v_{x}(t,h(t)) \Big)\leq \Big(\lim\limits_{t\rightarrow +\infty}\tilde{u}_{x}(t,h(t)), \lim\limits_{t\rightarrow +\infty}\tilde{v}_{x}(t,h(t))\Big)<(0,0).
\end{aligned}
\end{equation}
By using the seventh equation of model \eqref{Zhou-Lin}, one can obtain that
\begin{equation*}
\begin{aligned}
\lim\limits_{t\rightarrow +\infty}u_{x}(t,h(t))=-\frac{\mu_{2}}{\mu_{1}}\lim\limits_{t\rightarrow +\infty}v_{x}(t,h(t))>0,
\end{aligned}
\end{equation*}
which shows a contradiction with \eqref{3-15}. Therefore, the desired conclusion holds.
\end{proof}
\end{lemma}

\autoref{lemma 3-5} shows that if $h_{\infty}<+\infty$, then $-g_{\infty}<+\infty$, and vice versa. This lemma classifies the double moving fronts $g(t)$ and $h(t)$ as falling within two cases: simultaneous finite or simultaneous infinite. Next, we will explore the spatial distributions of infected individuals and
pathogenic bacteria in each of these two cases as time tends to infinity.
\begin{theorem}\label{theorem 3-2}
If $h_{\infty}-g_{\infty}<+\infty$, then $\lambda^{\vartriangle}_{1}\big(G'(0),(g_{\infty},h_{\infty})\big)\geq 0$ and
\begin{equation*}
\lim\limits_{t\rightarrow\infty}\|u(t, \cdot)\|_{\mathbb{C}[g(t), h(t)]}+\lim\limits_{t\rightarrow\infty}\|v(t, \cdot)\|_{\mathbb{C}[g(t), h(t)]}=0.
\end{equation*}
\begin{proof}
Arguing indirectly, suppose that $\lambda^{\vartriangle}_{1}\big(G'(0),(g_{\infty},h_{\infty})\big)<0$. According to \autoref{lemma 3-1}\textcolor{blue}{$(1, 3)$}, for any sufficiently small $\epsilon$, one can obtain that
\begin{equation*}
\lambda^{\vartriangle}_{1}\big(G'(0),(g_{\infty}+\epsilon, h_{\infty}-\epsilon)\big)<0.
\end{equation*}
Moreover, for any given $\epsilon$, there exists a sufficiently large positive integer $k_{0}$ dependent on $\epsilon$ such that
\begin{equation*}
\big(g_{\infty}, h_{\infty}-\epsilon\big)<\big(g(t),h(t)\big)<\big(g_{\infty}+\epsilon, h_{\infty})
\end{equation*}
when $t\geq k_{0}\tau$.

Now we consider the  problem
\begin{eqnarray*}
\left\{
\begin{array}{ll}
\underline{u}_{t}=d_{1}\Delta \underline{u}-a_{11}\underline{u}+a_{12}\underline{v},\; &\, t\in((k\tau)^{+}, (k+1)\tau], x\in\Omega_{1}, \\[2mm]
\underline{v}_{t}=d_{2}\Delta \underline{v}-a_{22}\underline{v}+f(\underline{u}),\; &\, t\in((k\tau)^{+}, (k+1)\tau], x\in\Omega_{1}, \\[2mm]
\underline{u}=0, \underline{v}=0,\; &\, t\in(k\tau, (k+1)\tau], ~~~x\in\partial\Omega_{1},\\[2mm]
\underline{u}((k\tau)^{+},x)=G(\underline{u}(k\tau,x)), \; &\, x\in\Omega_{1}, \\[2mm]
\underline{v}((k\tau)^{+},x)=\underline{v}((k\tau),x), \; &\,  x\in\Omega_{1}, \\[2mm]
\underline{u}(k_{0}\tau,x)=u(k_{0}\tau,x),\; &\,x\in\overline{\Omega}_{1}, \\
\underline{v}(k_{0}\tau,x)=v(k_{0}\tau,x),\; &\,x\in\overline{\Omega}_{1}, k=k_{0},k_{0}+1,\cdots,
\end{array} \right.
\end{eqnarray*}
where $\Omega_{1}=(g_{\infty}+\epsilon, h_{\infty}-\epsilon)$. With the help of the comparison principle, it follows that
\begin{equation}\label{3-16}
\big(\underline{u}, \underline{v}\big)\leq \big(u,v\big), ~~~~t\geq k_{0}\tau, ~x\in\overline{\Omega}_{1}.
\end{equation}
Noting that $\lambda^{\vartriangle}_{1}\big(G'(0),\Omega_{1}\big)<0$, it follows from \autoref{lemma 3-1}\textcolor{blue}{$(5)$} that
\begin{equation}\label{3-17}
\lim\limits_{m\rightarrow+\infty}\big(\underline{u}(t+m\tau,x),\underline{v}(t+m\tau,x)\big)=\big(\underline{U}^{\epsilon}(t,x), \underline{V}^{\epsilon}(t,x)\big)
\end{equation}
uniformly for $(t,x)\in[0,\tau]\times\overline{\Omega}_{1}$. Here $\big(\underline{U}^{\epsilon}, \underline{V}^{\epsilon}\big)$ denotes the unique positive periodic solution to the problem
\begin{eqnarray}\label{3-23}
\left\{
\begin{array}{ll}
(\underline{U}^{\epsilon})_{t}=d_{1} \Delta \underline{U}^{\epsilon}-a_{11}\underline{U}^{\epsilon}+a_{12}\underline{V}^{\epsilon},\; &\, t\in(0^{+}, \tau], x\in\Omega_{1}, \\[2mm]
(\underline{V}^{\epsilon})_{t}=d_{2} \Delta \underline{V}^{\epsilon}-a_{22}\underline{V}^{\epsilon}+f(\underline{U}^{\epsilon}),\; &\, t\in(0^{+}, \tau], x\in\Omega_{1}, \\[2mm]
\underline{U}^{\epsilon}=0, \underline{V}^{\epsilon}=0,\; &\, t\in[0, \tau], x\in\partial\Omega_{1}, \\[2mm]
\underline{U}^{\epsilon}(0,x)=\underline{U}^{\epsilon}(\tau,x), \underline{V}^{\epsilon}(0,x)=\underline{V}^{\epsilon}(\tau,x),\; &\,x\in\overline{\Omega}_{1},\\[2mm]
\underline{U}^{\epsilon}(0^{+},x)=G\big(\underline{U}^{\epsilon}(0,x)\big), \underline{V}^{\epsilon}(0^{+},x)=\underline{V}^{\epsilon}(0,x), \; &\, x\in\Omega_{1}.
\end{array} \right.
\end{eqnarray}
From  \eqref{3-16} and \eqref{3-17}, one can obtain that
\begin{equation}\label{3-19}
\liminf\limits_{m\rightarrow+\infty}\big(u(t+m\tau,x),v(t+m\tau,x)\big)\geq\big(\underline{U}^{\epsilon}(t,x), \underline{V}^{\epsilon}(t,x)\big)
\end{equation}
for $(t,x)\in[0,\tau]\times\overline{\Omega}_{1}$.

In what follows, we consider the problem
\begin{eqnarray*}
\left\{
\begin{array}{ll}
\overline{u}_{t}=d_{1}\Delta \overline{u}-a_{11}\overline{u}+a_{12}\overline{v},\; &\, t\in((k\tau)^{+}, (k+1)\tau], x\in\Omega_{2}, \\[2mm]
\overline{v}_{t}=d_{2}\Delta \overline{v}-a_{22}\overline{v}+f(\overline{u}),\; &\, t\in((k\tau)^{+}, (k+1)\tau], x\in\Omega_{2}, \\[2mm]
\big(\overline{u}, \overline{v}\big)=(0,0),\; &\, t\in(k\tau, (k+1)\tau], ~~~x\in\partial\Omega_{2},\\[2mm]
\overline{u}((k\tau)^{+},x)=G(\overline{u}(k\tau,x)), \; &\, x\in\Omega_{2}, \\[2mm]
\overline{v}((k\tau)^{+},x)=\overline{v}((k\tau),x), \; &\,  x\in\Omega_{2}, \\[2mm]
\big(\overline{u}(0,x), \overline{v}(0,x)\big)=\big(\overline{u}_{0}(x), \overline{v}_{0}(x)\big),\; &\,x\in\overline{\Omega}_{2}, ~~k=0,1,2,\cdots,
\end{array} \right.
\end{eqnarray*}
where $\Omega_{2}=(g_{\infty}, h_{\infty})$, and the initial functions satisfy
\begin{eqnarray*}
\overline{u}_{0}(x)=
\left\{
\begin{array}{ll}
u_{0}(x),\; &\, x\in[-h(0), h(0)], \\[2mm]
0,\; &\, x\in\overline{\Omega}_{2}\setminus[-h(0), h(0)],
\end{array} \right.
\overline{v}_{0}(x)=
\left\{
\begin{array}{ll}
v_{0}(x),\; &\, x\in[-h(0), h(0)], \\[2mm]
0,\; &\, x\in\overline{\Omega}_{2}\setminus[-h(0), h(0)].
\end{array} \right.
\end{eqnarray*}
With the help of the comparison principle, it follows that
\begin{equation}\label{3-20}
\big(u, v\big)\leq \big(\overline{u}, \overline{u}\big), ~~(t,x)\in[0, \infty)\times \overline{\Omega}_{2}.
\end{equation}
Due to the condition $\lambda^{\vartriangle}_{1}\big(G'(0),\Omega_{2}\big)<0$, one can obtain from \autoref{lemma 3-1}\textcolor{blue}{$(5)$} that
\begin{equation}\label{3-21}
\lim\limits_{m\rightarrow+\infty}\big(\overline{u}(t+m\tau,x),\overline{v}(t+m\tau,x)\big)=\big(U(t,x), V(t,x)\big)
\end{equation}
uniformly on $[0,\tau]\times\overline{\Omega}_{2}$. Here $\big(U, V\big)$ denotes the unique solution to the problem
\begin{eqnarray*}
\left\{
\begin{array}{ll}
U_{t}=d_{1} \Delta U-a_{11}U+a_{12}V,\; &\, t\in(0^{+}, \tau], x\in\Omega_{2}, \\[2mm]
V_{t}=d_{2} \Delta V-a_{22}V+f(U),\; &\, t\in(0^{+}, \tau], x\in\Omega_{2}, \\[2mm]
U(t,x)=0, V(t,x)=0,\; &\, t\in[0, \tau], x\in\partial\Omega_{2}, \\[2mm]
U(0,x)=U(\tau,x), V(0,x)=V(\tau,x),\; &\,x\in\overline{\Omega}_{2},\\[2mm]
U(0^{+},x)=G\big(U(0,x)\big), V(0^{+},x)=V(0,x), \; &\, x\in\Omega_{2}.
\end{array} \right.
\end{eqnarray*}
By \eqref{3-20} and \eqref{3-21}, it follows that
\begin{equation}\label{3-22}
\limsup\limits_{m\rightarrow+\infty}\big(u(t+m\tau,x),v(t+m\tau,x)\big)\leq\big(U(t,x), V(t,x)\big)
\end{equation}
for $(t,x)\in[0,\tau]\times\overline{\Omega}_{2}$.

Using standard regularity theory for partial differential equations for the problem \eqref{3-23} we find that for arbitrary $\alpha\in(0,1)$, we have that
\begin{equation*}
\big(\underline{U}^{\epsilon}(t,x), \underline{V}^{\epsilon}(t,x)\big)\rightarrow \big(U(t,x), V(t,x)\big)\text{~in~}\big[\mathbb{C}^{(1+\alpha)/2, 1+\alpha}_{\text{loc}}\big([0, \tau]\times\Omega_{2}\big)\big]^{2}\text{~as~}\epsilon\rightarrow0.
\end{equation*}
This combined with \eqref{3-19}, \eqref{3-22}, and the arbitrariness of $\epsilon$ yields that
\begin{equation*}
\lim\limits_{m\rightarrow+\infty}\big(u(t+m\tau,x),v(t+m\tau,x)\big)=\big(U(t,x), V(t,x)\big)
\end{equation*}
on $[0,\tau]\times\Omega_{2}$. Then it is easy to derive that
\begin{equation*}
\lim\limits_{m\rightarrow+\infty}\big(u_{x}(t+m\tau,h(t+m\tau)),v_{x}(t+m\tau,h(t+m\tau))\big)=\big(U_{x}(t,h_{\infty}), V_{x}(t,h_{\infty})\big)<(0,0)
\end{equation*}
and
\begin{equation*}
\lim\limits_{m\rightarrow+\infty}\big(u_{x}(t+m\tau,g(t+m\tau)),v_{x}(t+m\tau,g(t+m\tau))\big)=\big(U_{x}(t,g_{\infty}), V_{x}(t,g_{\infty})\big)>(0,0)
\end{equation*}
uniformly on $[0,\tau]$. Additionally, since
\begin{equation*}
h'(t+m\tau)=-\mu_{1}u_{x}(t+m\tau,h(t+m\tau))-\mu_{2}v_{x}(t+m\tau,h(t+m\tau))
\end{equation*}
and
\begin{equation*}
g'(t+m\tau)=-\mu_{1}u_{x}(t+m\tau,g(t+m\tau))-\mu_{2}v_{x}(t+m\tau,g(t+m\tau)),
\end{equation*}
one can obtain that
\begin{equation*}
\lim\limits_{t\rightarrow+\infty}h'(t)=\lim\limits_{t\rightarrow+\infty}-g'(t)>0,
\end{equation*}
which is a contradiction to $h_{\infty}-g_{\infty}<+\infty$. Therefore, $\lambda^{\vartriangle}_{1}\big(G'(0),(g_{\infty},h_{\infty})\big)\geq 0$.

With the help of $\lambda^{\vartriangle}_{1}\big(G'(0),(g_{\infty},h_{\infty})\big)\geq 0$ and \autoref{lemma 3-1}\textcolor{blue}{$(4)$}, it follows that
\begin{equation*}
\lim\limits_{t\rightarrow\infty}\big[\|\overline{u}(t, \cdot)\|_{\mathbb{C}(\overline{\Omega}_{2})}+\||\overline{v}(t, \cdot)\|_{\mathbb{C}(\overline{\Omega}_{2})}\big]=0.
\end{equation*}
This combined with \eqref{3-20} yields that
\begin{equation*}
\lim\limits_{t\rightarrow\infty}\big[\|u(t, \cdot)\|_{\mathbb{C}([g(t), h(t)])}+\|v(t, \cdot)\|_{\mathbb{C}([g(t), h(t)])}\big]=0.
\end{equation*}
This completes the proof.
\end{proof}
\end{theorem}
\autoref{theorem 3-2} tells us that the diseases transmitted by the faecal-oral route become extinct when the double moving fronts are simultaneously finite. Next, we explore the case when the double moving fronts are simultaneously infinite.

\begin{theorem}\label{theorem 3-3}
If $h_{\infty}-g_{\infty}=+\infty$ and $\lambda^{\vartriangle}_{1}\big(G'(0),(g_{\infty},h_{\infty})\big)<0$, then for any $\Omega_{t}\subset\subset [0, \tau]\times \mathbb{R}$
\begin{equation*}
\lim\limits_{m\rightarrow+\infty}u(t+m\tau,x)=U_{1}(t)~\text{~and~}~\lim\limits_{m\rightarrow+\infty}v(t+m\tau,x)=V_{1}(t)
\end{equation*}
for $(t,x)\in \Omega_{t}$. Here $(U_{1}(t), V_{1}(t))$ stands for the unique solution to problem
\begin{eqnarray}\label{3-24}
\left\{
\begin{array}{ll}
u_{t}=a_{12}v-a_{11}u,~~~~~~~~t\in(0^{+}, \tau], \\[2mm]
v_{t}=f(u)-a_{22}v, ~~~~~~~~t\in(0^{+}, \tau], \\[2mm]
u(0)=u(\tau), v(0)=v(\tau),\\[2mm]
u(0^{+})=G(u(0)), v(0^{+})=v(0).
\end{array} \right.
\end{eqnarray}
\begin{proof}
Due to the length of the proof  it is divided into three parts for easier reading.

\textbf{(1) The existence and uniqueness of the positive solution of problem \eqref{3-24}}

To begin with, introduce the following eigenvalue problem
\begin{eqnarray}\label{3-25}
\left\{
\begin{array}{ll}
\Phi'(t)=(\lambda-a_{11})\Phi(t)+a_{12}\Psi(t),\; &\, t\in(0^{+}, \tau], \\[2mm]
\Psi'(t)=f'(0)\Phi(t)+(\lambda-a_{22})\Phi(t),\; &\, t\in(0^{+}, \tau], \\[2mm]
\big(\Phi(0),\Psi(0)\big)=\big(\Phi(\tau), \Psi(\tau)\big),\; &\,\\[2mm]
\big(\Phi(0^{+}), \Psi(0^{+})\big)=\big(G'(0)\Phi(0),\Psi(0)\big).\; &\,
\end{array} \right.
\end{eqnarray}
By using the same method as in the proof of \autoref{lemma 3-3}, the periodic eigenvalue problem \eqref{3-25} has the principal eigenvalue $\tilde{\lambda}_{1}$ with a positive eigenfunction $(\Phi(t),\Psi(t) )$ for $t\in[0, \tau]$. By using $\lambda^{\vartriangle}_{1}\big(G'(0),(g_{\infty},\\h_{\infty})\big)<0$, a sufficiently small positive number $\sigma$ can be found such that $\tilde{\lambda}_{1}+\sigma<0$.

Now, construct the super-and sub-solutions of problem \eqref{3-24}. Let
\begin{eqnarray*}
\underline{u}(t)=
\left\{
\begin{array}{ll}
\epsilon \Phi(0),\; &\, t=0, \\[2mm]
\epsilon e^{(\tilde{\lambda}_{1}+\sigma)\tau} \Phi(0^{+}),\; &\, t=0^{+}, \\[2mm]
\epsilon e^{(\tilde{\lambda}_{1}+\sigma)(\tau-t)}\Phi(t),\; &\, t\in(0^{+}, \tau],
\end{array} \right.
\underline{v}(t)=
\left\{
\begin{array}{ll}
\epsilon \Psi(0),\; &\, t=0, \\[2mm]
\epsilon e^{(\tilde{\lambda}_{1}+\sigma)\tau} \Psi(0^{+}),\; &\, t=0^{+}, \\[2mm]
\epsilon e^{(\tilde{\lambda}_{1}+\sigma)(\tau-t)}\Psi(t),\; &\, t\in(0^{+}, \tau],
\end{array} \right.
\end{eqnarray*}
where $\epsilon$ is a sufficiently small positive constant. A simple calculation shows that $(\underline{u}(t),\underline{v}(t))$ is a lower solution of problem \eqref{3-24}. In addition, a simple calculation provides that  $(C_{2}, C_{3})$ is an super-solution of problem \eqref{3-24}, where $(C_{2}, C_{3})$ is defined in \autoref{lemma 3-2}. Therefore, problem \eqref{3-24} has a maximum positive periodic solution $(\overline{U}_{1}, \overline{V}_{1})$ and a minimum positive periodic solution $(\underline{U}_{1}, \underline{V}_{1})$ by the upper and lower solution method.

Next, we prove that the solution of problem \eqref{3-24} is unique. Arguing indirectly, assume that $(\overline{U}_{1}, \overline{V}_{1})$ and $(\underline{U}_{1}, \underline{V}_{1})$ are two distinct solutions of problem \eqref{3-24}. Define the following set
\begin{equation*}
S:=\Big\{s\in[0,1]:s(\underline{U}_{1}, \underline{V}_{1})\leq (\overline{U}_{1}, \overline{V}_{1}), t\in[0,\tau] \Big\}.
\end{equation*}
We now prove $1\in S$. On the contrary, suppose that $s_{0}=\sup S<1$. With loss of generality, assume that
\begin{equation}\label{3-26}
s_{0}=\sup \Big\{s\in[0,1]:s\underline{V}_{1}\leq \overline{V}_{1}, t\in[0,\tau]\Big\}.
\end{equation}
With the help of Assumption (A2), a standard calculation yields that
\begin{equation*}
\begin{aligned}
&\big(\overline{V}_{1}-s_{0}\underline{V}_{1}\big)_{t}+a_{22}\big(\overline{V}_{1}-s_{0}\underline{V}_{1}\big)=f\big(\overline{V}_{1}\big)
-s_{0}f\big(\underline{V}_{1}\big)\\
\geq& f\big(s_{0}\underline{V}_{1}\big)-s_{0}f\big(\underline{V}_{1}\big)\geq 0
\end{aligned}
\end{equation*}
for $t\in(0^{+}, \tau]$. Moreover, it is not hard to obtain that
\begin{equation*}
\overline{V}_{1}(0^{+})-s_{0}\underline{V}_{1}(0^{+})=\overline{V}_{1}(\tau)-s_{0}\underline{V}_{1}(\tau)\geq 0.
\end{equation*}
By applying the strong maximum principle to $\overline{V}_{1}-s_{0}\underline{V}_{1}$, one of the following two conclusions must hold:

(1) $\overline{V}_{1}-s_{0}\underline{V}_{1}>0$ holds for $t\in(0^{+},\tau]$. Note that $\overline{V}_{1}$ and $\underline{V}_{1}$ are $\tau$-periodic solutions,
that is, $\big(\overline{V}_{1}(0,x), \underline{V}_{1}(0,x)\big)=\big(\overline{V}_{1}(\tau,x), \underline{V}_{1}(\tau,x)\big)$. This implies that $\overline{V}_{1}-s_{0}\underline{V}_{1}>0$ holds when $t\in[0,\tau]$. Therefore, a sufficiently small $\epsilon_{0}$ can be found such that $\overline{V}_{1}-s_{0}\underline{V}_{1}>\epsilon_{0} \underline{V}_{1}$. This is in  contradiction to \eqref{3-26}.

(2) $\overline{V}_{1}-s_{0}\underline{V}_{1}\equiv0$ holds for $t\in[0,\tau]$. This implies that $s_{0}=\overline{V}_{1}/\underline{V}_{1}\geq 1$. This is a  contradiction to the assumption that $s_{0}=\sup S<1$.

To conclude, $1\in S$, that is, $\underline{U}_{1}\leq \overline{U}_{1}$ and $\underline{V}_{1}\leq \overline{V}_{1}$ on $[0,\tau]$. Similarly, the conclusion that $\overline{U}_{1}\leq \underline{U}_{1}$ and $\overline{V}_{1}\leq \underline{V}_{1}$ on $[0,\tau]$ can also be proved. Therefore, $(\underline{U}_{1}, \underline{V}_{1})=(\overline{U}_{1}, \overline{V}_{1})$ on $[0,\tau]$, that is, the solution of  problem \eqref{3-24} is unique. Let the unique solution of problem \eqref{3-24} be $(U_{1}, V_{1})$.

\textbf{(2)} \textbf{For any }\bm{$\Omega_{t}\subset\subset [0, \tau]\times \mathbf{R}$}, \bm{$\liminf\limits_{m\rightarrow+\infty}u(t+m\tau,x)
\geq U_{1}(t)$}\textbf{ and} \bm{$\liminf\limits_{m\rightarrow+\infty}v(t+m\tau,x)\geq V_{1}(t)$} \textbf{for} \bm{$(t,x)\in\Omega_{t}$}

Since $\lambda^{\vartriangle}_{1}\big(G'(0),(g_{\infty},h_{\infty})\big)<0$, it follows from  \autoref{lemma 3-1}\textcolor{blue}{$(1,3)$} that there exists a sufficiently large positive constant $S$ such that $\lambda^{\vartriangle}_{1}\big(G'(0),(-S,S)\big)<0$. For any given $S$, a sufficiently large positive integer $k_{2}$ can be found such that $h(t)-g(t)\geq 2S$ when $t\geq k_{2}\tau$.

Now, take into account the following
\begin{eqnarray*}
\left\{
\begin{array}{ll}
\underline{u}_{t}=d_{1}\Delta \underline{u}-a_{11}\underline{u}+a_{12}\underline{v},\; &\, t\in((k\tau)^{+}, (k+1)\tau], x\in\Omega_{3}, \\[2mm]
\underline{v}_{t}=d_{2}\Delta \underline{v}-a_{22}\underline{v}+f(\underline{u}),\; &\, t\in((k\tau)^{+}, (k+1)\tau], x\in\Omega_{3}, \\[2mm]
\underline{u}=0, \underline{v}=0,\; &\, t\in(k\tau, (k+1)\tau], ~~~x\in\partial\Omega_{3},\\[2mm]
\underline{u}((k\tau)^{+},x)=G(\underline{u}(k\tau,x)), \; &\, x\in\Omega_{3}, \\[2mm]
\underline{v}((k\tau)^{+},x)=\underline{v}((k\tau),x), \; &\,  x\in\Omega_{3}, \\[2mm]
\underline{u}(k_{2}\tau,x)=u(k_{2}\tau,x),\; &\,x\in\overline{\Omega}_{3}, \\
\underline{v}(k_{2}\tau,x)=v(k_{2}\tau,x),\; &\,x\in\overline{\Omega}_{3}, k=k_{2},k_{2}+1,\cdots,
\end{array} \right.
\end{eqnarray*}
where $\Omega_{3}=(-S, S)$. Due to  \autoref{lemma 3-2}, it follows that
\begin{equation}\label{3-27}
\begin{aligned}
\big(\underline{u}, \underline{v}\big)\leq \big(u,v\big), ~~(t,x)\in[k_{2}\tau, \infty)\times\overline{\Omega}_{3}.
\end{aligned}
\end{equation}
By the choice of $S_{0}$, it follows from \autoref{lemma 3-1}\textcolor{blue}{$(5)$} that
\begin{equation}\label{3-28}
\lim\limits_{m\rightarrow+\infty}\big(\underline{u}(t+m\tau,x),\underline{v}(t+m\tau,x)\big)=\big(U^{S}(t,x), V^{S}(t,x)\big)
\end{equation}
uniformly on $[0,\tau]\times\overline{\Omega}_{3}$. Here, $\big(U^{S}, V^{S}\big)$ stands for the unique solution to the problem
\begin{eqnarray*}
\left\{
\begin{array}{ll}
U^{S}_{t}=d_{1} \Delta U^{S}-a_{11}U^{S}+a_{12}V^{S},\; &\, t\in(0^{+}, \tau], x\in\Omega_{3}, \\[2mm]
V^{S}_{t}=d_{2} \Delta V^{S}-a_{22}V^{S}+f(U^{S}),\; &\, t\in(0^{+}, \tau], x\in\Omega_{3}, \\[2mm]
U^{S}(t,x)=0, V^{S}(t,x)=0,\; &\, t\in[0, \tau], x\in\partial\Omega_{3}, \\[2mm]
U^{S}(0,x)=U^{S}(\tau,x), V^{S}(0,x)=V^{S}(\tau,x),\; &\,x\in\overline{\Omega}_{3},\\[2mm]
U^{S}(0^{+},x)=G\big(U^{S}(0,x)\big), V^{S}(0^{+},x)=V^{S}(0,x), \; &\, x\in\Omega_{3}.
\end{array} \right.
\end{eqnarray*}
From  \eqref{3-27} and \eqref{3-28}, it follows that
\begin{equation*}
\liminf\limits_{m\rightarrow+\infty}u(t+m\tau,x)\geq U^{S}(t,x) \text{~and~}\liminf\limits_{m\rightarrow+\infty}v(t+m\tau,x)\geq V^{S}(t,x)
\end{equation*}
uniformly on $[0,\tau]\times\overline{\Omega}_{3}$. By standard methods, one can obtain that for any $\Omega_{t}\subset\subset [0, \tau]\times\mathbb{R}$
\begin{equation*}
\lim\limits_{S\rightarrow+\infty}U^{S}(t,x)=U_{1}(t)\text{~and~}\lim\limits_{S\rightarrow+\infty}V^{S}(t,x)=V_{1}(t)
\end{equation*}
for $(t,x)\in \Omega_{t}$. Therefore, we have that for any $\Omega_{t}\subset\subset [0, \tau]\times\mathbb{R}$
\begin{equation*}
\liminf\limits_{m\rightarrow+\infty}u(t+m\tau,x)\geq U_{1}(t) \text{~and~}\liminf\limits_{m\rightarrow+\infty}v(t+m\tau,x)\geq V_{1}(t)
\end{equation*}
for $(t,x)\in \Omega_{t}$.

\textbf{(3)} \bm{$\limsup\limits_{m\rightarrow+\infty}u(t+m\tau,x)\leq U_{1}(t)$} \textbf{ and }\bm{$\limsup\limits_{m\rightarrow+\infty}v(t+m\tau,x)\leq V_{1}(t)$} \textbf{ for }  \bm{$(t,x)\in[0, \tau]\times\mathbf{R}$}

Denote the solution to the problem
\begin{eqnarray*}
\left\{
\begin{array}{ll}
\overline{u}_{t}=a_{12}\overline{v}-a_{11}\overline{u},~~~~~~~~~~~~~ t\in((k\tau)^{+}, (k+1)\tau], \\[2mm]
\overline{v}_{t}=f(\overline{u})-a_{22}\overline{v},~~~~~~~~~~~~~ t\in((k\tau)^{+}, (k+1)\tau], \\[2mm]
\overline{u}((k\tau)^{+})=G(\overline{u}(k\tau)), \overline{v}((k\tau)^{+})=\overline{v}(k\tau), \\[2mm]
\overline{u}(0)=\|u_{0}(\cdot)\|_{\mathbb{C}[-h_{0}, h_{0}]}, \overline{v}(0)=\|v_{0}(\cdot)\|_{\mathbb{C}[-h_{0}, h_{0}]}
\end{array} \right.
\end{eqnarray*}
by $(\overline{u}(t), \overline{v}(t))$. With the help of \autoref{lemma 3-2}, it follows that
\begin{equation}\label{3-29}
u(t,x)\leq \overline{u}(t)\text{~and~} v(t,x)\leq \overline{v}(t)
\end{equation}
for $t\geq 0$ and $x\in [g(t), h(t)]$.

Now, introduce the following iteration rule
\begin{eqnarray*}
\left\{
\begin{array}{ll}
\overline{u}_{t}^{(i)}+m_{1}\overline{u}^{(i)}=a_{12}\overline{v}^{(i-1)}+a_{11}\overline{u}^{(i-1)},~~~ t\in(0^{+}, \tau],\\[2mm]
\overline{v}_{t}^{(i)}+m_{2}\overline{v}^{(i)}=f(\overline{u}^{(i-1)})+a_{22}\overline{v}^{(i-1)},~~~ t\in(0^{+}, \tau], \\[2mm]
\big(\overline{u}^{(i)}(0),\overline{v}^{(i)}(0)\big)=\big(\overline{u}^{(i-1)}(\tau), \overline{v}^{(i-1)}(\tau)\big),\\[2mm]
\big(\overline{u}^{(i)}(0^{+}), \overline{v}^{(i)}(0^{+})\big)=\big(G\big(\overline{u}^{(i-1)}(\tau)\big), \overline{v}^{(i-1)}(\tau)\big), \\
\end{array} \right.
\end{eqnarray*}
where $m_{1}=2a_{11}$ and $m_{2}=2a_{22}$. Then, a monotonically decreasing iterative sequence $\big\{(\overline{u}^{(i)},
\overline{v}^{(i)})\big\}_{i=0}^{n}$ can be obtained, with initial iteration $\big\{(\overline{u}^{(0)}, \overline{v}^{(0)})\big\}=(C_{2}, C_{3})$, where $(C_{2}, C_{3})$ is defined in \autoref{lemma 3-2}.

Next we will prove that
\begin{equation}\label{3-30}
\big(\overline{u}(t+m\tau), \overline{v}(t+m\tau)\big)\leq\big(\overline{u}^{(m)}(t), \overline{v}^{(m)}(t)\big)~\text{for}~t\in[0,\tau].
\end{equation}
Noting that $(C_{2}, C_{3})$ is an upper solution of $(\overline{u}(t), \overline{v}(t))$, it follows that
\begin{equation*}
\big(\overline{u}(t), \overline{v}(t)\big)\leq(C_{2}, C_{3})=\big(\overline{u}^{(0)}(t), \overline{v}^{(0)}(t)\big)~\text{for}~t\in[0,\tau].
\end{equation*}
Assume that \eqref{3-30} holds for $m=n$, that is,
 \begin{equation}\label{3-31}
u(t+n\tau)\leq \overline{u}^{(n)}(t)\text{~and~} v(t+n\tau)\big)\leq \overline{v}^{(n)}(t)
\end{equation}
for $t\in[0,\tau]$. Now we prove that \eqref{3-30} holds for $m=n+1$. In virtue of \eqref{3-31}, one can obtain that
\begin{equation*}
u((n+1)\tau)\leq \overline{u}^{(n+1)}(0) \text{~and~} v((n+1)\tau)\leq\overline{v}^{(n+1)}(0).
\end{equation*}
This together with the assumption that $G'(u)\geq 0$ for $u\geq 0$ in (A3) yields that
\begin{equation*}
u(((n+1)\tau)^{+})\leq \overline{u}^{(n+1)}(0^{+}) \text{~and~}v(((n+1)\tau)^{+})\leq \overline{v}^{(n+1)}(0^{+}).
\end{equation*}
This combined with the classical comparison principle finally yields that
\begin{equation*}
\big(u(t+(n+1)\tau), v(t+(n+1)\tau)\big)\leq\big(\overline{u}^{(n+1)}(t), \overline{v}^{(n+1)}(t)\big)~\text{for}~t\in[0,\tau].
\end{equation*}
Therefore, \eqref{3-30} holds by mathematical induction.

Now the desired conclusion (3)  holds from \eqref{3-29}, \eqref{3-30}, and the conclusions in (1). The proof is completed.
\end{proof}
\end{theorem}
Based on \textcolor{blue}{Theorems} \ref{theorem 3-1}-\ref{theorem 3-3}, and \autoref{lemma 3-5}, the following spreading-vanishing dichotomy can be obtained.
\begin{theorem}\label{theorem 3-4}
The infected individuals and pathogenic bacteria are either spreading or vanishing.
\end{theorem}

\autoref{theorem 3-4} does not provide us with a criteria for whether the infected individuals and pathogenic bacteria are spreading or not.
This issue will be investigated in the next subsection.
\subsection{\bf Spreading-vanishing criteria}\label{Section-3.3}
With the help of $\mu_{1}$, $\mu_{2}$, $u_{0}(x)$, $v_{0}(x)$, $\lambda^{\vartriangle}_{1}(h_{0})$, and $\lambda^{\vartriangle}_{1}(\infty)$,
this subsection provides the sharp criteria for spreading and vanishing.

When $\lambda^{\vartriangle}_{1}(\infty)\geq 0$, the following theorem can be obtained directly by applying \autoref{theorem 3-1} and \autoref{lemma 3-4}.
\begin{theorem}\label{theorem 3-5}
If $\lambda^{\vartriangle}_{1}(\infty)\geq 0$, then the diseases are vanishing.
\end{theorem}
Next, when $\lambda^{\vartriangle}_{1}(\infty)<0$, we separately discuss the dynamical behaviours for the two cases: $\lambda^{\vartriangle}_{1}(h_{0})\leq 0$ and $\lambda^{\vartriangle}_{1}(h_{0})>0$.
\begin{theorem}\label{theorem 3-6}
Assume $\lambda^{\vartriangle}_{1}(\infty)<0$. If $\lambda^{\vartriangle}_{1}(h_{0})\leq 0$, then the diseases are spreading.
\begin{proof}
Arguing indirectly, assume that $-g_{\infty}<+\infty$ or $h_{\infty}<+\infty$. Then, $h_{\infty}-g_{\infty}<+\infty$ can be obtained with
the help of \autoref{lemma 3-5}. \autoref{theorem 3-2} further yields $\lambda^{\vartriangle}_{1}(G'(0),(g_{\infty},h_{\infty}))\geq0$. Since $h_{\infty}-g_{\infty}>2h_{0}$, one can obtain from \autoref{lemma 3-1}\textcolor{blue}{(1,3)} that
\begin{equation*}
0\leq \lambda^{\vartriangle}_{1}(G'(0),(g_{\infty},h_{\infty}))<\lambda^{\vartriangle}_{1}(h_{0}),
\end{equation*}
which show a contradiction with the condition  $\lambda^{\vartriangle}_{1}(h_{0})\leq 0$. This ends the proof.
\end{proof}
\end{theorem}
Before discussing the dynamical behaviours of model \eqref{Zhou-Lin} when $\lambda^{\vartriangle}_{1}(h_{0})>0$, we first give the following lemmas.
\begin{lemma}\label{lemma 3-6}
Assume $\lambda^{\vartriangle}_{1}(G'(0),(g_{\infty},h_{\infty}))< 0$, then $-g_{\infty}=h_{\infty}=\infty$.
\begin{proof}
By applying the assumption that $\lambda^{\vartriangle}_{1}(G'(0),(g_{\infty},h_{\infty}))< 0$, one can obtain that
$\lambda^{\vartriangle}_{1}(G'(0),(-\infty,+\infty))< 0$ with the help of \autoref{lemma 3-1}\textcolor{blue}{(1,3)}.
Arguing indirectly, assume that either $-g_{\infty}<\infty$ or $h_{\infty}<\infty$. Then, \autoref{lemma 3-5}
yields that $h_{\infty}-g_{\infty}<\infty$. This combined with \autoref{theorem 3-2} implies that $\lambda^{\vartriangle}_{1}\big(G'(0),(g_{\infty},h_{\infty})\big)\geq 0$, which is in contradiction to $\lambda^{\vartriangle}_{1}\big(G'(0),(g_{\infty},h_{\infty})\big)< 0$. The proof is completed.
\end{proof}
\end{lemma}
\begin{lemma}\label{lemma 3-7}
Assume that $\lambda^{\vartriangle}_{1}(\infty)< 0$ and $\lambda^{\vartriangle}_{1}(h_{0})> 0$, then a positive constant $l^{*}$ can be found such that when $h_{\infty}-g_{\infty}>l^{*}$, we have that $\lambda^{\vartriangle}_{1}\big(G'(0),(g_{\infty},h_{\infty})\big)< 0$.
\begin{proof}
By using the assumptions $\lambda^{\vartriangle}_{1}(\infty)< 0$ and $\lambda^{\vartriangle}_{1}(h_{0})> 0$, it follows from \autoref{lemma 3-1}\textcolor{blue}{(1,3)} that there exists a positive constant $l^{*}$ such that $\lambda^{\vartriangle}_{1}(G'(0),(-\frac{l^{*}}{2},\frac{l^{*}}{2}))=0$.
Then,  from \autoref{lemma 3-1}\textcolor{blue}{(1,3)} we have that $\lambda^{\vartriangle}_{1}\big(G'(0),(g_{\infty},h_{\infty})\big)< 0$ when $h_{\infty}-g_{\infty}>l^{*}$. The proof is completed.
\end{proof}
\end{lemma}
Next we will show that whether spreading or vanishing happens depends on the expansion capacities $\mu_{1}$ and $\mu_{2}$.
\begin{theorem}\label{theorem 3-7}
Assume that $\lambda^{\vartriangle}_{1}(\infty)< 0$ and $\lambda^{\vartriangle}_{1}(h_{0})> 0$, then  positive constants $\overline{\mu}$ and $\underline{\mu}$ that depend on the initial and impulsive functions can be found, such that
\begin{enumerate}
\item[$(1)$]
if $\mu_{1}+\mu_{2}\geq \overline{\mu}$, the diseases are spreading,
\item[$(2)$]
if $\mu_{1}+\mu_{2}\leq \underline{\mu}$, the diseases are vanishing.
\end{enumerate}
\begin{proof}
(1) To begin with consider the following problem
\begin{eqnarray}\label{3-32}
\left\{
\begin{array}{ll}
\underline{u}_{t}=d_{1}\Delta \underline{u}-a\underline{u},\; &\, (t,x)\in((k\tau)^{+}, (k+1)\tau]\times(\underline{g}(t),\underline{h}(t)), \\[2mm]
\underline{v}_{t}=d_{2}\Delta \underline{v}-a\underline{v},\; &\, (t,x)\in((k\tau)^{+}, (k+1)\tau]\times(\underline{g}(t),\underline{h}(t)), \\[2mm]
\underline{u}=0, \underline{v}=0,\; &\, (t,x)\in(k\tau, (k+1)\tau]\times\{\underline{g}(t),\underline{h}(t)\},\\[2mm]
\underline{u}((k\tau)^{+},x)=G(\underline{u}(k\tau,x)), \; &\, x\in(\underline{g}(t),\underline{h}(t)), \\[2mm]
\underline{v}((k\tau)^{+},x)=\underline{v}(k\tau,x), \; &\,  x\in(\underline{g}(t),\underline{h}(t)), \\[2mm]
\underline{g}'(t)=-\mu_{1}\underline{u}_{x}(t,\underline{g}(t))-\mu_{2}\underline{v}_{x}(t,\underline{g}(t)),\; &\, t\in(k\tau, (k+1)\tau],\\[2mm]
\underline{h}'(t)=-\mu_{1}\underline{u}_{x}(t,\underline{h}(t))-\mu_{2}\underline{v}_{x}(t,\underline{h}(t)),\; &\, t\in(k\tau, (k+1)\tau],\\[2mm]
\underline{g}(0)=-h_{0}, \underline{u}(0,x)=u_{0}(x),\; &\,x\in[-h_{0},h_{0}],\\[2mm]
\underline{h}(0)=h_{0},~~~ \underline{v}(0,x)=v_{0}(x),\; &\,x\in[-h_{0},h_{0}], ~k=0,1,2,3,4,\cdots,
\end{array} \right.
\end{eqnarray}Problem 3.2 has a unique global nonnegative classical solution,
where $a:=\max\{a_{11}, a_{22}\}$. By using a similar proof method as in proving \autoref{lemma 2-2} and \autoref{theorem 3-1},
it can be obtained that problem \eqref{3-32} possesses a unique global nonnegative classical solution $(\overline{u}, \overline{v}, \overline{g}, \overline{h})$, where $-\overline{g}(t), \overline{h}(t)>0$ for $t\in\mathbb{R}^{+}$. Due to  \autoref{Remark 3-1}, it follows that
\begin{equation}\label{3-33}
\big(\underline{u}, \underline{v}\big)\leq \big(u,v\big)\text{~and~}\big(\underline{h}, -\underline{g}\big)\leq \big(h,-g\big)
\end{equation}
for $(t,x)\in [0, +\infty)\times [\underline{g}(t), \underline{h}(t)]$.

Consider the following linear parabolic problem
\begin{eqnarray*}
\left\{
\begin{array}{ll}
\zeta_{t}=d_{1}\Delta \zeta-a\zeta,\; &\, t\in(0^{+}, \tau], x\in(\phi(t),\varphi(t)), \\[2mm]
\eta_{t}=d_{2}\Delta \eta-a\eta,\; &\, t\in(0^{+}, \tau], x\in(\phi(t),\varphi(t)), \\[2mm]
\zeta=0, \eta=0,\; &\, t\in(0, \tau], ~~x\in\{\phi(t),\varphi(t)\},\\[2mm]
\zeta(0^{+},x)=G(\zeta(0,x)), \; &\, x\in(\phi(t),\varphi(t)), \\[2mm]
\eta(0^{+},x)=\eta(0,x), \; &\,  x\in(\phi(t),\varphi(t)), \\[2mm]
\phi'(t)=-\mu_{1}\zeta_{x}(t,\phi(t))-\mu_{2}\eta_{x}(t,\phi(t)),\; &\, t\in(0, \tau],\\[2mm]
\varphi'(t)=-\mu_{1}\zeta_{x}(t,\varphi(t))-\mu_{2}\eta_{x}(t,\varphi(t)),\; &\, t\in(0, \tau],\\[2mm]
\phi(0)=-\frac{h_{0}}{2}, \zeta(0,x)=\zeta_{0}(x),\; &\,x\in[-\frac{h_{0}}{2},\frac{h_{0}}{2}],\\[2mm]
\varphi(0)=\frac{h_{0}}{2},~~ \eta(0,x)=\eta_{0}(x),\; &\,x\in[-\frac{h_{0}}{2},\frac{h_{0}}{2}],
\end{array} \right.
\end{eqnarray*}
where the initial functions satisfy
\begin{eqnarray*}
\left\{
\begin{array}{ll}
\zeta_{0}(\pm\frac{h_{0}}{2})=0, ~\zeta'_{0}(-\frac{h_{0}}{2})>0, ~\zeta'_{0}(\frac{h_{0}}{2})<0,\; &\ \\[2mm]
\eta_{0}(\pm\frac{h_{0}}{2})=0, ~\eta'_{0}(-\frac{h_{0}}{2})>0, ~\eta'_{0}(\frac{h_{0}}{2})<0,\; &\ \\[2mm]
(0,0)<(\zeta_{0}(x),\eta_{0}(x))\leq (u_{0}(x), v_{0}(x)),\; &\
\end{array} \right.
\end{eqnarray*}
for $x\in(-\frac{h_{0}}{2},\frac{h_{0}}{2})$. Then, using the Hopf boundary lemma for cooperative systems (see \cite[Theorem 14 in Chapter 3]{Protter-Weinberger}),
one can obtain that
\begin{equation*}
\zeta_{x}(t,\phi(t)), \eta_{x}(t, \phi(t)), -\zeta_{x}(t,\varphi(t)), -\eta_{x}(t,\varphi(t))<0
\end{equation*}
for $t\in(0,\tau]$. Therefore, a positive constant $\overline{\mu}$ can be found such that when $\mu_{1}+\mu_{2}\geq \overline{\mu}$, we have that
\begin{equation*}
\phi(\tau)\leq -\frac{l^{*}}{2}\text{~and~}\varphi(\tau)\geq \frac{l^{*}}{2},
\end{equation*}
where $l^{*}$ is defined in \autoref{lemma 3-7}. When $\mu_{1}+\mu_{2}\geq \overline{\mu}$, it follows from the comparison principle that
\begin{equation*}
g(\tau)\leq -\frac{l^{*}}{2}\text{~and~}h(\tau)\geq \frac{l^{*}}{2}.
\end{equation*}
This means that $h_{\infty}-g_{\infty}>l^{*}$. By using \autoref{lemma 3-7}, one can obtain that $\lambda^{\vartriangle}_{1}\big(G'(0),(g_{\infty},h_{\infty})\big)< 0$. This combined with \autoref{lemma 3-6} yields that $h_{\infty}=-g_{\infty}=\infty$. Finally, \autoref{theorem 3-3} gives that the diseases are spreading.

(2) Because the procedure of the proof resembles \autoref{theorem 3-1}, only a sketch is given here. First construct the same functions as in \eqref{3-4}. Then, select the corresponding coefficients. Select
\begin{equation*}
h:=h_{0} \text{~and~}\gamma:=\frac{\lambda^{\vartriangle}_{1}(G'(0),(-h_{0}, h_{0}))}{2},
\end{equation*}
where the positive constant $\sigma$ can be chosen such that
\begin{equation*}
1-\kappa^{-2}<\frac{\lambda^{\vartriangle}_{1}(G'(0),(-h_{0}, h_{0}))}{2\lambda_{0}(h_{0})\max\big\{d_{1}, d_{2}\big\}},
\end{equation*}
and  the sufficiently large $M$ can be chosen such that
\begin{equation*}
\big(\overline{u}(0, x), \overline{v}(0, x)\big)\geq \big(u_{0}(x), v_{0}(x) \big).
\end{equation*}
Next, when
\begin{equation*}
\mu_{1}+\mu_{2}\leq \underline{\mu}:=\frac{\gamma \sigma(2+\sigma)h^{2}_{0} }{2\pi M\min\{\alpha_{2}, \beta_{2}\}},
\end{equation*}
the constructed quadruple $(\overline{u},\overline{v},-\Gamma(t),\Gamma(t))$ is a super-solution of model \eqref{Zhou-Lin}. Finally, that
the diseases are vanishing can be obtained by the same proof as in \autoref{theorem 3-1}. For more details, the reader can refer
to the proof of \autoref{theorem 3-1}. This proof is completed.
\end{proof}
\end{theorem}
Based on \autoref{theorem 3-7}, we can now give the following sharp criteria governing the spreading and vanishing of the diseases when $-\lambda^{\vartriangle}_{1}(h_{0}), \lambda^{\vartriangle}_{1}(\infty)< 0$.
\begin{theorem}\label{theorem 3-8}
Assume $\lambda^{\vartriangle}_{1}(\infty)< 0$ and $\lambda^{\vartriangle}_{1}(h_{0})>0$. For any given positive constant $\mu_{1}$,
there exists a positive constant $\mu_{0}$ that depends on the initial and impulse functions such that spreading occurs as $\mu_{2}>\mu_{0}$,
and vanishing occurs as $0<\mu_{2}\leq\mu_{0}$.
\begin{proof}
The proof process is the same as in \cite[Theorem 3.10]{Wang-Du}, we omit it here.
\end{proof}
\end{theorem}
The criteria provided in \autoref{theorem 3-8} are governed by the expansion capacity $\mu_{1}$ of pathogenic bacteria and the expansion capacity $\mu_{2}$ of the infected individuals. Next, we will deduce  criteria governed by the initial functions $u_{0}(x)$ and $v_{0}(x)$. As a start we prove the following lemma, which has a significant role in the construction of the lower solution of model \eqref{Zhou-Lin}.
\begin{lemma}\label{lemma 3-8}
Let $\mu_{0}$ be the principal eigenvalue of
\begin{eqnarray}\label{lei-1}
\left\{
\begin{array}{ll}
d\Delta\varphi(x)+\frac{1}{2}\varphi'(x)+\mu\varphi(x)=0,\; &\, 0<x<1,\\[2mm]
\varphi'(0)=\varphi(1)=0,\; &\,
\end{array} \right.
\end{eqnarray}
then there exists a principal eigenfunction $\varphi_{0}(x)\in\mathbb{C}^{\infty}([0,1])$,  satisfying that $\varphi_{0}(x)>0$ and $\varphi'_{0}(x)<0$ in $(0,1)$, and $\|\varphi_{0}(\cdot)\|_{L^{\infty}[0,1]}=1$.

\begin{proof}
Let
\begin{equation*}
\alpha:=-\frac{1}{4d}, \beta:=\frac{\sqrt{4d\mu-1/4}}{2d}~\text{and}~\beta_{0}:=\frac{\sqrt{4d\mu_{0}-1/4}}{2d}.
\end{equation*}
From standard eigenvalue and eigenfunction theory, the principal eigenvalue $\mu_{0}$ of problem \eqref{lei-1} is the minimum solution of
\begin{equation*}
\tan \beta=\beta/\alpha,
\end{equation*}
and the principal eigenfunction can be represented as
\begin{equation}\label{lei-2}
\varphi(x)=-\alpha e^{\alpha x}\sin(\beta_{0}x)+\beta_{0}e^{\alpha x}\cos(\beta_{0}x).
\end{equation}
By standard computation, \eqref{lei-2} can be further simplified into
\begin{equation}\label{lei-3}
\begin{aligned}
\varphi(x)&=e^{\alpha x}\sqrt{\alpha ^{2}+\beta_{0}^{2}}\sin\Big[\beta_{0}x+\arctan(\frac{\beta_{0}}{-\alpha})\Big]\\
&=e^{\alpha x}\sqrt{\alpha ^{2}+\beta_{0}^{2}}\sin\Big[\beta_{0}x+\pi-\beta_{0}\Big]\\
&=-e^{\alpha x}\sqrt{\alpha ^{2}+\beta_{0}^{2}}\sin\big[\beta_{0}(x-1)\big].
\end{aligned}
\end{equation}
Noting that
\begin{equation*}
-\pi<-\beta_{0}<\beta_{0}(x-1)<0,
\end{equation*}
it follows that $\varphi_{0}(x)>0$ in $(0,1)$.

From  \eqref{lei-3} we find that
\begin{equation*}
\begin{aligned}
\varphi'(x)&=-\alpha e^{\alpha x}\sqrt{\alpha ^{2}+\beta_{0}^{2}}\sin\big[\beta_{0}(x-1)\big]-\beta_{0}e^{\alpha x}\sqrt{\alpha ^{2}+\beta_{0}^{2}}\cos\big[\beta_{0}(x-1)\big].\\
&=-e^{\alpha x}\sqrt{\alpha ^{2}+\beta_{0}^{2}}\Big[\alpha\sin\big(\beta_{0}(x-1)\big)+\beta_{0}\cos\big(\beta_{0}(x-1)\big)\Big]\\
&=-e^{\alpha x}(\alpha ^{2}+\beta_{0}^{2})\cos\Big[\beta_{0}(x-1)-\arctan(\frac{\alpha}{\beta_{0}})\Big]\\
&=-e^{\alpha x}(\alpha ^{2}+\beta_{0}^{2})\cos\Big[\beta_{0}(x-1)+\frac{\pi}{2}-\arctan(\frac{\beta_{0}}{-\alpha})\Big]\\
&=-e^{\alpha x}(\alpha ^{2}+\beta_{0}^{2})\cos\Big[\beta_{0}(x-1)+\frac{\pi}{2}-(\pi-\beta_{0})\Big]\\
&=-e^{\alpha x}(\alpha ^{2}+\beta_{0}^{2})\cos\Big[\frac{\pi}{2}-\beta_{0}x\Big]\\
&=-e^{\alpha x}(\alpha ^{2}+\beta_{0}^{2})\sin(\beta_{0}x).\\
\end{aligned}
\end{equation*}
Noting that $0<\beta_{0}x<\beta_{0}<\pi$, it follows that $\varphi'_{0}(x)<0$ in $(0,1)$. Finally, let
\begin{equation*}
\varphi_{0}(x)=\frac{\varphi(x)}{\|\varphi(\cdot)\|_{L^{\infty}[0,1]}}.
\end{equation*}
The ends this proof.
\end{proof}
\end{lemma}
We are now ready to provide the criteria governed by the initial functions $u_{0}(x)$ and $v_{0}(x)$ by using \autoref{lemma 3-8}.
\begin{theorem}\label{theorem 3-9}
Assume $\lambda^{\vartriangle}_{1}(\infty)< 0$ and $\lambda^{\vartriangle}_{1}(h_{0})>0$. Then, the following statements are valid.
\begin{enumerate}
\item[$(1)$]
If $u_{0}(x)$ and $v_{0}(x)$ are sufficiently small, the diseases are vanishing,
\item[$(2)$]
If $u_{0}(x)$ and $v_{0}(x)$ are sufficiently large and the impulsive function is linear, the diseases are spreading.
\end{enumerate}
\begin{proof}
(1) To begin with, construct the following functions
\begin{eqnarray*}
\begin{array}{l}
~~~~~~~~~~~~ ~~~~~\Gamma(t):=h_{0}\kappa(t),~\kappa(t):=1+\sigma-\frac{\sigma}{2}e^{-\gamma t}, ~~~~t\in[0,+\infty), \\[2mm]
\overline{u}(t,x):=Me^{-\gamma t}\phi(t, \frac{h_{0}x}{\Gamma(t)}), \overline{v}(t,x):=Me^{-\gamma t}\psi(t,\frac{h_{0}x}{\Gamma(t)}),~(t,x)\in[0,+\infty)\times[-\Gamma(t), \Gamma(t)],
\end{array}
\end{eqnarray*}
where the values of the positive constants $\gamma$, $\sigma$, and $M$ are chosen subsequently, and $(\phi, \psi)$ is a strongly positive principal eigenfunction pair for the periodic eigenvalue problem \eqref{3-2} with $(s_{1}, s_{2})=(-h_{0}, h_{0})$. Using the same method as in the proof of \autoref{theorem 3-1}, $\gamma$ and $\sigma$ can be found such that
\begin{eqnarray*}
\left\{
\begin{array}{ll}
\overline{u}_{t}=d_{1}\Delta \overline{u}-a_{11}\overline{u}+a_{12}\overline{v},\; &\, t\in((k\tau)^{+}, (k+1)\tau], x\in(-\Gamma(t),\Gamma(t)), \\[2mm]
\overline{v}_{t}=d_{2}\Delta \overline{v}-a_{22}\overline{v}+f(\overline{u}),\; &\, t\in((k\tau)^{+}, (k+1)\tau], x\in(-\Gamma(t),\Gamma(t)), \\[2mm]
\overline{u}=0, \overline{v}=0,\; &\, t\in(k\tau, (k+1)\tau], ~~~~~x\in\{-\Gamma(t), \Gamma(t)\},\\[2mm]
\overline{u}((k\tau)^{+},x)=G(\overline{u}(k\tau,x)), \; &\, x\in(-\Gamma(k\tau), \Gamma(k\tau)), \\[2mm]
\overline{v}((k\tau)^{+},x)=\overline{v}(k\tau,x), \; &\,  x\in(-\Gamma(k\tau), \Gamma(k\tau)), ~k=0,1,2,\cdots.
\end{array} \right.
\end{eqnarray*}
Next, choose
\begin{equation*}
M:=\frac{\gamma\sigma h^{2}_{0}(2+\sigma)}{2\pi(\mu_{1}\alpha_{2}+\mu_{2}\beta_{2})},
\end{equation*}
where $\alpha_{2}$ and $\beta_{2}$ are defined in \autoref{lemma 3-3}. Then careful calculations yield
\begin{eqnarray*}
\left\{
\begin{array}{ll}
\Gamma'(t)\geq\mu_{1}\overline{u}_{x}(t,-\Gamma(t))+\mu_{2}\overline{v}_{x}(t,-\Gamma(t)),\; &\, t\in\mathbb{R}^{+},\\[2mm]
\Gamma'(t)\geq-\mu_{1}\overline{u}_{x}(t,\Gamma(t))-\mu_{2}\overline{v}_{x}(t,\Gamma(t)),\; &\, t\in\mathbb{R}^{+}.
\end{array} \right.
\end{eqnarray*}
At the end, take sufficiently small initial function $(u_{0}(x), v_{0}(x))$ such that
\begin{equation*}
\big(u_{0}(y), v_{0}(y)\big)\leq \big(\overline{u}(0,x) , \overline{v}(0,x)\big)\text{~for~}x\in[-\Gamma(0), \Gamma(0)],
\end{equation*}
where $y=\frac{h_{0}x}{\Gamma(0)}$. Therefore, the quadruple $(\overline{u},\overline{v},-\Gamma(t),\Gamma(t))$ is a super-solution of model \eqref{Zhou-Lin}. The remaining of the proof resembles \autoref{theorem 3-1} and therefore is omitted here.

(2) Let $\big(\mu_{1}, \varphi_{1}\big)$ and $\big(\mu_{2}, \varphi_{2}\big)$ be the principal eigenpairs of problem \eqref{lei-1}
with diffusion coefficients $d_{1}$ and $d_{2}$, respectively, and $(\Phi(t), \Psi(t))$ be the unique positive solution of
\eqref{SFC-8} with $\lambda$ replaced by $\min\{\mu_{1},\mu_{2}\}$. First of all, construct the following functions
\begin{eqnarray*}
\begin{array}{l}
~~~~~~~~~~~~~\underline{h}(t):=\sqrt{t+\delta},~~~~~~~~~t\geq 0,\\
\underline{u}(t, x):=\frac{M}{(t+\delta)^{k}}\Phi(t)\varphi_{1}(\frac{x}{\sqrt{t+\delta}}),~~t\geq 0, ~x\in\big[0, \underline{h}(t)\big],\\
\underline{v}(t, x):=\frac{M}{(t+\delta)^{k}}\Psi(t)\varphi_{2}(\frac{x}{\sqrt{t+\delta}}),~~t\geq 0, ~x\in\big[0, \underline{h}(t)\big],\\
\end{array}
\end{eqnarray*}
where $\delta$, $k$, and $M$ are positive constants to be chosen later. By the assumptions $\lambda^{\vartriangle}_{1}(\infty)< 0$ and $\lambda^{\vartriangle}_{1}(h_{0})>0$, it follows from \autoref{lemma 3-1}\textcolor{blue}{(1,3)} that there exists a positive integer $k_{0}$ such that
\begin{equation*}
\lambda^{\vartriangle}_{1}\big(G'(0),\big(-\underline{h}(k_{0}\tau), \underline{h}(k_{0}\tau)\big)\big)<0.
\end{equation*}
Since $u_{0}(x)$ and $v_{0}(x)$ are sufficiently large, there exist sufficiently large functions $\overline{u}_{0}(x)$ and $\overline{v}_{0}(x)$ such that
\begin{equation*}
\big(\overline{u}_{0}(x), \overline{v}_{0}(x)\big)\leq \big(u_{0}(x),  v_{0}(x)\big)\text{~and~}\big(\overline{u}_{0}(x),\overline{v}_{0}(x)\big) =\big(\overline{u}_{0}(-x), \overline{v}_{0}(-x)\big).
\end{equation*}
for $x\in[-h_{0},h_{0}]$.

Next, we proceed to show that the triple $(\underline{u},\underline{v}, \underline{h}(t))$ is a lower solution of model \eqref{Zhou-Lin} in half line $(0, h(t))$.
 A calculation yields that, for $(t,x)\in(0, k_{0}\tau]\times\big[0, \underline{h}(t)\big]$,
\begin{equation*}
\begin{aligned}
&\underline{u}_{t}-d_{1}\Delta \underline{u}+a_{11}\underline{u}-a_{12}\underline{v}\\
\leq&~-\frac{M}{(t+\delta)^{k+1}}\Big[\big(k-a_{11}(t+\delta))\Phi(t)\varphi_{1}+\frac{x}{2\sqrt{t+\delta}}\Phi(t)\varphi'_{1}
-(t+\delta)\Phi_{t}(t)\varphi_{1}+d_{1}\Phi(t)\varphi''_{1}\Big] \\
\leq&~-\frac{M}{(t+\delta)^{k+1}}\Big[k\Phi(t)\varphi_{1}+\frac{1}{2}\Phi(t)\varphi'_{1}-a_{12}(t+\delta)\Psi(t)\varphi_{1}
+d_{1}\Phi(t)\varphi''_{1}+(\lambda-d_{1}\lambda_{0})(t+\delta)\Phi(t)\varphi_{1}\Big] \\
\leq&~-\frac{M\Phi(t)\varphi_{1}}{(t+\delta)^{k+1}}\Big[k-\mu_{1}-a_{12}C(k_{0}\tau+\delta)+(\lambda-d_{1}\lambda_{0})(t+\delta)\Big]\leq 0
\end{aligned}
\end{equation*}
and
\begin{equation*}
\begin{aligned}
&\underline{v}_{t}-d_{2}\Delta \underline{v}+a_{22}\underline{v}-f(\underline{u})\\
\leq&~-\frac{M}{(t+\delta)^{k+1}}\Big[\big(k-a_{22}(t+\delta))\Psi(t)\varphi_{2}+\frac{x}{2\sqrt{t+\delta}}\Psi(t)\varphi'_{2}
-(t+\delta)\Psi_{t}(t)\varphi_{2}+d_{2}\Psi(t)\varphi''_{2}\Big] \\
\leq&~-\frac{M}{(t+\delta)^{k+1}}\Big[k\Psi(t)\varphi_{2}+\frac{1}{2}\Psi(t)\varphi'_{2}-f'(0)(t+\delta)\Phi(t)\varphi_{2}
+d_{2}\Psi(t)\varphi''_{2}+(\lambda-d_{2}\lambda_{0})(t+\delta)\Psi(t)\varphi_{2}\Big] \\
\leq&~-\frac{M\Psi(t)\varphi_{2}}{(t+\delta)^{k+1}}\Big[k-\mu_{2}-\frac{a_{12}(k_{0}\tau+\delta)}{C}+(\lambda-d_{2}\lambda_{0})(t+\delta)\Big]\leq 0
\end{aligned}
\end{equation*}
provided that $0<\delta\leq\max\{1, h^{2}_{0}\}$ and
\begin{equation*}
k>\max\limits_{t\in[0,k_{0}\tau]}\Big\{\mu_{1}+a_{12}C(k_{0}\tau+\delta)-(\lambda-d_{1}\lambda_{0})(t+\delta), \mu_{2}+\frac{a_{12}(k_{0}\tau+\delta)}{C}-(\lambda-d_{2}\lambda_{0})(t+\delta)\Big\}
\end{equation*}
where $\lambda$, $\lambda_{0}$ and $C$ are defined in \eqref{3-2}, \eqref{SFC-8} and \eqref{3-13}, respectively. For the boundary condition,
it is easy to obtain that
\begin{equation*}
\underline{u}_{x}(t, 0)=\underline{v}_{x}(t, 0)=0 \text{~and~}\underline{u}(t,\underline{h}(t))=\underline{v}(t, \underline{h}(t))=0\text{~for~}t>0.
\end{equation*}
From the assumption that the impulsive function is linear, one can obtain that
\begin{equation*}
\begin{aligned}
&\underline{u}\big((k\tau)^{+}, x\big)- G\big(\underline{u}(k\tau,x)\big)=\frac{M}{(k\tau+\delta)^{k}}\varphi_{1}\big(\frac{x}{\sqrt{k\tau+\delta}}\big)\Big[\Phi((k\tau)^{+})-G'(0)\Phi(k\tau)\Big] \\
\leq&\frac{M}{(k\tau+\delta)^{k}}\varphi_{1}\big(\frac{x}{\sqrt{k\tau+\delta}}\big)\Phi(0)\Big[G'(0)-G'(0)\Big]\leq 0 \\
\end{aligned}
\end{equation*}
and
\begin{equation*}
\begin{aligned}
&\underline{v}\big((k\tau)^{+}, x\big)-\underline{v}(k\tau,x)=\frac{M}{(k\tau+\delta)^{k}}\varphi_{2}\big(\frac{x}{\sqrt{k\tau+\delta}}\big)\Big[\Psi((k\tau)^{+})-\Psi(k\tau)\Big] \\
\leq&\frac{M}{(k\tau+\delta)^{k}}\varphi_{1}\big(\frac{x}{\sqrt{k\tau+\delta}}\big)\Big[\Psi((0)^{+})-\Psi(0)\Big]\leq 0 \\
\end{aligned}
\end{equation*}
hold for $x\in[0, \underline{h}(t)]$. The sufficiently large $M$ can be chosen such that
\begin{equation*}
-\frac{\mu_{1}\varphi'_{1}(1)\min\limits_{t\in[0,\tau]}\Phi(t)+\mu_{2}\varphi'_{2}(1)\min\limits_{t\in[0,\tau]}\Psi(t)}{(k_{0}\tau+1)^{k}}\geq \frac{1}{2}.
\end{equation*}
This implies that
\begin{equation*}
\begin{aligned}
&\underline{h}'(t)+\mu_{1}\underline{u}_{x}(t,\underline{h}(t))-\mu_{2}\underline{v}_{x}(t,\underline{h}(t))\\
\leq& \frac{1}{2\sqrt{t+\delta}}+\frac{\mu_{1}M\varphi'_{1}(1)\Phi(t)}{(t+\delta)^{k+1/2}}+\frac{\mu_{2}M\varphi'_{2}(1)\Psi(t)}{(t+\delta)^{k+1/2}} \\
\leq&\frac{1}{\sqrt{t+\delta}}\Big[\frac{1}{2}+\frac{\mu_{1}\varphi'_{1}(1)\Phi(t)+\mu_{2}\varphi'_{2}(1)\Psi(t)}{(k_{0}\tau+1)^{k}}M\Big] \\
\leq&0.
\end{aligned}
\end{equation*}
Finally, if $\underline{u}(x,0)\leq \frac{M}{\delta^{k}}\Phi(0)\varphi_{1}(\frac{x}{\sqrt{\delta}})\leq \overline{u}_{0}(x)$ and $\underline{v}(x,0)\leq \frac{M}{\delta^{k}}\Psi(0)\varphi_{2}(\frac{x}{\sqrt{\delta}})\leq \overline{v}_{0}(x)$ in $[0, \sqrt{\delta}]$, then the triple $(\underline{u},\underline{v}, \underline{h}(t))$ is a lower solution of model \eqref{Zhou-Lin} on the half line $[0, h(t)]$.

By the analysis above, there exists $z,w\in\big[\mathbb{PC}^{1,2}_{t,x}\cap \mathbb{PC}_{t,x}\big]^{2}$ such that the quadruple $(z,w, \underline{h}(t),-\underline{h}(t))$ is a lower solution of model \eqref{Zhou-Lin} in $[-h(t), h(t)]$. Now, from  \autoref{Remark 3-1}  it follows that $h(k_{0}\tau)\geq \underline{h}(k_{0}\tau)$. This together with \autoref{lemma 3-1}\textcolor{blue}{(1,3)} implies that
\begin{equation*}
\begin{aligned}
&\lambda^{\vartriangle}_{1}\big(G'(0),(g_{\infty},h_{\infty})\big)<\lambda^{\vartriangle}_{1}\big(G'(0),(g(k_{0}\tau),h(k_{0}\tau)\big)\\
\leq& \lambda^{\vartriangle}_{1}\big(G'(0),(-\underline{h}(k_{0}\tau),\underline{h}(k_{0}\tau)\big)\leq0.
\end{aligned}
\end{equation*}
Then, $h_{\infty}=-g_{\infty}=+\infty$ can be obtained by using \autoref{lemma 3-6}. With the help of \autoref{theorem 3-3} we conclude that the diseases are spreading.
This proof is completed.
\end{proof}
\end{theorem}
Based on \autoref{theorem 3-9}, we now can give the following sharp criteria governing the spreading and vanishing of the diseases when $-\lambda^{\vartriangle}_{1}(G'(0),(-h_{0},h_{0})), \lambda^{\vartriangle}_{1}(G'(0),(-\infty,+\infty))< 0$.
\begin{theorem}\label{theorem 3-10}
Suppose that the impulsive function is linear, that $\lambda^{\vartriangle}_{1}(\infty)< 0$ and $\lambda^{\vartriangle}_{1}(h_{0})>0$, and that $(\upsilon(x),v_{0}(x))$ satisfies Assumption (A1). Let $(u,v,g,h)$ be a solution of model \eqref{Zhou-Lin} with $(u_{0}(x),v_{0}(x))=(\kappa\upsilon(x),v_{0}(x))$ for some $\kappa>0$. Then, for any given $v_{0}(x)$, there exists $\kappa_{0}\in(0,+\infty)$ such that spreading occurs when $\kappa>\kappa_{0}$,
and vanishing occurs when $0<\kappa\leq\kappa_{0}$.
\begin{proof}
The proof process is the same as \cite[Theorem 3.5]{Li-Lin}, we omit it here.
\end{proof}
\end{theorem}
\section{\bf Numerical simulation}\label{Section-4}
The theoretical analyses have been completed in the previous sections. In this section, two numerical examples will be given to demonstrate the theoretical results and more specifically show the influence of the implementation of impulse intervention and the movement of infected region on the spread of the diseases.

Noting that \eqref{Zhou-Lin} is a cooperative model, this section only simulates the evolution of pathogenic bacteria with time and space in order to save space. In order to judge the sign of the principal eigenvalue conveniently, we regard the simulated $h(t)$ as the real one. In all simulations, we choose the growth function to be Beverton-Holt functions, that is
\begin{equation*}
f(u)=\frac{mu}{a+u},
\end{equation*}
and the initial infected environment is chosen to be $[-h_{0},h_{0}]=[-2,2]$. Additionally, the initial functions are as follows:
\begin{equation*}
u_{0}(x)=0.3\cos\Big(\frac{\pi x}{4}\Big),~v_{0}(x)=0.1\cos\Big(\frac{\pi x}{4}\Big), ~x\in[-2, 2].
\end{equation*}
Other parameters of model \eqref{Zhou-Lin} are given in each of the following subsections to show the different dynamical behaviours.
\subsection{The influence of impulse intervention}
To observe the impact of pulse intervention on the development of the diseases, this subsection performs numerical simulations when the expansion capabilities of pathogenic bacteria and infected individuals are unchanged.
\begin{exm}\label{exm1}
Fix $d_{1}=0.1$, $d_{2}=0.4$, $a_{11}=0.3$, $a_{12}=0.5$,  $a_{22}=0.1$, $m=1$, $a=10$, $\mu_{1}=10$, $\mu_{2}=15$, and $\tau=5$.
The periodic impulse function $g$ is taken to be $u$ and $\frac{0.5u}{10+u}$, respectively.
\end{exm}
\begin{figure}[!ht]
\centering
\subfigure{ {
\includegraphics[width=0.40\textwidth]{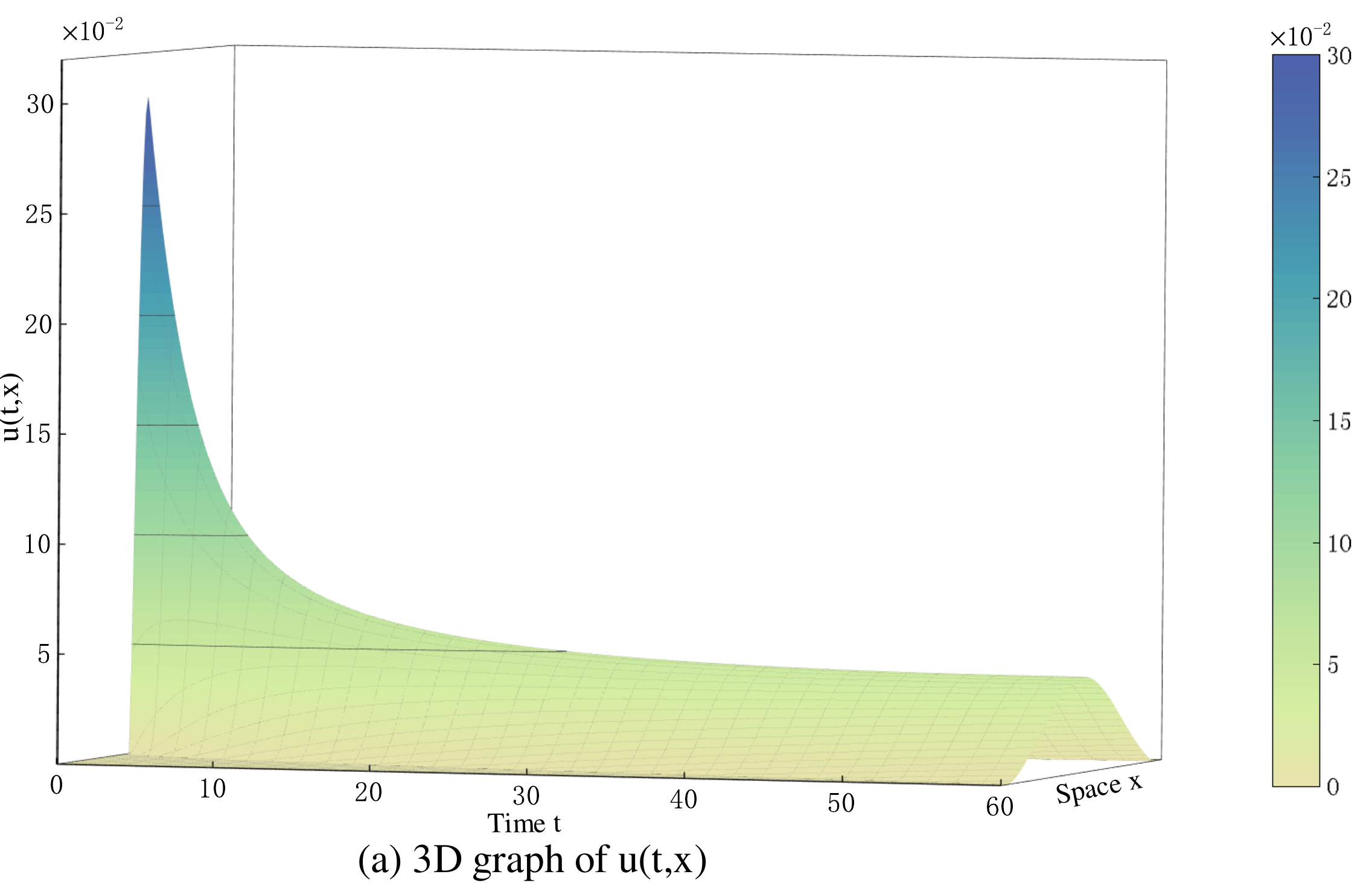}
} }
\subfigure{ {
\includegraphics[width=0.40\textwidth]{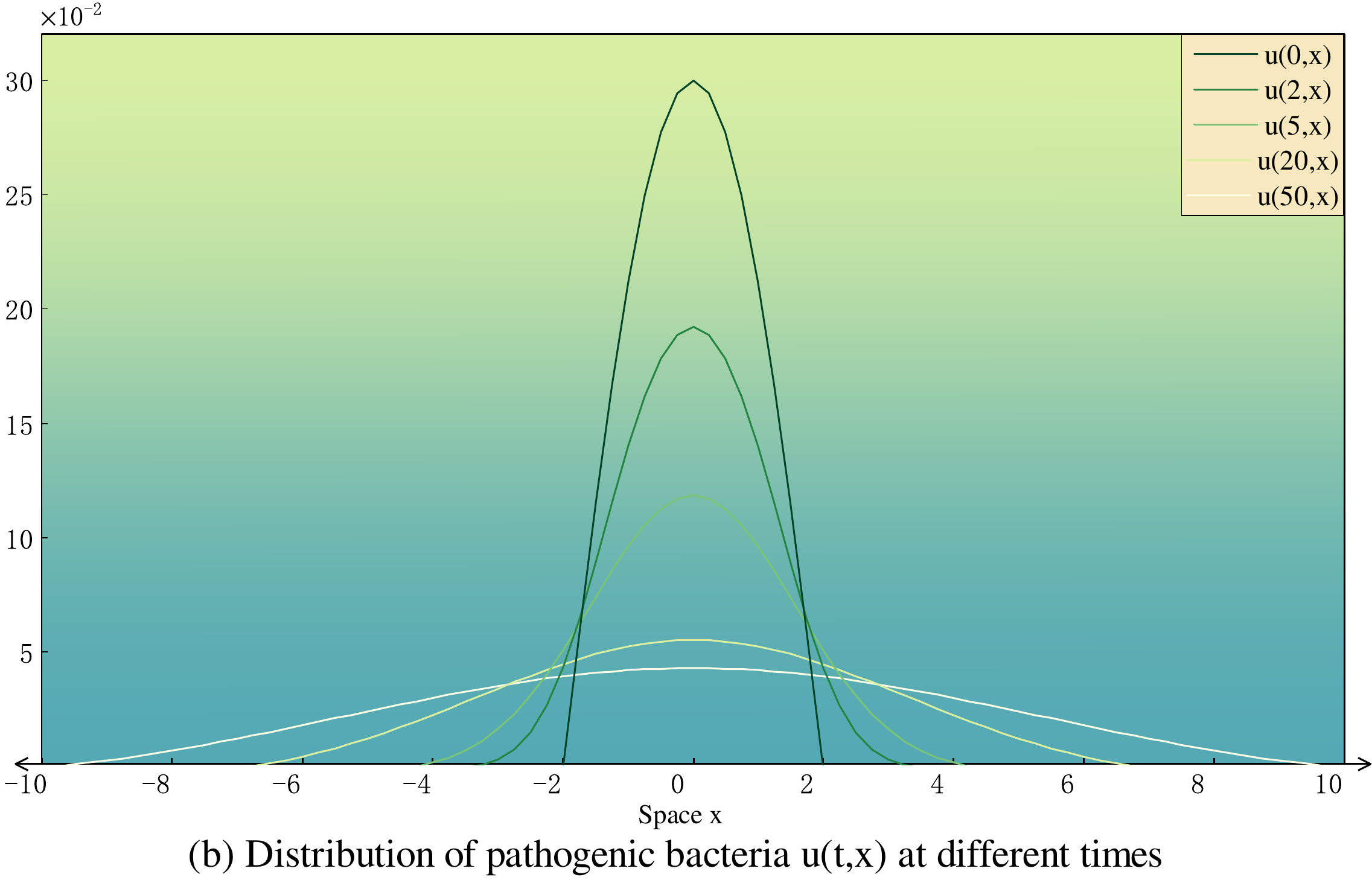}
} }
\subfigure{ {
\includegraphics[width=0.40\textwidth]{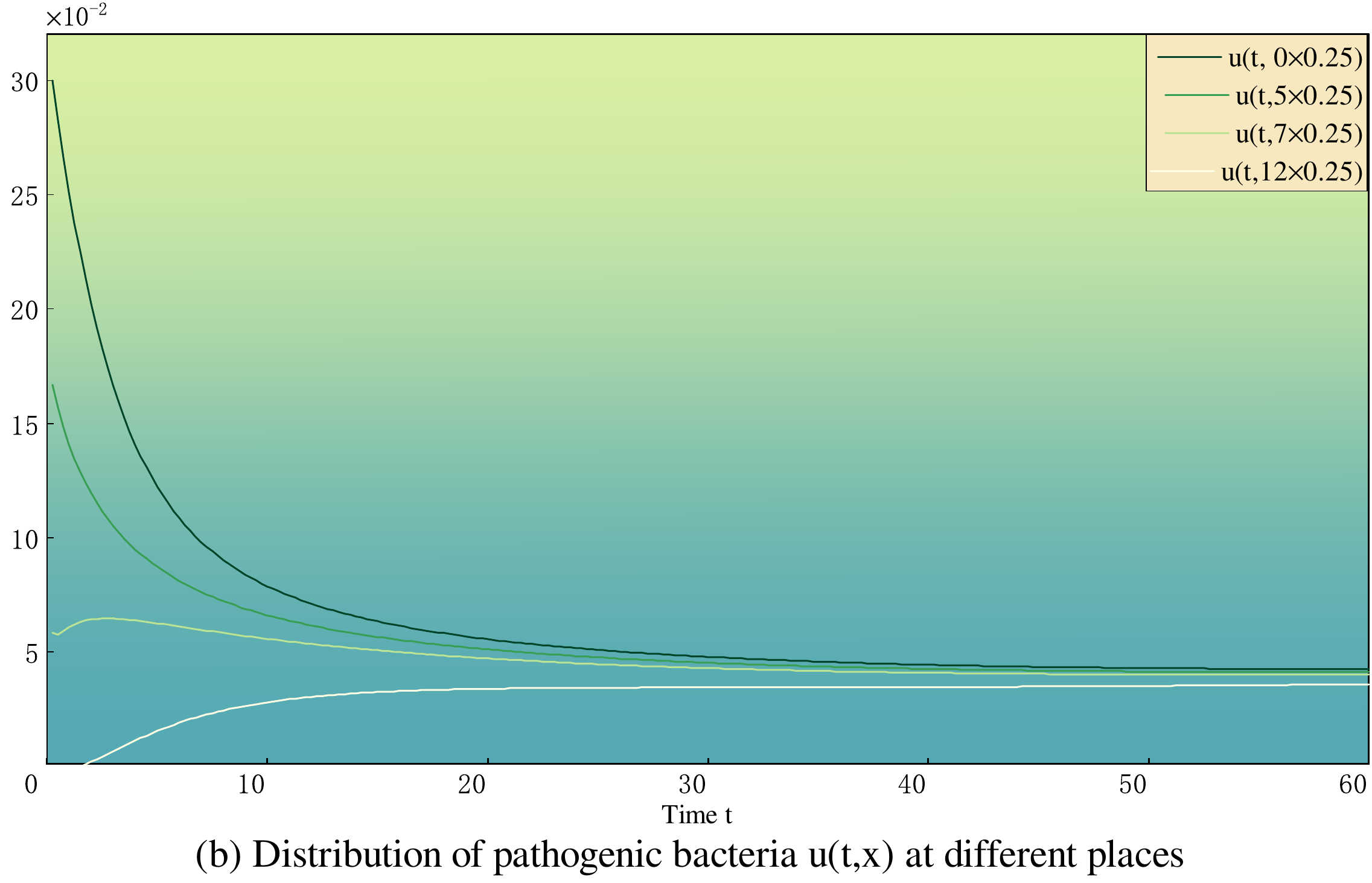}
} }
\subfigure{ {
\includegraphics[width=0.4\textwidth]{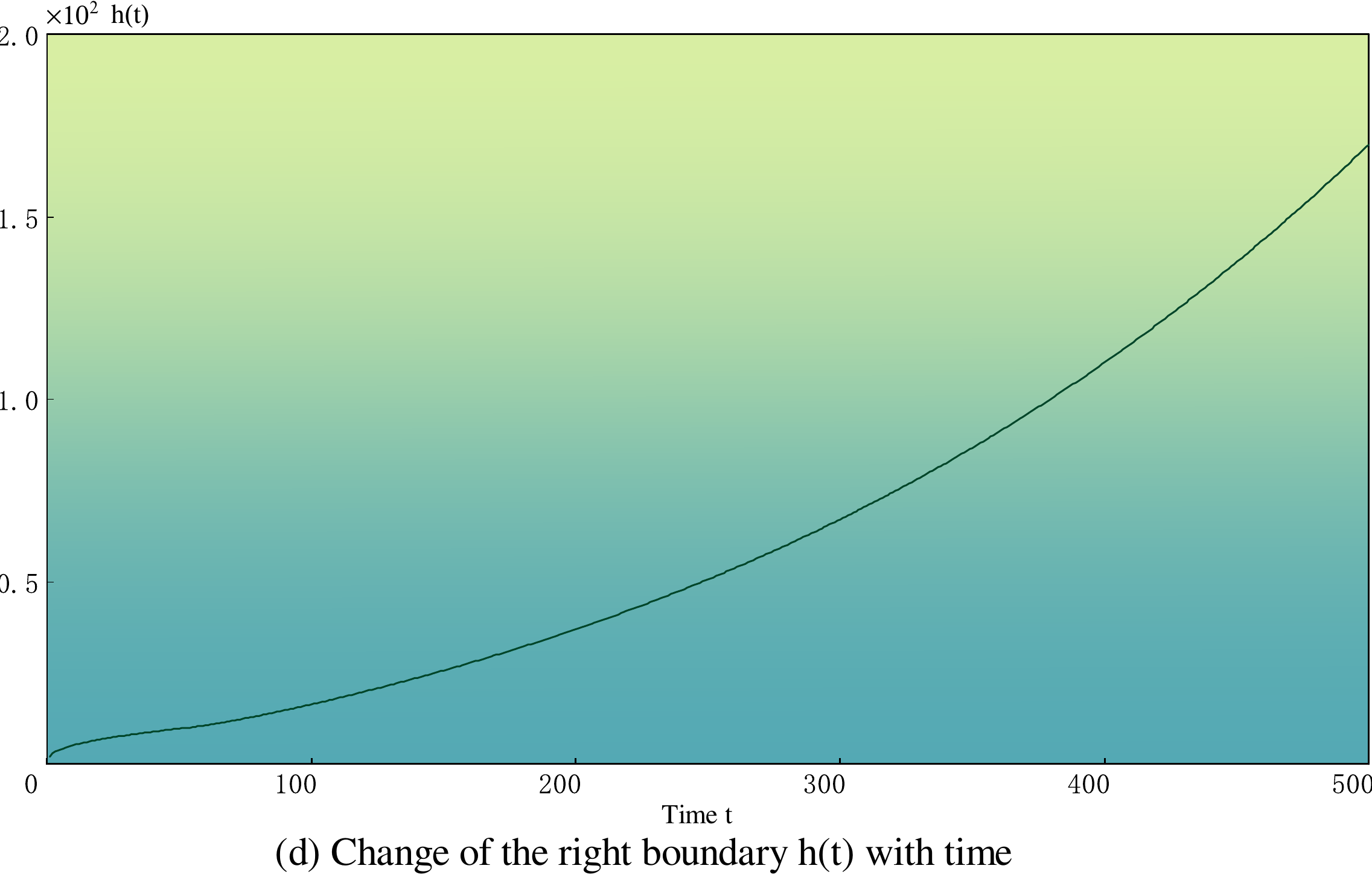}
} }
\caption{\footnotesize{When $G(u)=u$(without impulse), graphs (a)-(d) show that $u$ converges to a steady state and $h_{\infty}\geq 150$}.}
\label{A}
\end{figure}
\begin{figure}[!ht]
\centering
\subfigure{ {
\includegraphics[width=0.40\textwidth]{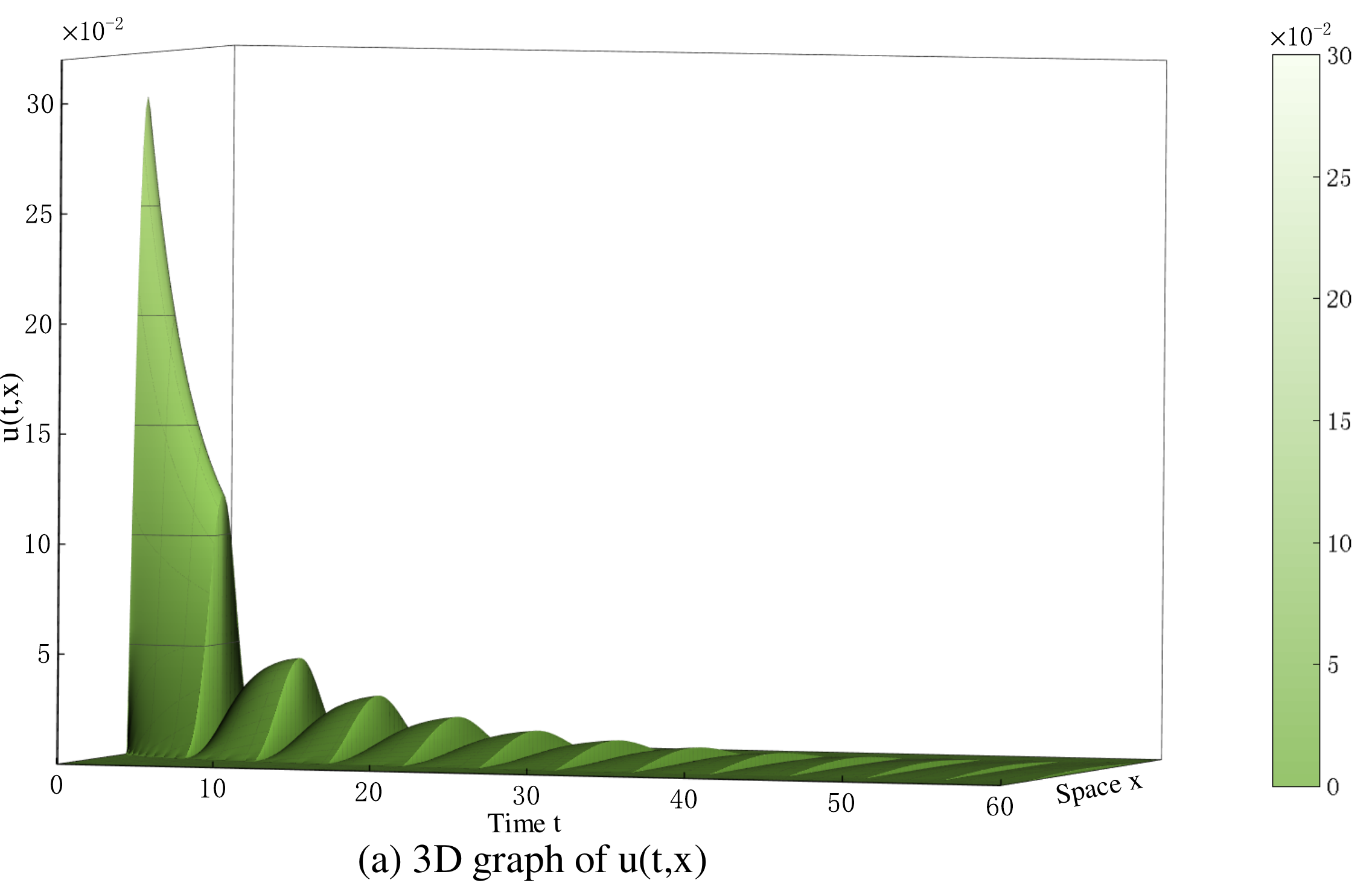}
} }
\subfigure{ {
\includegraphics[width=0.40\textwidth]{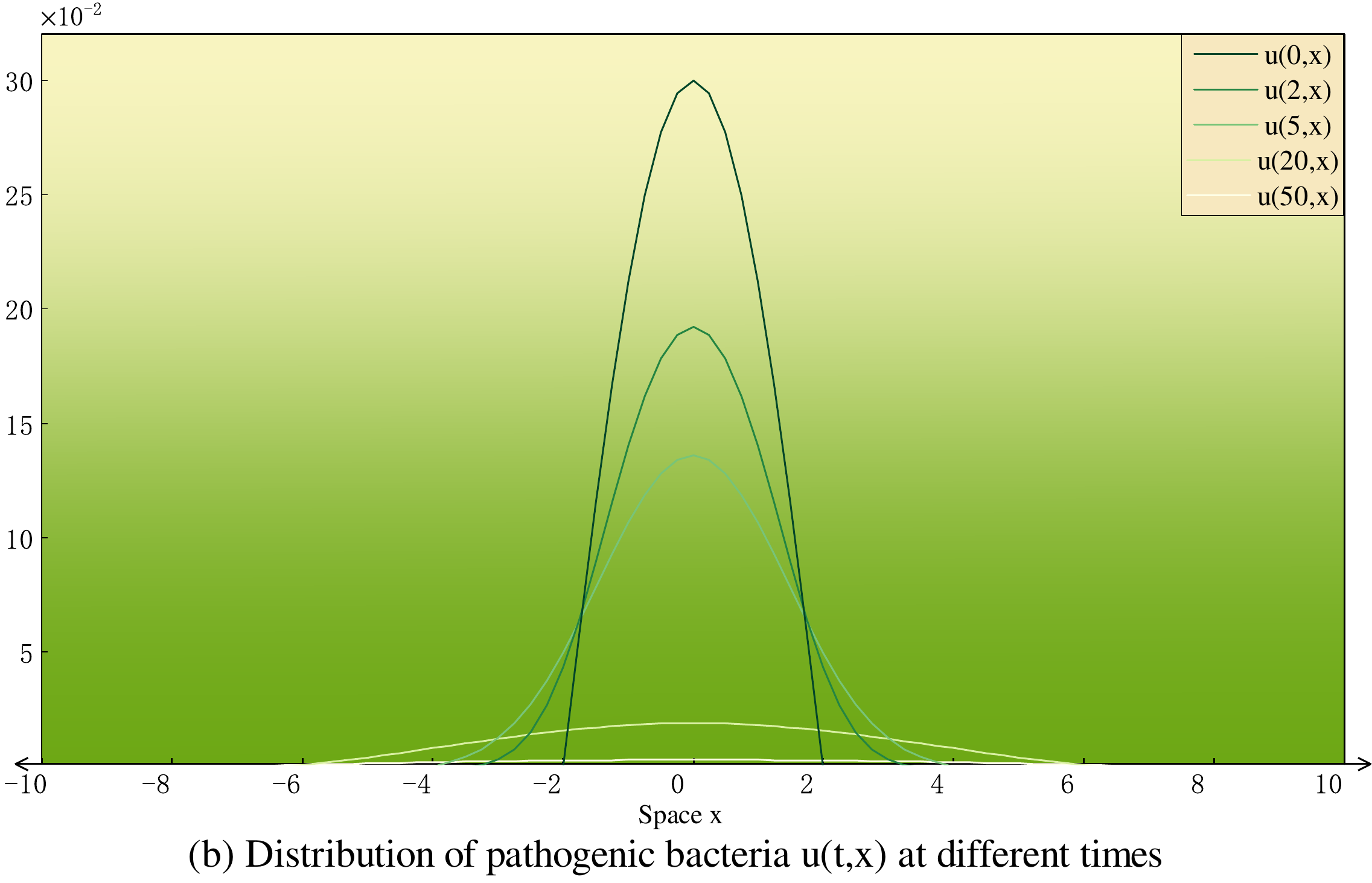}
} }
\subfigure{ {
\includegraphics[width=0.40\textwidth]{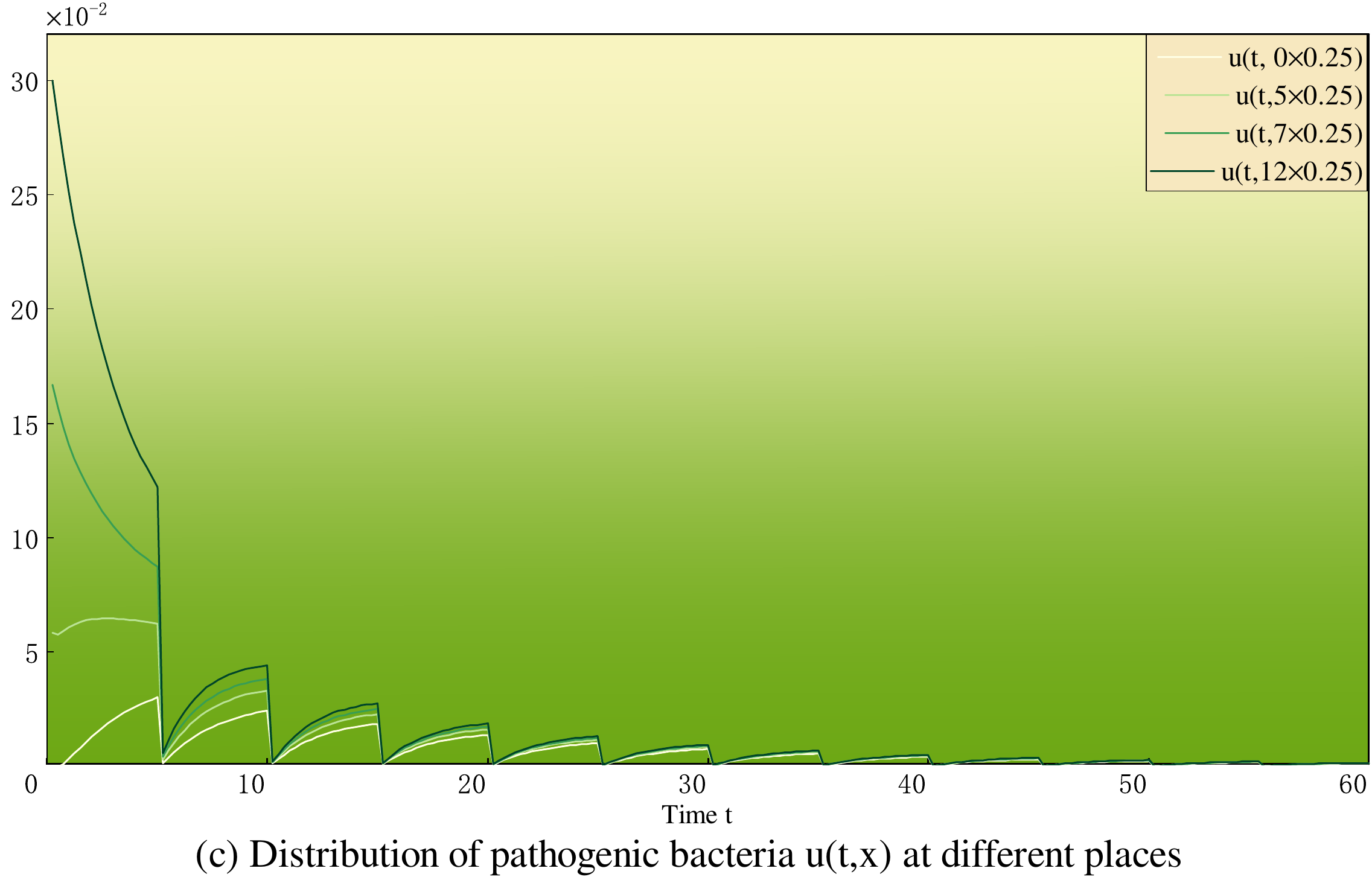}
} }
\subfigure{ {
\includegraphics[width=0.40\textwidth]{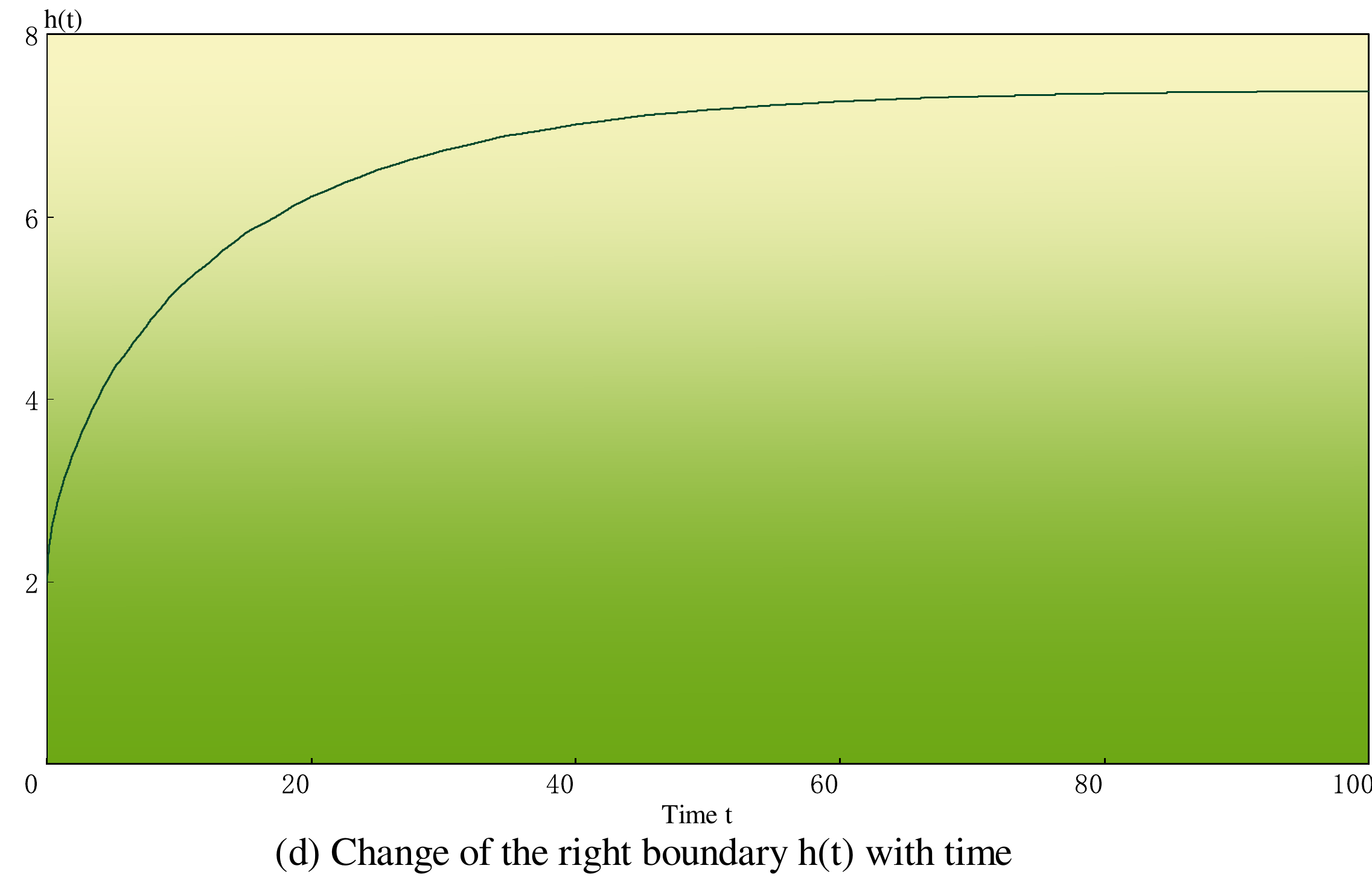}
} }
\caption{When $G(u)=\frac{0.5u}{10+u}$(with impulse), graphs (a)-(d) show $u$ decays to 0 and $h_{\infty}\leq 8$.}
\label{B}
\end{figure}
In the absence of impulse intervention in model \eqref{Zhou-Lin}, \autoref{A}\textcolor[rgb]{0.00,0.00,1.00}{(d)}, \autoref{lemma 3-1}\textcolor{blue}{(1,3)}, and \cite[Theorem 5.3(1)]{Zhou-Lin-Santos} give
\begin{equation*}
\lambda^{\vartriangle}_{1}(G'(0),(g_{\infty},h_{\infty}))\leq \lambda^{\vartriangle}_{1}(1,(-150,150))=-0.012<0.
\end{equation*}
Then, \autoref{lemma 3-6} yields that $h_{\infty}=-g_{\infty}=\infty$. This combined with \autoref{theorem 3-3} shows that the steady state of model \eqref{Zhou-Lin} is positive and unique, and it is also global asymptotically stable. In fact, it can be seen from \autoref{A}\textcolor[rgb]{0.00,0.00,1.00}{(a-c)} that the pathogenic bacteria converge to a positive distribution over time. This coincides with the conclusion of  \autoref{theorem 3-3}.

In the presence of impulsive intervention in model \eqref{Zhou-Lin}, it follows from \autoref{B}\textcolor[rgb]{0.00,0.00,1.00}{(d)} that $h_{\infty}-g_{\infty}<16$. This combined with \autoref{theorem 3-2} yields that $(0,0)$ is globally asymptotically stable. Actually, observing \autoref{B}\textcolor[rgb]{0.00,0.00,1.00}{(a-c)} gives us that the pathogenic bacteria eventually tends to zero. This result is the same as given by \autoref{theorem 3-2}.

It can be derived from comparing \textcolor[rgb]{0.00,0.00,1.00}{Figures}~\ref{A} and \ref{B} that the periodic pulse not only results in the long-lasting diseases becoming extinct, but also shrinks the region where the diseases are prevalent. Therefore, the pulse intervention is beneficial for the prevention and control of the diseases that are spread by faecal-oral pathway.
\subsection{The influence of expansion capacity}
To observe the impact of the expansion capability of infected individuals on the development of the diseases, this subsection performs numerical simulations without impulsive intervention.
\begin{exm}\label{exm2}
Fix $d_{1}=0.1$, $d_{2}=0.4$, $a_{11}=0.3$, $a_{12}=0.5$,  $a_{22}=0.1$, $m=1$, $a=10$, $G(u)=1$, and $\mu_{1}=10$. The expansion capability $\mu_{2}$ of infected individuals is taken to be $1$ and $10$, respectively.
\end{exm}
\begin{figure}[!ht]
\centering
\subfigure{ {
\includegraphics[width=0.40\textwidth]{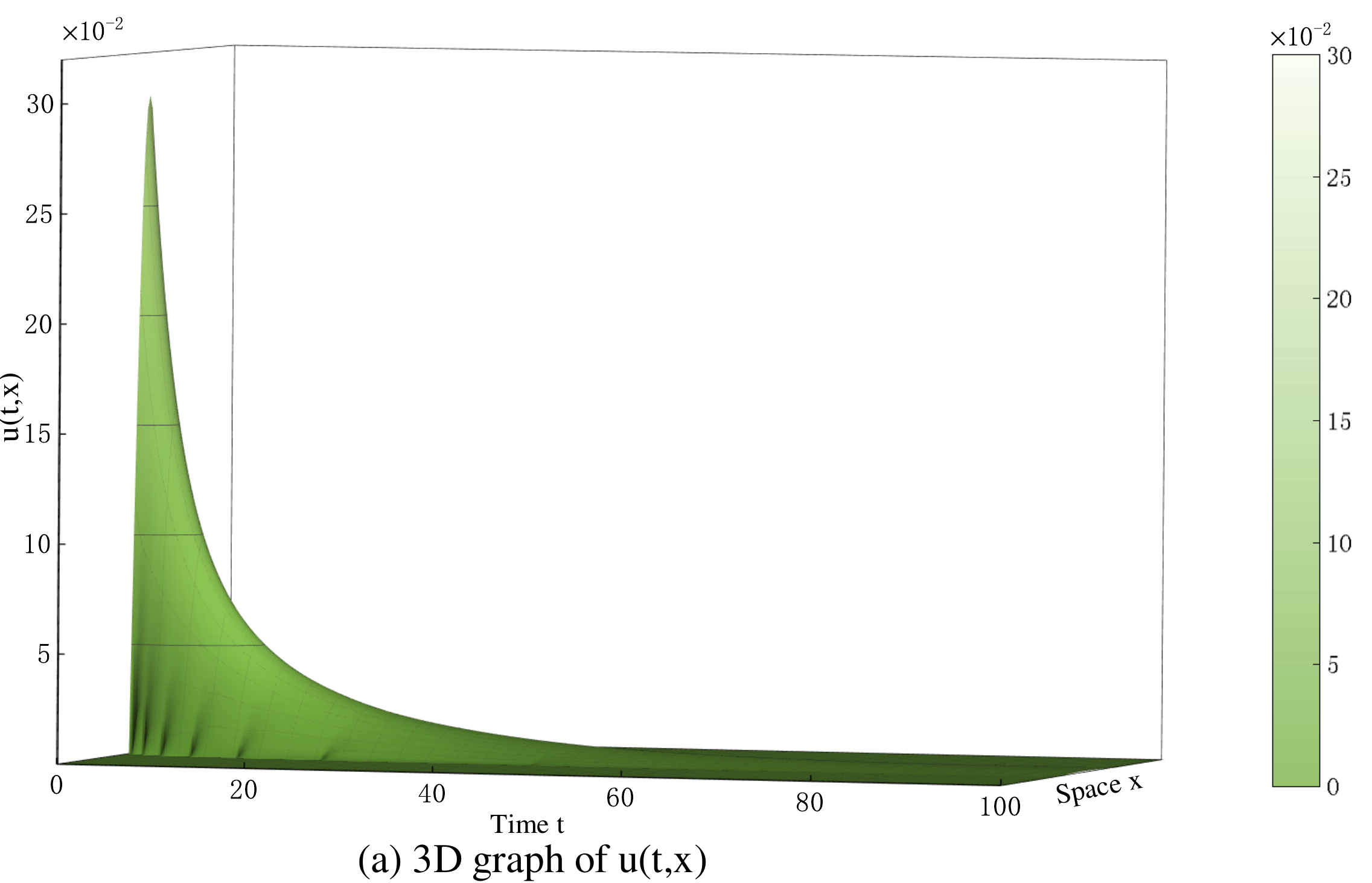}
} }
\subfigure{ {
\includegraphics[width=0.40\textwidth]{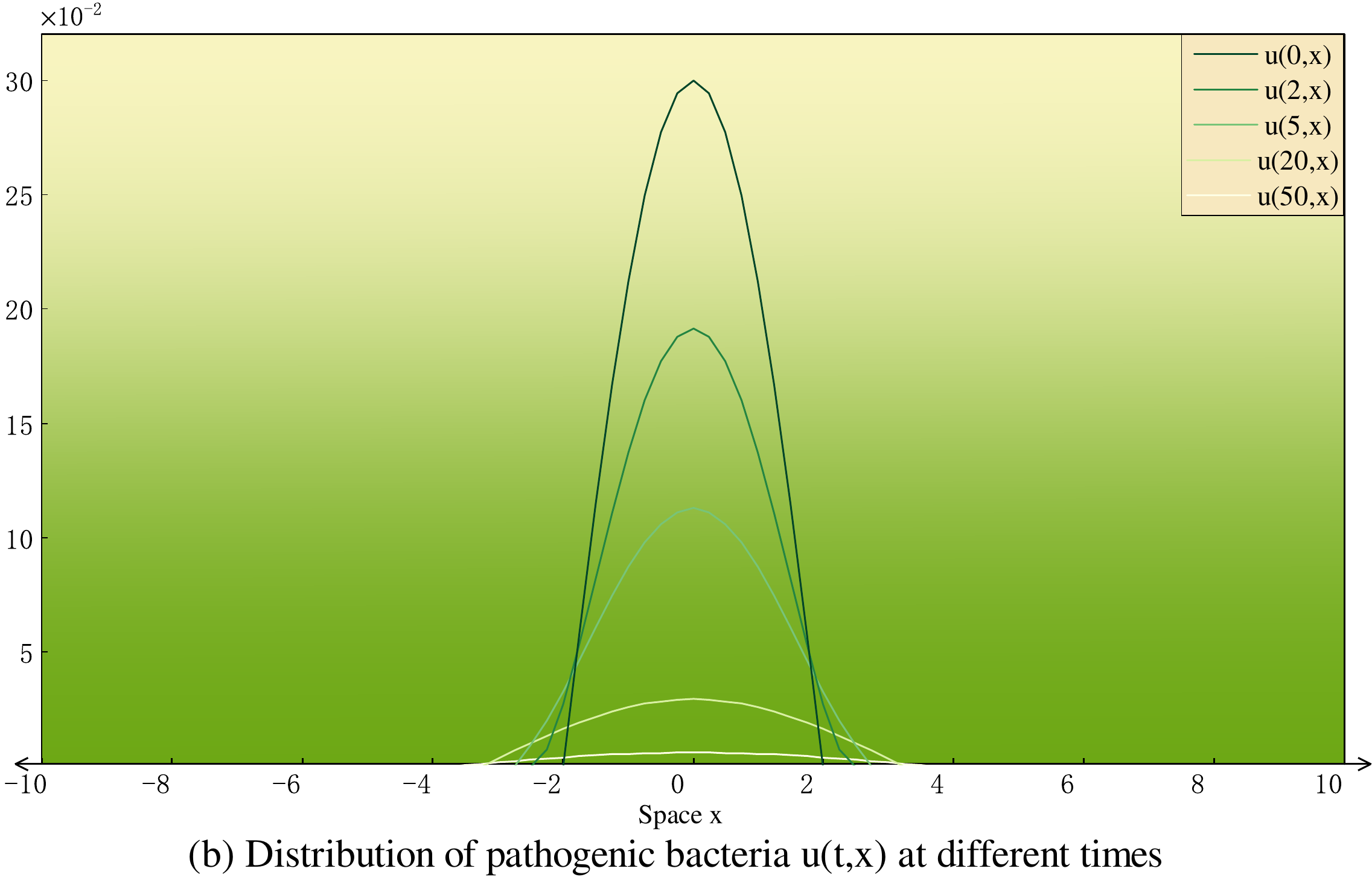}
} }
\subfigure{ {
\includegraphics[width=0.40\textwidth]{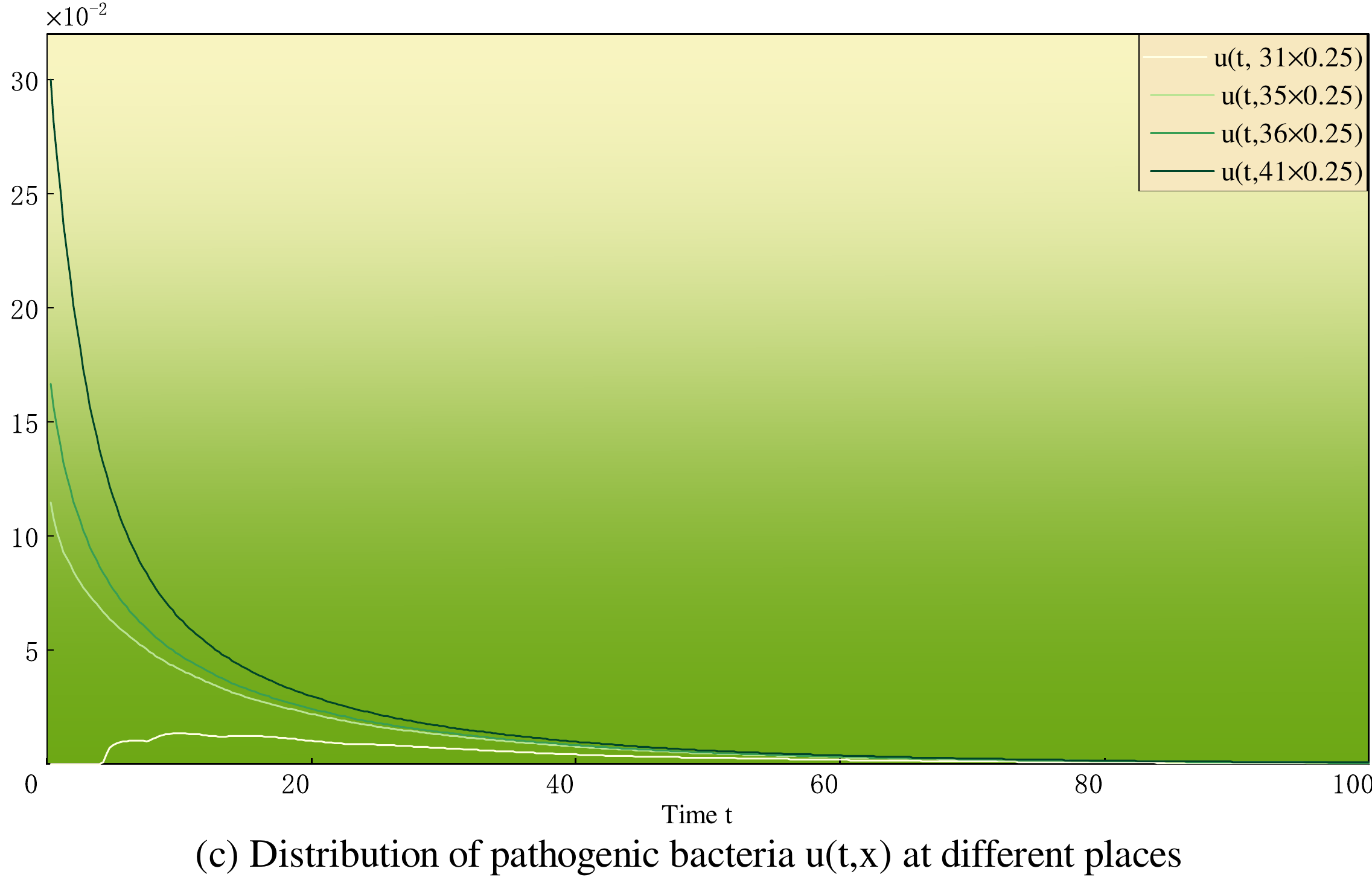}
} }
\subfigure{ {
\includegraphics[width=0.4\textwidth]{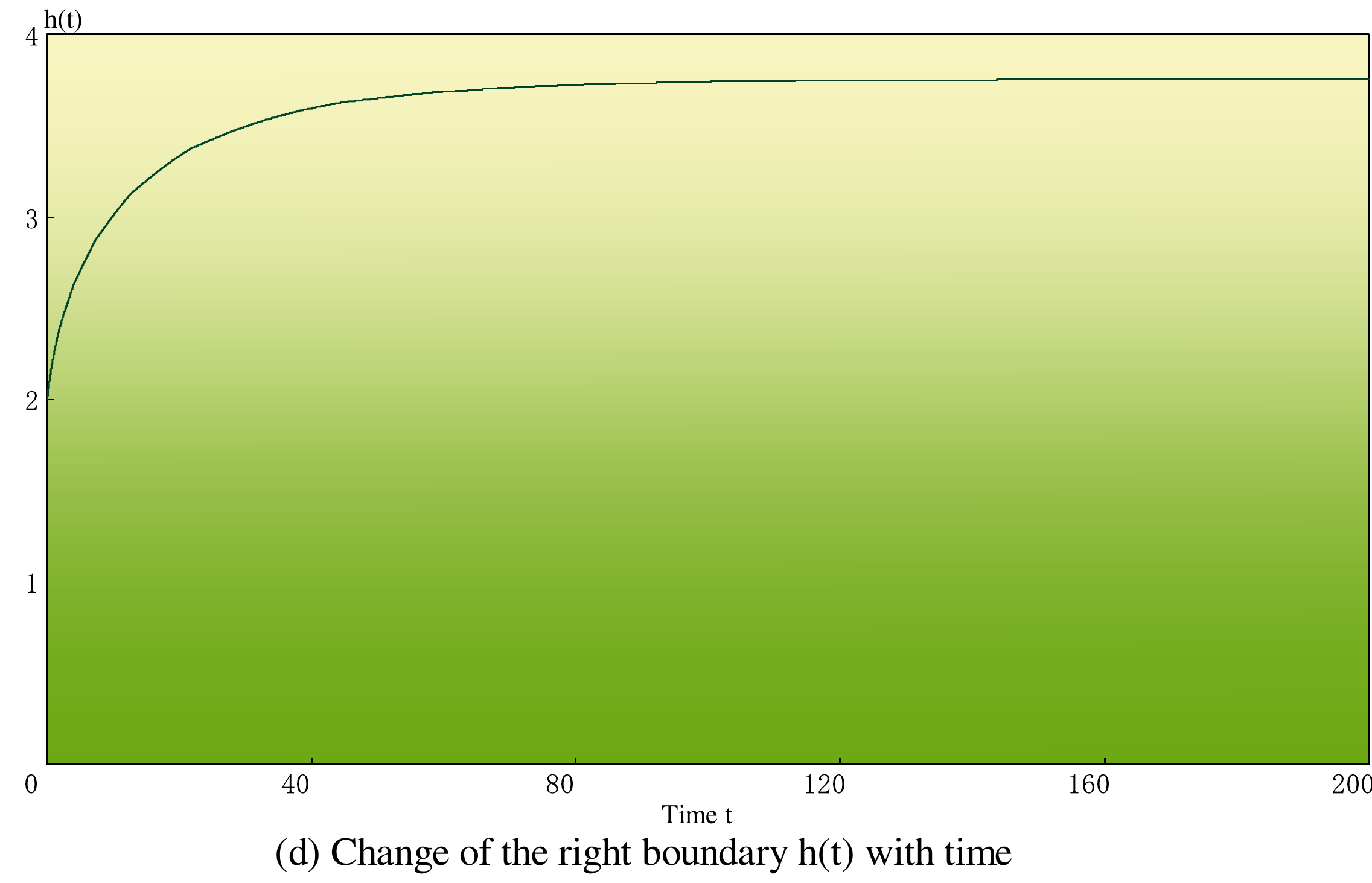}
} }
\caption{When $\mu_{2}=1$(weak expansion), graphs (a)-(d) exhibit $u$ decays to 0 and $h_{\infty}\leq 4$.}
\label{C}
\end{figure}
\begin{figure}[!ht]
\centering
\subfigure{ {
\includegraphics[width=0.40\textwidth]{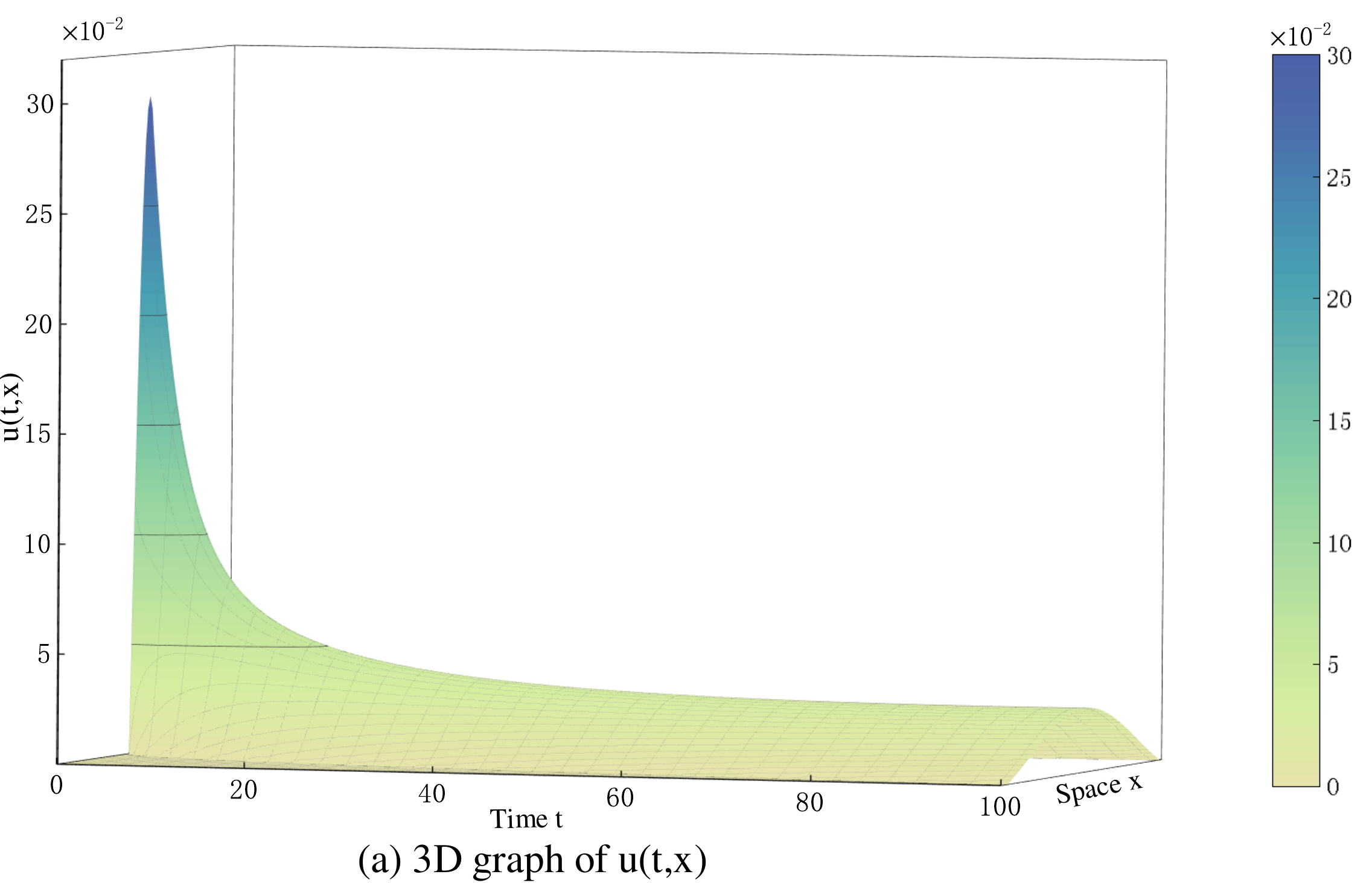}
} }
\subfigure{ {
\includegraphics[width=0.40\textwidth]{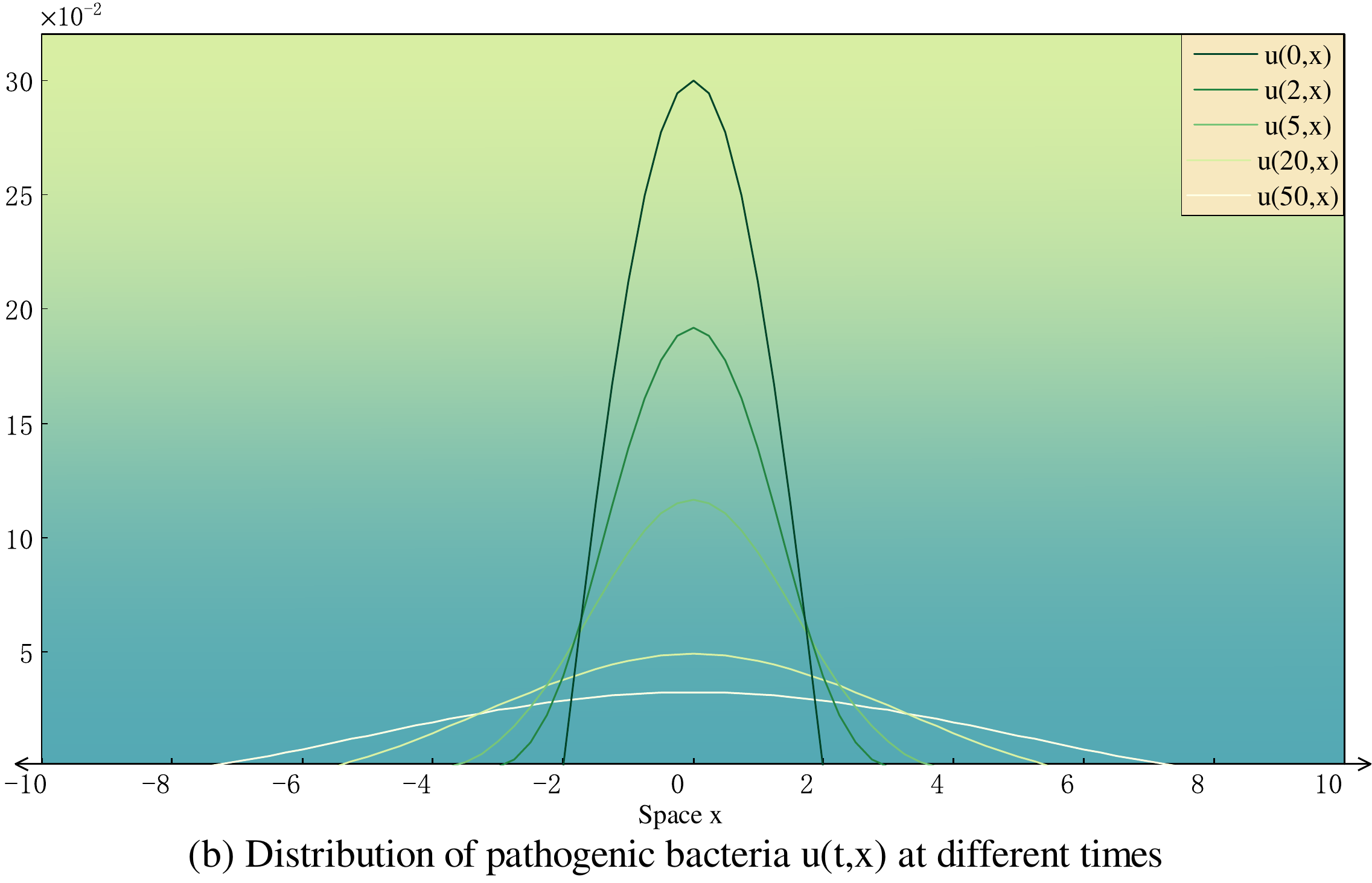}
} }
\subfigure{ {
\includegraphics[width=0.40\textwidth]{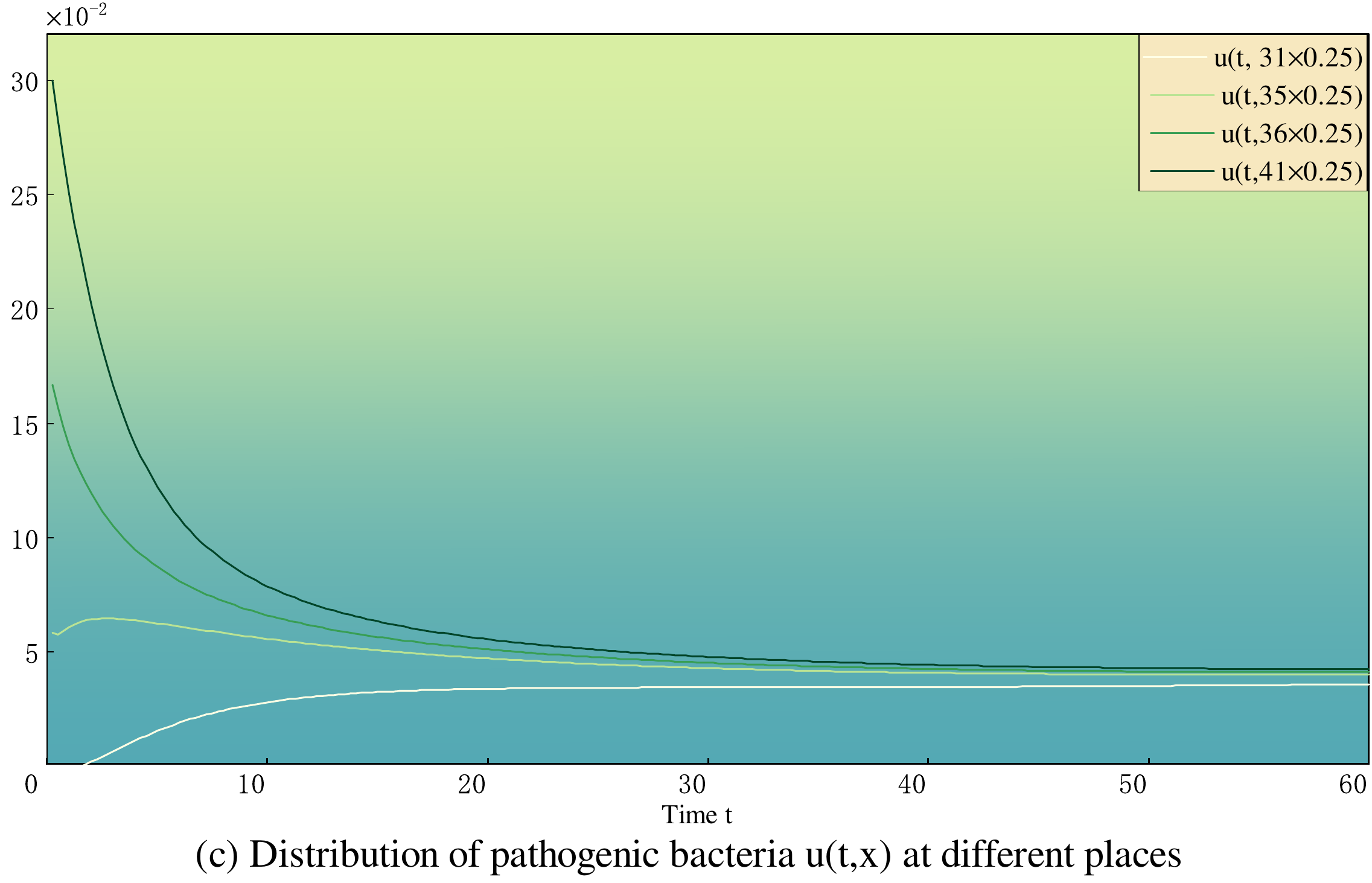}
} }
\subfigure{ {
\includegraphics[width=0.40\textwidth]{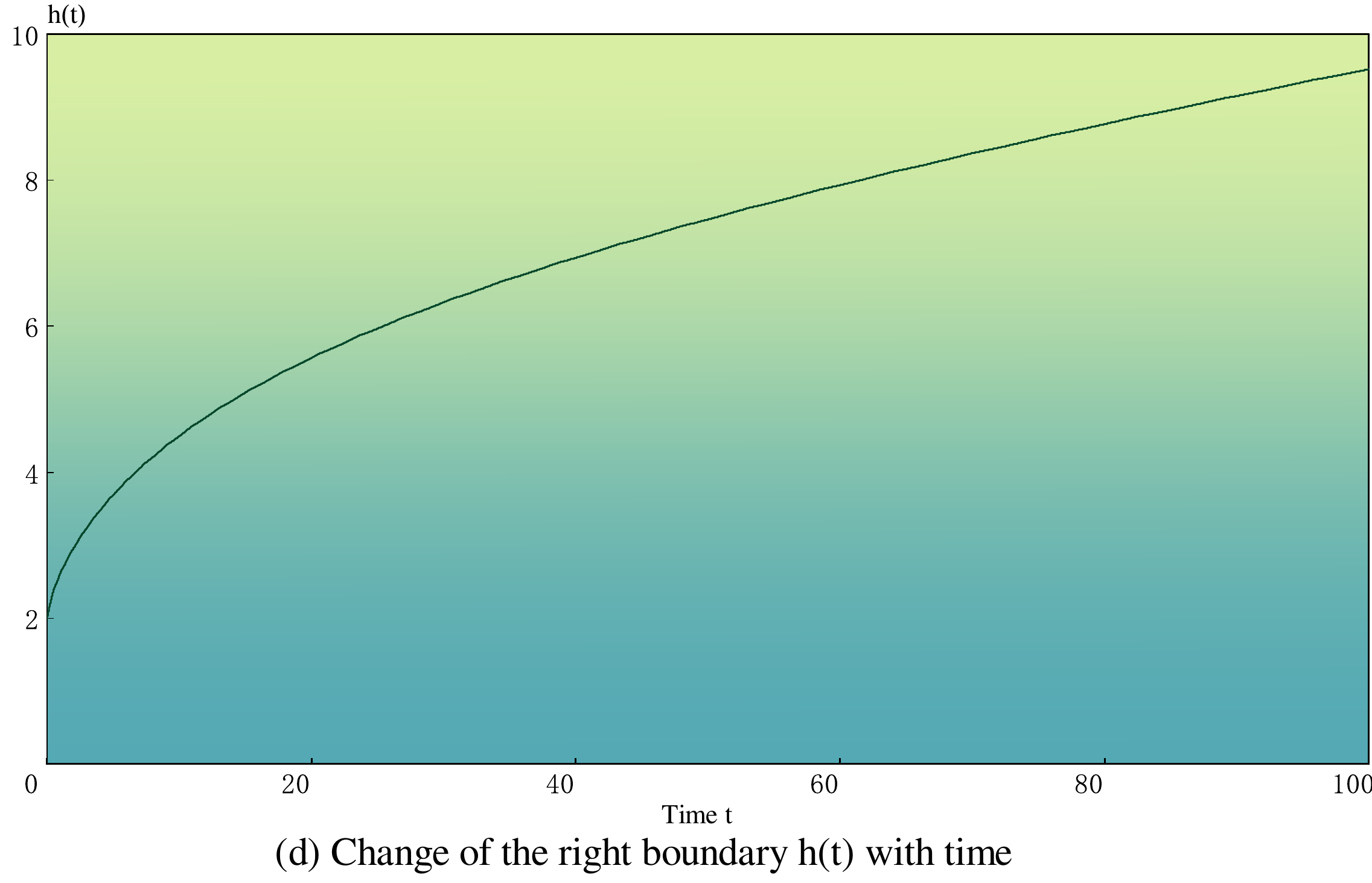}
} }
\caption{\small{When $\mu_{2}=10$(strong expansion), graphs (a)-(d) show $u$ converges to a steady state and $h_{\infty}\geq9$}.}
\label{D}
\end{figure}
As the expansion capacity of the infected individuals is weak, \autoref{C}\textcolor[rgb]{0.00,0.00,1.00}{(d)} implies that $h_{\infty}-g_{\infty}<\infty$. With the help of \autoref{theorem 3-2}, it can then be obtained that both $u$ and $v$ uniformly converge to zero when time tends to infinity regardless of the choice of initial functions. In fact, it can be seen in \autoref{C}\textcolor[rgb]{0.00,0.00,1.00}{(a-c)} that the number of infected individuals at each location tends to zero over time. This agrees with the conclusion offered by \autoref{theorem 3-2}.

As the expansion capacity of infected individuals is strong, it follows from \autoref{D}\textcolor[rgb]{0.00,0.00,1.00}{(d)} that $h_{\infty}-g_{\infty}\geq 18$. With the help of \autoref{lemma 3-1}\textcolor{blue}{(1,3)} and \cite[Theorem 5.3(1)]{Zhou-Lin-Santos}, it can then be obtained that
\begin{equation*}
\lambda^{\vartriangle}_{1}(G'(0),(g_{\infty},h_{\infty}))\leq \lambda^{\vartriangle}_{1}(1,(-9,9))=-0.002<0.
\end{equation*}
This combined with \autoref{lemma 3-6} yields that $h_{\infty}-g_{\infty}=+\infty$. Next, \autoref{theorem 3-3} implies that model \eqref{Zhou-Lin} has only one  steady state solution, which is globally asymptotically stable. In fact, it can be obtained from observing \autoref{D}\textcolor[rgb]{0.00,0.00,1.00}{(a-c)} that the pathogenic bacteria gradually converge to a positive distribution as time progresses, which is consistent with the finding of \autoref{theorem 3-3}.

The comparison of \textcolor[rgb]{0.00,0.00,1.00}{Figures}~\ref{C} and \ref{D} reveals that the weak expansion capacity of the infected individual reduces the speed of the expansion front of the diseases, and enables the diseases  existing in the entire region to be contained  in a limited region, and extincts the persistent diseases. This comparison suggests that the diseases can be prevented and controlled by limiting the expansion capacity of the infected individual.
\section{\bf Conclusion and discussion}\label{Section-5}
Systematic environmental disinfection is an effective measure for the prevention and control of  infectious diseases transmitted by faecal-oral route. This measure usually leads to a dramatic reduction in the number of pathogenic bacteria in a short period of time. In order to capture this phenomenon, this paper develops an impulsive faecal-oral model with free boundary by incorporating the pulse intervention into the model formulated in \cite{Wang-Du}. Due to the introduction of the intervention, the theoretical analysis is completely different from that of \cite{Wang-Du}. More precisely, the solution no longer has the half-flow property, thus some energy methods are no longer available; the impulsive condition should also be considered when constructing the upper and lower solutions; the steady state is governed by a more complex parabolic problem.

The main content of this paper are summarised as follows. First of all, \autoref{theorem 2-1} checks that the model has a unique globally nonnegative classical solution. The spreading-vanishing dichotomy for the model is then proved in \autoref{theorem 3-4} by using the comparison principle and the upper and lower solutions method. To study the long-time dynamical behaviours, we define two principal eigenvalues $\lambda^{\vartriangle}_{1}(h_{0})$ and $\lambda^{\vartriangle}_{1}(\infty)$. Next, the spreading-vanishing criteria is given with the help of these two values. Specifically, the following conclusions are valid:
\begin{itemize}
\item{Assume $\lambda^{\vartriangle}_{1}(\infty)>0$. Then the diseases are vanishing (see \autoref{theorem 3-1});}
\item{Assume $\lambda^{\vartriangle}_{1}(\infty)=0$. Then the diseases are extinct (see \autoref{lemma 3-4});}
\item{Assume $\lambda^{\vartriangle}_{1}(\infty)\leq 0$ and $\lambda^{\vartriangle}_{1}(h_{0})\leq 0$. Then the diseases are spreading (see \autoref{theorem 3-6});}
\item{Assume $\lambda^{\vartriangle}_{1}(\infty)\leq 0$ and $\lambda^{\vartriangle}_{1}(h_{0})>0$. Then,
   \begin{itemize}
     \item[$\bullet$]{for any given positive constant $\mu_{1}$, there exists a positive constant $\mu_{0}$ such that spreading happens as $\mu_{2}\in( \mu_{0},+\infty)$, and vanishing happens as $\mu_{2}\in(0, \mu_{0})$ (see \autoref{theorem 3-8}), }
      \item[$\bullet$]{or for any given initial function $v_{0}(x)$, there exists a positive constant $\kappa_{0}$ such that spreading happens as $\kappa\in(\kappa_{0}, +\infty)$, and vanishing happens as $\kappa\in(0,\kappa_{0})$ provided the impulsive function is linear and the initial function $(u_{0}(x),v_{0}(x))=(\kappa\upsilon(x),v_{0}(x))$ for some $\kappa>0$ (see \autoref{theorem 3-10}).}
    \end{itemize}}
\end{itemize}
At the end, the numerical simulation verifies the correctness of the theoretical results and demonstrates more intuitively the effect of the pulse intervention and the expansion capacity on the evolution of the diseases.

Our theoretical results and numerical simulations simultaneously show that both increasing the impulse intensity and decreasing the expansion capacity of infected individuals are not favourable to the evolution of the diseases. This coincides with the findings observed in real life. Therefore, we suggest to the preventive and control authority of infectious diseases that the prevalence of the diseases transmitted by the faecal-oral route can be counteracted by restricting the expansion of the infected individuals and increasing the intensity of the regular disinfection.

There are still two shortcomings in this paper. One is that when $\lambda^{\vartriangle}_{1}(\infty)=0$, we only obtain that the diseases are extinct and not  whether the infected area is finite or infinite. The other is that \autoref{theorem 3-10} requires the additional assumption that the impulsive function is linear. These two disadvantages will be the direction of our future research.

\section*{\bf Declaration of competing interest}
 The authors have no conflict of interest.

\end{document}